\documentclass[reqno]{amsart}
\usepackage{amssymb}
\usepackage{euscript}
\usepackage{graphicx}
\usepackage[left]{lineno}
\usepackage{placeins}
\usepackage{float}
\usepackage{xcolor,bbm}
\usepackage{enumitem}
\usepackage{todonotes}
\usepackage{mathtools}
\usepackage{relsize}
\usepackage{tikz}
\usetikzlibrary{arrows.meta,fadings}

\newtheorem{theorem}{Theorem}[section]
\newtheorem{lemma}[theorem]{Lemma}

\newtheorem{proposition}[theorem]{Proposition}
\newtheorem{corollary}[theorem]{Corollary}
\newtheorem{definition}[theorem]{Definition}

\newtheorem{example}[theorem]{Example}
\theoremstyle{remark}
\newtheorem{remark}[theorem]{Remark}
\numberwithin{equation}{section}

\usepackage{algpseudocode}

\floatstyle{ruled}
\newfloat{algorithm}{tbp}{loa}
\floatname{algorithm}{Algorithm}

\usepackage[margin=2.34cm]{geometry}
\newcommand{\R}{\mathbb{R}}

\DeclareMathOperator*{\argmin}{arg\,min}

\DeclareMathOperator*{\id}{id}

\title{Decay rates for monotone operator flows in the space of probability measures}

\author{Radu Ioan Bo{\c t}}
\address{Faculty of Mathematics, University of Vienna,
Oskar-Morgenstern-Platz 1, 1090 Vienna, Austria}
\email{radu.bot@univie.ac.at}

\author{Alberto Dom\'inguez Corella}
\address{Departamento de Matem\'aticas, Universidad de Sonora,
Blvd. Luis Encinas y Rosales s/n, Col. Centro,
83000 Hermosillo, Sonora, M\'exico}
\email{alberto.corella@unison.mx, alberto.of.sonora@gmail.com}

\author{Tr\'{\i} Minh L\^e}
\address{Faculty of Mathematics, University of Vienna,
Oskar-Morgenstern-Platz 1, 1090 Vienna, Austria}
\email{tri.minh.le@univie.ac.at}

    \keywords{Wasserstein flows, monotone operators, long-time behavior, decay rates, Tikhonov regularization}

\subjclass[2020]{Primary 47H05; Secondary 47J35, 49Q22, 47J26}

\thanks{The work of the third author was supported by  Austrian Science Fund (FWF) [10.55776/STA223].}

\begin{document}

\begin{abstract}
Motivated by the mean-field dynamics of large populations of interacting agents, we study the long-time behavior of flows driven by monotone demicontinuous operators on the Wasserstein space. For cocoercive operators, we prove that the flow converges to a zero of the operator and that the operator residual and the metric speed decay as $o(t^{-1/2})$, so that the population comes to rest (the agents asymptotically stop moving). We also consider Tikhonov regularized flows of merely monotone demicontinuous operators and show that the operator residual and the metric speed tend to zero, with rate $O(1/t)$ when the Tikhonov parameter is $\alpha/t$ with $\alpha\geq 1$, and that the trajectory converges to the zero of minimal second moment, which is a Dirac mass, so that the agents reach consensus (all agents reach the same point). Finally, we study a discrete analogue for nonexpansive operators and show that the generated sequence converges to the fixed point closest to a reference Dirac mass; this fixed point is again a Dirac mass, so that consensus is also reached under this scheme. The main  arguments rely on a specially tailored superposition principle, which represents the flow as a curve in a fixed Hilbert space of random trajectories,  where the crucial estimates are derived.
\end{abstract}

\maketitle
\setcounter{tocdepth}{1}
\tableofcontents

\section{Introduction, motivation, and contributions}
   We use the introduction to motivate the paper and to comment informally on its contributions. The main results are stated rigorously in the next section, and their proofs are given in Section~\ref{pruebas}.
\subsection{On zeroes of monotone operators and mean-field dynamics}

\subsubsection*{Monotone operators and their zeroes}
Monotone operators are one of the main unifying ideas of nonlinear analysis and optimization. An operator $A:\mathbb R^d\to\mathbb R^d$ is monotone if
\[
  \langle A(x)-A(y),x-y\rangle\geq 0\quad\text{for all } x,y\in \mathbb R^d .
\]
Though the main examples of monotone operators are gradients of convex functions, there are many other examples of interest. For instance, if $d=d_1+d_2$ and $L:\mathbb R^{d_1}\times\mathbb R^{d_2}\to\mathbb R$ is a differentiable convex--concave function, then the so-called \textit{saddle} operator $A(x_1,x_2):=\bigl(\nabla_{x_1}L(x_1,x_2),-\nabla_{x_2}L(x_1,x_2)\bigr)$ is monotone. Operators arising from variational inequalities and fixed-point problems are often also monotone.
\smallbreak\noindent
Finding zeroes of a monotone operator, i.e., points $x\in\mathbb R^d$ satisfying $A(x)=0$, is a core problem that arises in minimization, min--max problems and equilibria alike. In the case of gradients of convex functions, finding their zeroes amounts to finding the minimizers of the corresponding functions. For the saddle operator above, the zeroes are the saddle points of $L$, i.e., the equilibria of the zero-sum game in which one player chooses $x_1$ to minimize $L$ and the other chooses $x_2$ to maximize it.
\smallbreak\noindent
Many classical methods for finding zeroes are discretizations of the associated evolution equation
\begin{align}\label{eveq}
    \dot x(t)+A\bigl(x(t)\bigr)=0.
\end{align}
Its explicit Euler discretization with step size $\tau>0$ is the forward iteration $x_{k+1}:=x_k-\tau A(x_k)$, which is gradient descent when $A$ is a gradient and simultaneous gradient descent--ascent when $A$ is a saddle operator, whereas its implicit Euler discretization $x_{k+1}:=x_k-\tau A(x_{k+1})$ is the proximal point algorithm. The analysis of the equation (\ref{eveq}) is usually meaningful because it is the common limit of these schemes as $\tau\longrightarrow0^+$ and helps to gain insight into them. In continuous time there are no step sizes, the effects of the operator are not mixed with those of the discretization, and the analysis typically reduces to showing that a suitable Lyapunov function, such as the distance to a zero or the operator residual $| A(x(t)) |$, decreases along trajectories, ideally at an explicit rate. Arguments of this kind very often survive discretization, and the flow serves as a guide for the design and the analysis of algorithms.
\subsubsection*{Mean-field dynamics and operators on the Wasserstein space}
In some applications, there is interest in considering many particles or agents at the same time. Since this can be computationally expensive or intractable, one considers their population as a whole, through a probability measure describing the global behavior. We may think of $N$ agents (particles of a gas, individuals in a crowd, neurons of a wide neural network, players of a large game), each described by a state $x_i\in\mathbb R^d$ and each adjusting its state according to its own position and the positions of all the others. Consider the dynamics
\begin{equation}\label{particles}
  \dot x_i(t)=-A_{\mu^N_t}\bigl(x_i(t)\bigr),\qquad i \in \{ 1,\cdots,N\},
\end{equation}
in which each agent follows the negative of the operator, but the operator is evaluated against the whole population, encoded by its empirical measure $\mu^N_t:=\frac1N\sum_{j=1}^N\delta_{x_j(t)}$. Each agent thus perceives the crowd only through $\mu^N_t$, to which it contributes a mere $1/N$; this is the so-called \emph{mean-field} regime.
%%%%%%%%%%%%%%%%%%%%
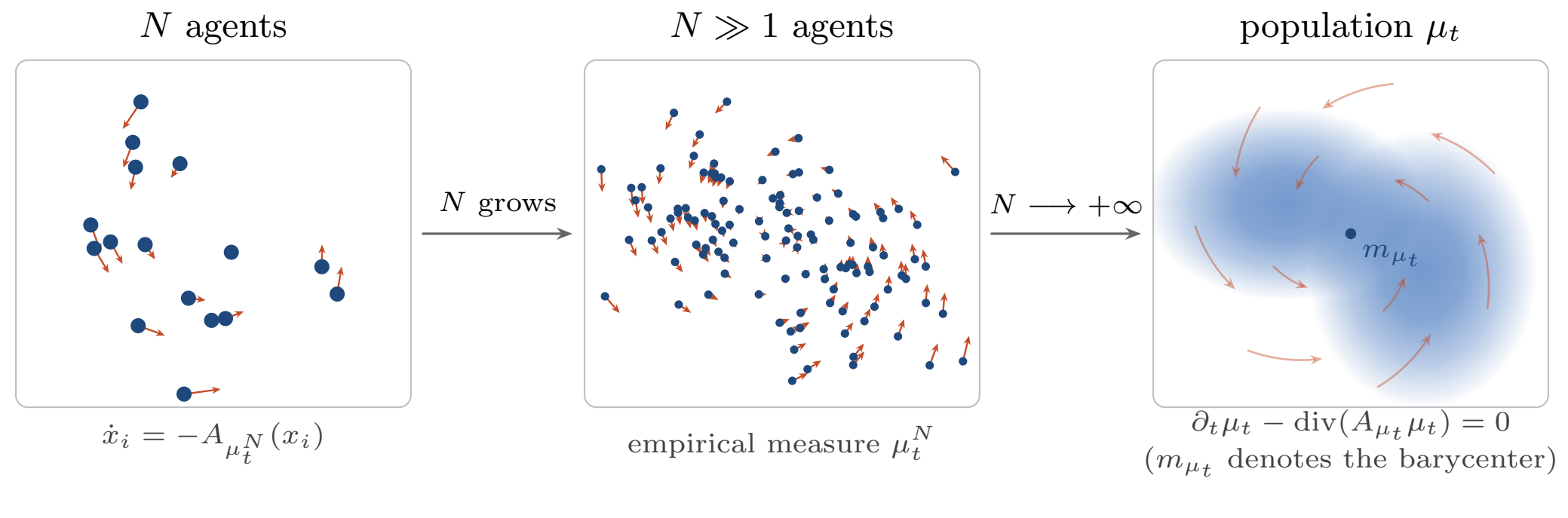
\begin{figure}[H]
\centering
\definecolor{agentcol}{RGB}{31,73,125}
\definecolor{velcol}{RGB}{196,78,42}
\definecolor{denscol}{RGB}{70,120,190}
\pgfdeclareradialshading{mfbumpshade}{\pgfpointorigin}{%
  color(0bp)=(pgftransparent!0); color(9bp)=(pgftransparent!12);
  color(17bp)=(pgftransparent!50); color(25bp)=(pgftransparent!100);
  color(50bp)=(pgftransparent!100)}%
\pgfdeclarefading{mfbump}{\pgfuseshading{mfbumpshade}}
\begin{tikzpicture}[
  scale=1.05, line cap=round,
  frame/.style={draw=black!25, rounded corners=3pt, line width=0.4pt},
  paneltitle/.style={font=\small},
  panelsub/.style={font=\scriptsize, text=black!75},
  agent/.style={fill=agentcol},
  vel/.style={-{Stealth[length=2.2pt, width=1.8pt]}, draw=velcol, line width=0.45pt},
  stream/.style={draw=velcol, line width=0.5pt, opacity=0.45, -{Stealth[length=3pt, width=2.4pt]}},
  dens/.style={fill=denscol},
  link/.style={-{Stealth[length=4pt]}, line width=0.6pt, draw=black!60}
]
  \draw[frame] (-1.650,-1.45) rectangle (1.650,1.45);
  \node[paneltitle] at (0.000,1.72) {$N$ agents};
  \node[panelsub] at (0.000,-1.75) {$\dot x_i=-A_{\mu^N_t}(x_i)$};
  \draw[vel] (-0.650,0.556) -- ++(-0.040,-0.191);
  \draw[vel] (-0.859,-0.069) -- ++(0.101,-0.183);
  \draw[vel] (-0.996,-0.123) -- ++(0.123,-0.206);
  \draw[vel] (-0.605,1.101) -- ++(-0.152,-0.226);
  \draw[vel] (-1.024,0.073) -- ++(0.086,-0.227);
  \draw[vel] (-0.279,0.585) -- ++(-0.075,-0.120);
  \draw[vel] (-0.570,-0.091) -- ++(0.083,-0.124);
  \draw[vel] (-0.673,0.763) -- ++(-0.079,-0.212);
  \draw[vel] (0.101,-0.709) -- ++(0.152,0.059);
  \draw[vel] (-0.245,-1.340) -- ++(0.306,0.041);
  \draw[vel] (-0.208,-0.539) -- ++(0.143,-0.016);
  \draw[vel] (0.150,-0.155) -- ++(0.038,0.025);
  \draw[vel] (1.034,-0.504) -- ++(0.037,0.229);
  \draw[vel] (-0.628,-0.769) -- ++(0.223,-0.081);
  \draw[vel] (0.906,-0.276) -- ++(0.002,0.185);
  \draw[vel] (-0.015,-0.724) -- ++(0.165,0.038);
  \fill[agent] (-0.650,0.556) circle (1.8pt);
  \fill[agent] (-0.859,-0.069) circle (1.8pt);
  \fill[agent] (-0.996,-0.123) circle (1.8pt);
  \fill[agent] (-0.605,1.101) circle (1.8pt);
  \fill[agent] (-1.024,0.073) circle (1.8pt);
  \fill[agent] (-0.279,0.585) circle (1.8pt);
  \fill[agent] (-0.570,-0.091) circle (1.8pt);
  \fill[agent] (-0.673,0.763) circle (1.8pt);
  \fill[agent] (0.101,-0.709) circle (1.8pt);
  \fill[agent] (-0.245,-1.340) circle (1.8pt);
  \fill[agent] (-0.208,-0.539) circle (1.8pt);
  \fill[agent] (0.150,-0.155) circle (1.8pt);
  \fill[agent] (1.034,-0.504) circle (1.8pt);
  \fill[agent] (-0.628,-0.769) circle (1.8pt);
  \fill[agent] (0.906,-0.276) circle (1.8pt);
  \fill[agent] (-0.015,-0.724) circle (1.8pt);
  \draw[frame] (3.100,-1.45) rectangle (6.400,1.45);
  \node[paneltitle] at (4.750,1.72) {$N\gg1$ agents};
  \node[panelsub] at (4.750,-1.75) {empirical measure $\mu^N_t$};
  \draw[vel] (3.662,-0.057) -- ++(0.054,-0.117);
  \draw[vel] (4.783,-0.057) -- ++(0.005,0.006);
  \draw[vel] (4.182,0.586) -- ++(-0.039,-0.088);
  \draw[vel] (3.879,0.208) -- ++(0.015,-0.105);
  \draw[vel] (4.261,0.274) -- ++(-0.009,-0.066);
  \draw[vel] (3.526,0.279) -- ++(0.023,-0.147);
  \draw[vel] (4.947,-0.338) -- ++(0.028,0.037);
  \draw[vel] (4.673,0.295) -- ++(-0.029,-0.022);
  \draw[vel] (3.847,1.010) -- ++(-0.071,-0.144);
  \draw[vel] (4.619,-0.206) -- ++(0.028,-0.005);
  \draw[vel] (4.241,0.469) -- ++(-0.029,-0.077);
  \draw[vel] (4.096,0.510) -- ++(-0.027,-0.094);
  \draw[vel] (4.163,0.504) -- ++(-0.030,-0.087);
  \draw[vel] (4.991,-0.007) -- ++(-0.010,0.027);
  \draw[vel] (4.312,0.074) -- ++(0.011,-0.051);
  \draw[vel] (4.270,-0.201) -- ++(0.043,-0.044);
  \draw[vel] (3.881,0.175) -- ++(0.019,-0.103);
  \draw[vel] (4.694,0.685) -- ++(-0.073,-0.036);
  \draw[vel] (3.472,-0.052) -- ++(0.062,-0.138);
  \draw[vel] (4.556,-0.507) -- ++(0.064,0.001);
  \draw[vel] (3.945,0.213) -- ++(0.012,-0.098);
  \draw[vel] (4.891,0.512) -- ++(-0.063,-0.007);
  \draw[vel] (4.020,0.110) -- ++(0.020,-0.085);
  \draw[vel] (4.062,0.829) -- ++(-0.061,-0.112);
  \draw[vel] (3.965,0.135) -- ++(0.020,-0.092);
  \draw[vel] (4.394,0.204) -- ++(-0.007,-0.048);
  \draw[vel] (4.091,-0.173) -- ++(0.048,-0.065);
  \draw[vel] (4.194,0.081) -- ++(0.016,-0.065);
  \draw[vel] (4.841,0.498) -- ++(-0.059,-0.012);
  \draw[vel] (4.187,0.504) -- ++(-0.031,-0.084);
  \draw[vel] (4.013,0.650) -- ++(-0.039,-0.110);
  \draw[vel] (4.197,0.472) -- ++(-0.028,-0.082);
  \draw[vel] (3.490,0.382) -- ++(0.013,-0.155);
  \draw[vel] (3.271,-0.523) -- ++(0.123,-0.140);
  \draw[vel] (4.033,-0.092) -- ++(0.042,-0.075);
  \draw[vel] (4.290,1.103) -- ++(-0.101,-0.099);
  \draw[vel] (3.743,0.013) -- ++(0.043,-0.111);
  \draw[vel] (4.313,0.437) -- ++(-0.029,-0.067);
  \draw[vel] (4.103,0.172) -- ++(0.010,-0.079);
  \draw[vel] (4.586,0.448) -- ++(-0.042,-0.038);
  \draw[vel] (3.631,0.220) -- ++(0.025,-0.133);
  \draw[vel] (4.219,-0.151) -- ++(0.040,-0.052);
  \draw[vel] (4.343,-0.076) -- ++(0.026,-0.041);
  \draw[vel] (4.735,0.323) -- ++(-0.035,-0.016);
  \draw[vel] (4.249,0.025) -- ++(0.019,-0.056);
  \draw[vel] (4.135,-0.509) -- ++(0.083,-0.045);
  \draw[vel] (3.578,0.388) -- ++(0.009,-0.146);
  \draw[vel] (3.240,0.538) -- ++(0.007,-0.190);
  \draw[vel] (3.735,0.546) -- ++(-0.015,-0.136);
  \draw[vel] (4.272,-0.334) -- ++(0.058,-0.038);
  \draw[vel] (4.887,0.798) -- ++(-0.094,-0.020);
  \draw[vel] (4.164,0.146) -- ++(0.010,-0.071);
  \draw[vel] (4.112,-0.121) -- ++(0.041,-0.065);
  \draw[vel] (4.804,0.044) -- ++(-0.007,0.004);
  \draw[vel] (4.172,-0.051) -- ++(0.031,-0.061);
  \draw[vel] (3.856,-0.236) -- ++(0.065,-0.088);
  \draw[vel] (4.891,0.191) -- ++(-0.027,0.007);
  \draw[vel] (4.731,0.255) -- ++(-0.027,-0.013);
  \draw[vel] (3.818,0.126) -- ++(0.027,-0.108);
  \draw[vel] (5.098,-0.296) -- ++(0.017,0.051);
  \draw[vel] (5.181,-0.465) -- ++(0.032,0.068);
  \draw[vel] (4.730,-0.744) -- ++(0.083,0.030);
  \draw[vel] (6.094,-0.669) -- ++(0.014,0.177);
  \draw[vel] (4.876,-0.114) -- ++(0.007,0.019);
  \draw[vel] (5.983,-1.100) -- ++(0.067,0.184);
  \draw[vel] (5.256,-0.648) -- ++(0.049,0.084);
  \draw[vel] (4.558,0.104) -- ++(-0.003,-0.026);
  \draw[vel] (5.339,-0.261) -- ++(0.003,0.076);
  \draw[vel] (5.011,-0.050) -- ++(-0.006,0.031);
  \draw[vel] (5.107,-0.379) -- ++(0.026,0.056);
  \draw[vel] (4.851,-0.969) -- ++(0.102,0.054);
  \draw[vel] (5.951,-0.579) -- ++(0.011,0.158);
  \draw[vel] (5.481,-0.319) -- ++(0.003,0.094);
  \draw[vel] (5.151,-0.579) -- ++(0.046,0.070);
  \draw[vel] (5.634,-0.466) -- ++(0.012,0.118);
  \draw[vel] (5.282,-0.288) -- ++(0.008,0.071);
  \draw[vel] (5.879,0.074) -- ++(-0.058,0.121);
  \draw[vel] (5.522,-0.610) -- ++(0.033,0.112);
  \draw[vel] (4.728,0.222) -- ++(-0.023,-0.012);
  \draw[vel] (5.785,-0.377) -- ++(-0.004,0.130);
  \draw[vel] (5.594,0.130) -- ++(-0.051,0.087);
  \draw[vel] (5.724,0.207) -- ++(-0.066,0.098);
  \draw[vel] (5.145,0.533) -- ++(-0.076,0.020);
  \draw[vel] (4.789,0.174) -- ++(-0.021,-0.003);
  \draw[vel] (5.572,0.180) -- ++(-0.056,0.083);
  \draw[vel] (6.196,0.516) -- ++(-0.120,0.136);
  \draw[vel] (4.835,-1.229) -- ++(0.131,0.063);
  \draw[vel] (5.718,-0.858) -- ++(0.052,0.144);
  \draw[vel] (5.344,0.162) -- ++(-0.044,0.058);
  \draw[vel] (5.467,-0.276) -- ++(-0.001,0.091);
  \draw[vel] (5.239,-0.279) -- ++(0.009,0.066);
  \draw[vel] (4.963,-1.132) -- ++(0.115,0.073);
  \draw[vel] (5.275,-0.834) -- ++(0.069,0.094);
  \draw[vel] (4.610,-0.022) -- ++(0.009,-0.014);
  \draw[vel] (5.322,-0.076) -- ++(-0.017,0.066);
  \draw[vel] (4.905,-0.662) -- ++(0.066,0.046);
  \draw[vel] (4.900,-0.788) -- ++(0.080,0.051);
  \draw[vel] (5.438,-0.731) -- ++(0.050,0.108);
  \draw[vel] (5.510,-0.113) -- ++(-0.021,0.089);
  \draw[vel] (6.261,-1.066) -- ++(0.051,0.213);
  \draw[vel] (5.750,-0.349) -- ++(-0.006,0.125);
  \draw[vel] (5.344,-1.097) -- ++(0.095,0.114);
  \draw[vel] (5.143,0.109) -- ++(-0.029,0.038);
  \draw[vel] (5.313,-0.255) -- ++(0.003,0.073);
  \draw[vel] (5.219,0.335) -- ++(-0.057,0.037);
  \draw[vel] (3.887,-0.592) -- ++(0.103,-0.069);
  \draw[vel] (5.950,-0.274) -- ++(-0.023,0.144);
  \draw[vel] (4.822,-0.817) -- ++(0.087,0.044);
  \draw[vel] (5.859,-0.213) -- ++(-0.025,0.131);
  \draw[vel] (5.372,-0.329) -- ++(0.009,0.083);
  \draw[vel] (5.367,0.143) -- ++(-0.043,0.062);
  \draw[vel] (5.599,-0.181) -- ++(-0.017,0.101);
  \draw[vel] (4.881,-0.019) -- ++(-0.004,0.015);
  \draw[vel] (5.042,0.302) -- ++(-0.046,0.019);
  \draw[vel] (4.778,-0.376) -- ++(0.040,0.020);
  \draw[vel] (5.347,-1.029) -- ++(0.087,0.111);
  \fill[agent] (3.662,-0.057) circle (1.0pt);
  \fill[agent] (4.783,-0.057) circle (1.0pt);
  \fill[agent] (4.182,0.586) circle (1.0pt);
  \fill[agent] (3.879,0.208) circle (1.0pt);
  \fill[agent] (4.261,0.274) circle (1.0pt);
  \fill[agent] (3.526,0.279) circle (1.0pt);
  \fill[agent] (4.947,-0.338) circle (1.0pt);
  \fill[agent] (4.673,0.295) circle (1.0pt);
  \fill[agent] (3.847,1.010) circle (1.0pt);
  \fill[agent] (4.619,-0.206) circle (1.0pt);
  \fill[agent] (4.241,0.469) circle (1.0pt);
  \fill[agent] (4.096,0.510) circle (1.0pt);
  \fill[agent] (4.163,0.504) circle (1.0pt);
  \fill[agent] (4.991,-0.007) circle (1.0pt);
  \fill[agent] (4.312,0.074) circle (1.0pt);
  \fill[agent] (4.270,-0.201) circle (1.0pt);
  \fill[agent] (3.881,0.175) circle (1.0pt);
  \fill[agent] (4.694,0.685) circle (1.0pt);
  \fill[agent] (3.472,-0.052) circle (1.0pt);
  \fill[agent] (4.556,-0.507) circle (1.0pt);
  \fill[agent] (3.945,0.213) circle (1.0pt);
  \fill[agent] (4.891,0.512) circle (1.0pt);
  \fill[agent] (4.020,0.110) circle (1.0pt);
  \fill[agent] (4.062,0.829) circle (1.0pt);
  \fill[agent] (3.965,0.135) circle (1.0pt);
  \fill[agent] (4.394,0.204) circle (1.0pt);
  \fill[agent] (4.091,-0.173) circle (1.0pt);
  \fill[agent] (4.194,0.081) circle (1.0pt);
  \fill[agent] (4.841,0.498) circle (1.0pt);
  \fill[agent] (4.187,0.504) circle (1.0pt);
  \fill[agent] (4.013,0.650) circle (1.0pt);
  \fill[agent] (4.197,0.472) circle (1.0pt);
  \fill[agent] (3.490,0.382) circle (1.0pt);
  \fill[agent] (3.271,-0.523) circle (1.0pt);
  \fill[agent] (4.033,-0.092) circle (1.0pt);
  \fill[agent] (4.290,1.103) circle (1.0pt);
  \fill[agent] (3.743,0.013) circle (1.0pt);
  \fill[agent] (4.313,0.437) circle (1.0pt);
  \fill[agent] (4.103,0.172) circle (1.0pt);
  \fill[agent] (4.586,0.448) circle (1.0pt);
  \fill[agent] (3.631,0.220) circle (1.0pt);
  \fill[agent] (4.219,-0.151) circle (1.0pt);
  \fill[agent] (4.343,-0.076) circle (1.0pt);
  \fill[agent] (4.735,0.323) circle (1.0pt);
  \fill[agent] (4.249,0.025) circle (1.0pt);
  \fill[agent] (4.135,-0.509) circle (1.0pt);
  \fill[agent] (3.578,0.388) circle (1.0pt);
  \fill[agent] (3.240,0.538) circle (1.0pt);
  \fill[agent] (3.735,0.546) circle (1.0pt);
  \fill[agent] (4.272,-0.334) circle (1.0pt);
  \fill[agent] (4.887,0.798) circle (1.0pt);
  \fill[agent] (4.164,0.146) circle (1.0pt);
  \fill[agent] (4.112,-0.121) circle (1.0pt);
  \fill[agent] (4.804,0.044) circle (1.0pt);
  \fill[agent] (4.172,-0.051) circle (1.0pt);
  \fill[agent] (3.856,-0.236) circle (1.0pt);
  \fill[agent] (4.891,0.191) circle (1.0pt);
  \fill[agent] (4.731,0.255) circle (1.0pt);
  \fill[agent] (3.818,0.126) circle (1.0pt);
  \fill[agent] (5.098,-0.296) circle (1.0pt);
  \fill[agent] (5.181,-0.465) circle (1.0pt);
  \fill[agent] (4.730,-0.744) circle (1.0pt);
  \fill[agent] (6.094,-0.669) circle (1.0pt);
  \fill[agent] (4.876,-0.114) circle (1.0pt);
  \fill[agent] (5.983,-1.100) circle (1.0pt);
  \fill[agent] (5.256,-0.648) circle (1.0pt);
  \fill[agent] (4.558,0.104) circle (1.0pt);
  \fill[agent] (5.339,-0.261) circle (1.0pt);
  \fill[agent] (5.011,-0.050) circle (1.0pt);
  \fill[agent] (5.107,-0.379) circle (1.0pt);
  \fill[agent] (4.851,-0.969) circle (1.0pt);
  \fill[agent] (5.951,-0.579) circle (1.0pt);
  \fill[agent] (5.481,-0.319) circle (1.0pt);
  \fill[agent] (5.151,-0.579) circle (1.0pt);
  \fill[agent] (5.634,-0.466) circle (1.0pt);
  \fill[agent] (5.282,-0.288) circle (1.0pt);
  \fill[agent] (5.879,0.074) circle (1.0pt);
  \fill[agent] (5.522,-0.610) circle (1.0pt);
  \fill[agent] (4.728,0.222) circle (1.0pt);
  \fill[agent] (5.785,-0.377) circle (1.0pt);
  \fill[agent] (5.594,0.130) circle (1.0pt);
  \fill[agent] (5.724,0.207) circle (1.0pt);
  \fill[agent] (5.145,0.533) circle (1.0pt);
  \fill[agent] (4.789,0.174) circle (1.0pt);
  \fill[agent] (5.572,0.180) circle (1.0pt);
  \fill[agent] (6.196,0.516) circle (1.0pt);
  \fill[agent] (4.835,-1.229) circle (1.0pt);
  \fill[agent] (5.718,-0.858) circle (1.0pt);
  \fill[agent] (5.344,0.162) circle (1.0pt);
  \fill[agent] (5.467,-0.276) circle (1.0pt);
  \fill[agent] (5.239,-0.279) circle (1.0pt);
  \fill[agent] (4.963,-1.132) circle (1.0pt);
  \fill[agent] (5.275,-0.834) circle (1.0pt);
  \fill[agent] (4.610,-0.022) circle (1.0pt);
  \fill[agent] (5.322,-0.076) circle (1.0pt);
  \fill[agent] (4.905,-0.662) circle (1.0pt);
  \fill[agent] (4.900,-0.788) circle (1.0pt);
  \fill[agent] (5.438,-0.731) circle (1.0pt);
  \fill[agent] (5.510,-0.113) circle (1.0pt);
  \fill[agent] (6.261,-1.066) circle (1.0pt);
  \fill[agent] (5.750,-0.349) circle (1.0pt);
  \fill[agent] (5.344,-1.097) circle (1.0pt);
  \fill[agent] (5.143,0.109) circle (1.0pt);
  \fill[agent] (5.313,-0.255) circle (1.0pt);
  \fill[agent] (5.219,0.335) circle (1.0pt);
  \fill[agent] (3.887,-0.592) circle (1.0pt);
  \fill[agent] (5.950,-0.274) circle (1.0pt);
  \fill[agent] (4.822,-0.817) circle (1.0pt);
  \fill[agent] (5.859,-0.213) circle (1.0pt);
  \fill[agent] (5.372,-0.329) circle (1.0pt);
  \fill[agent] (5.367,0.143) circle (1.0pt);
  \fill[agent] (5.599,-0.181) circle (1.0pt);
  \fill[agent] (4.881,-0.019) circle (1.0pt);
  \fill[agent] (5.042,0.302) circle (1.0pt);
  \fill[agent] (4.778,-0.376) circle (1.0pt);
  \fill[agent] (5.347,-1.029) circle (1.0pt);
  \draw[frame] (7.850,-1.45) rectangle (11.150,1.45);
  \node[paneltitle] at (9.500,1.72) {population $\mu_t$};
  \node[panelsub, align=center] at (9.500,-1.75) {$\partial_t\mu_t-\mathrm{div}(A_{\mu_t}\mu_t)=0$ \\ ($m_{\mu_t}$ denotes the barycenter) };
  \begin{scope}
  \clip (7.850,-1.45) rectangle (11.150,1.45);
  \fill[dens, path fading=mfbump, fill opacity=0.75] (8.950,0.25) ellipse (1.150 and 0.800);
  \fill[dens, path fading=mfbump, fill opacity=0.75] (10.100,-0.3) ellipse (0.950 and 1.150);
  \draw[stream] (10.145,0.273) -- (10.118,0.299) -- (10.091,0.323) -- (10.063,0.345) -- (10.035,0.365) -- (10.006,0.384) -- (9.977,0.401) -- (9.947,0.416) -- (9.918,0.429) -- (9.888,0.441) -- (9.859,0.451);
  \draw[stream] (9.227,0.645) -- (9.201,0.618) -- (9.177,0.591) -- (9.155,0.563) -- (9.135,0.535) -- (9.116,0.506) -- (9.099,0.477) -- (9.084,0.447) -- (9.071,0.418) -- (9.059,0.388) -- (9.049,0.359);
  \draw[stream] (8.855,-0.273) -- (8.882,-0.299) -- (8.909,-0.323) -- (8.937,-0.345) -- (8.965,-0.365) -- (8.994,-0.384) -- (9.023,-0.401) -- (9.053,-0.416) -- (9.082,-0.429) -- (9.112,-0.441) -- (9.141,-0.451);
  \draw[stream] (9.773,-0.645) -- (9.799,-0.618) -- (9.823,-0.591) -- (9.845,-0.563) -- (9.865,-0.535) -- (9.884,-0.506) -- (9.901,-0.477) -- (9.916,-0.447) -- (9.929,-0.418) -- (9.941,-0.388) -- (9.951,-0.359);
  \draw[stream] (10.697,0.506) -- (10.648,0.554) -- (10.598,0.599) -- (10.546,0.640) -- (10.493,0.678) -- (10.440,0.713) -- (10.386,0.744) -- (10.331,0.772) -- (10.276,0.797) -- (10.221,0.819) -- (10.167,0.837);
  \draw[stream] (9.851,1.252) -- (9.782,1.243) -- (9.716,1.232) -- (9.651,1.217) -- (9.589,1.199) -- (9.528,1.179) -- (9.470,1.156) -- (9.414,1.131) -- (9.361,1.104) -- (9.310,1.075) -- (9.261,1.043);
  \draw[stream] (8.740,1.055) -- (8.704,0.996) -- (8.672,0.937) -- (8.643,0.877) -- (8.618,0.817) -- (8.596,0.757) -- (8.577,0.698) -- (8.562,0.638) -- (8.550,0.579) -- (8.541,0.521) -- (8.535,0.464);
  \draw[stream] (8.202,0.063) -- (8.225,-0.001) -- (8.251,-0.064) -- (8.280,-0.123) -- (8.311,-0.180) -- (8.344,-0.235) -- (8.379,-0.286) -- (8.416,-0.335) -- (8.455,-0.381) -- (8.495,-0.425) -- (8.536,-0.465);
  \draw[stream] (8.641,-0.976) -- (8.706,-0.998) -- (8.771,-1.016) -- (8.836,-1.031) -- (8.900,-1.042) -- (8.963,-1.050) -- (9.025,-1.055) -- (9.086,-1.056) -- (9.146,-1.055) -- (9.205,-1.051) -- (9.263,-1.044);
  \draw[stream] (9.727,-1.280) -- (9.785,-1.243) -- (9.840,-1.203) -- (9.892,-1.162) -- (9.940,-1.119) -- (9.986,-1.075) -- (10.029,-1.029) -- (10.068,-0.982) -- (10.104,-0.934) -- (10.138,-0.886) -- (10.168,-0.836);
  \draw[stream] (10.642,-0.621) -- (10.649,-0.552) -- (10.653,-0.485) -- (10.653,-0.418) -- (10.650,-0.353) -- (10.643,-0.290) -- (10.634,-0.228) -- (10.622,-0.168) -- (10.607,-0.110) -- (10.590,-0.054) -- (10.570,0.001);
  \end{scope}
  \fill[agent] (9.500,0) circle (1.3pt) node[below right=-1pt, font=\scriptsize, text=agentcol] {$m_{\mu_t}$};
  \draw[link] (1.750,0) -- (3.000,0) node[midway, above=3pt, font=\scriptsize, fill=white, inner sep=1pt] {$N$ grows};
  \draw[link] (6.500,0) -- (7.750,0) node[midway, above=3pt, font=\scriptsize, fill=white, inner sep=1pt] {$N\longrightarrow+\infty$};
\end{tikzpicture}
\caption{Each agent moves with the velocity $-A_{\mu^N_t}$ (arrows), which is evaluated against the empirical measure of the whole population, as in \eqref{particles}. As $N\longrightarrow+\infty$, the empirical measures are replaced by a probability measure $\mu_t$ (shaded) transported by the velocity field $-A_{\mu_t}$ it generates itself (streamlines).}
\label{fig:mean_field}
\end{figure}
%%%%%%%%%%%%%%%%%%%%
Each empirical measure can be proven to satisfy the continuity equation
\[
  \partial_t\mu^N_t-\mathrm{div}\bigl(A_{\mu^N_t}\mu^N_t\bigr)=0
\]
in the sense of distributions. Therefore, as $N\longrightarrow+\infty$, one may expect, and this is the case under suitable assumptions, that the empirical measures converge to a probability measure $\mu_t$ that is transported by the velocity field $-A_{\mu_t}$ that it itself generates, i.e., to a solution of
\begin{equation}\label{mf_flow}
  \partial_t\mu_t-\mathrm{div}\bigl(A_{\mu_t}\mu_t\bigr)=0.
\end{equation}
This continuity equation is also understood in the sense of distributions. Figure~\ref{fig:mean_field} illustrates the convergence process. The natural ambient space for \eqref{mf_flow} is the so-called Wasserstein space $(\mathcal P_2(\mathbb R^d),\mathcal W_2)$. This is mainly because its formal Riemannian structure allows one to read continuity equations as ordinary differential equations on a curved space of measures; see \cite[Chapter~8]{AGS_2008} and \cite[Section~1]{O_2001}. In this idealization, \eqref{mf_flow} plays the role of the evolution equation \eqref{eveq} in the classical single-particle case, with the difference that the space in which the velocity field lives changes with the measure.
\iffalse
\smallbreak\noindent
The classical evolution equation is contained in the mean-field framework. If $A_\mu\equiv M$, with $M:\mathbb R^d\to\mathbb R^d$ a continuous monotone vector field, then the agents do not interact, and the solution of \eqref{mf_flow} is obtained by letting each agent follow, independently of the others, the classical flow generated by $-M$. Restricting \eqref{mf_flow} to Dirac masses then recovers the classical Cauchy problem for a monotone operator on $\mathbb R^d$, so that every difficulty of the classical problem is already in the mean-field regime.
\fi 
\smallbreak\noindent
As in the single-particle case, the interest lies in the long-time behavior of \eqref{mf_flow}. One expects $\mu_t$ to approach a zero of the operator, i.e., a measure $\mu$ with $A_\mu=0$ $\mu$-almost everywhere, and the questions are under which notion of monotonicity this happens and at which rate. Monotonicity on $\mathcal P_2(\mathbb R^d)$ is understood through optimal transport plans, i.e., we say an operator $A$ is \emph{monotone} if
\[
  \int\langle A_\mu(x)-A_\nu(y),x-y\rangle\,d\gamma(x,y)\geq0
  \quad\text{for all }\mu,\nu\in\mathcal P_2(\mathbb R^d)\text{ and }\gamma\in\Gamma_o(\mu,\nu).
\]
Here $\Gamma_o(\mu, \nu)$ denotes the set of optimal couplings between $\mu$ and $\nu$.
% We denote by $\mathrm{zer}(A)$ the set of zeroes of the operator $A$. 
In the mean-field context, we may regard a zero as a state of the population in which every agent is at rest. For the Wasserstein gradient of a displacement convex energy, the zeroes are the minimizers of the energy, while for the saddle operators of zero-sum games with convex--concave payoffs, played in mixed strategies, they are the mixed Nash equilibria.
\smallbreak\noindent
To quantify at which rate $\mu_t$ comes to rest, one needs an analogue of the residual $|A(x(t))|$ from the single-particle case. In the mean-field setting there are two natural quantities that measure the remaining motion along the flow. 
The first one is the operator residual \(\|A_{\mu_t}\|_{L^2(\mu_t)^d}\), which measures the average magnitude of the agents' velocities and vanishes exactly when \(\mu_t\in\mathrm{zer}(A)\), where $\mathrm{zer}(A)$ is the set of zeroes of the operator $A$. 
The second one is the metric speed $|\mu'|(t)$, which is the speed at which the population as a whole moves in $(\mathcal P_2(\mathbb R^d),\mathcal W_2)$. For a single particle both quantities reduce to $|\dot x(t)|=|A(x(t))|$, but in general they differ. The metric speed does not exceed the residual for almost every $t$, and the inequality can be strict, since the agents may move while their distribution stays the same. 
\smallbreak\noindent
When the operator is the Wasserstein gradient of a displacement convex energy, the theory of gradient flows in Wasserstein spaces yields the decay of the residual whenever the energy attains its minimum; see \cite[Section~11.2]{AGS_2008}. For a general monotone operator, however, there might be no energy to rely on, and the residual and the metric speed need not vanish in the long run.
\smallbreak\noindent
Non-gradient measure-valued evolutions have been studied in \cite{AC_2014,BF_2021,BF_2025_corr,CCDP_2021,CSS_2023,CSS_2026,CHS_2026,P_2019}. In mean field games, the monotonicity condition above is known as displacement monotonicity \cite{A_2016,GM_2022,GMMZ_2022,MM_2024}, which is distinct from Lasry--Lions monotonicity \cite{GM_2023,LL_2007}. Wasserstein dynamics of this kind arise in descent--ascent schemes and multispecies games \cite{AL_2025,CHMR_2024,CHMR_2025,IS_2025,WC_2026}. Some of the aforementioned works are discussed further in Section~\ref{s2}.

\subsubsection*{Failure of convergence to zeroes}
It can be seen that the Wasserstein distance from $\mu_t$ to any zero of a monotone operator $A$ is nonincreasing along \eqref{mf_flow}; in other words, the trajectory is Fej\'er monotone with respect to $\mathrm{zer}(A)$. Although monotonicity provides some control over the dynamics, it forces neither the trajectory to converge to a zero nor the residual to vanish. This issue is already present in the classical single-particle case. The bilinear function $L(x_1,x_2):=x_1x_2$ on $\mathbb R\times\mathbb R$ gives rise to the saddle operator $A(x)=Jx$, with
\[
  J:=\begin{pmatrix}0&1\\-1&0\end{pmatrix}.
\]
This operator is monotone, because $J\in\mathbb R^{2\times 2}$ is skew-symmetric, and its only zero is the origin. Yet the trajectories of the associated evolution equation circle the origin forever, at constant distance and with constant speed, so that neither the distance to the zero nor the residual decays. The explicit Euler scheme behaves even worse. Since $|x-\tau Jx|^2=(1+\tau^2)|x|^2$, its iterates spiral away from the origin for every step size $\tau>0$. This is the well-known cycling of gradient descent--ascent on bilinear games. Passing to continuous time removes the outward drift, but the trajectories still fail to converge.
\smallbreak\noindent
The same phenomenon occurs in the mean-field regime. Denoting by $m_\mu\in\mathbb R^2$ the barycenter of a probability measure $\mu\in\mathcal P_2(\mathbb R^2)$, consider the operator
\begin{equation}\label{rotation}
  A_\mu(x):=J\,m_\mu\quad\text{for } \mu\in\mathcal P_2(\mathbb R^2)\text{ and } x\in\mathbb R^2 .
\end{equation}
This operator is the mean-field version of the previous saddle operator. Here the two marginals of $\mu$ are understood as the mixed strategies of the minimizing and the maximizing player, and $A$ is obtained from the expected payoff, in the same way as the saddle operator is obtained from $L$. Thus \eqref{mf_flow} is a mean-field form of gradient descent--ascent. It is not hard to check that the operator is monotone and that its zeroes are exactly the measures with barycenter at the origin.
\smallbreak\noindent
Along the flow all agents move with the same velocity, so the population is simply translated, and its barycenter follows the rotation above, circling the origin at a constant distance.
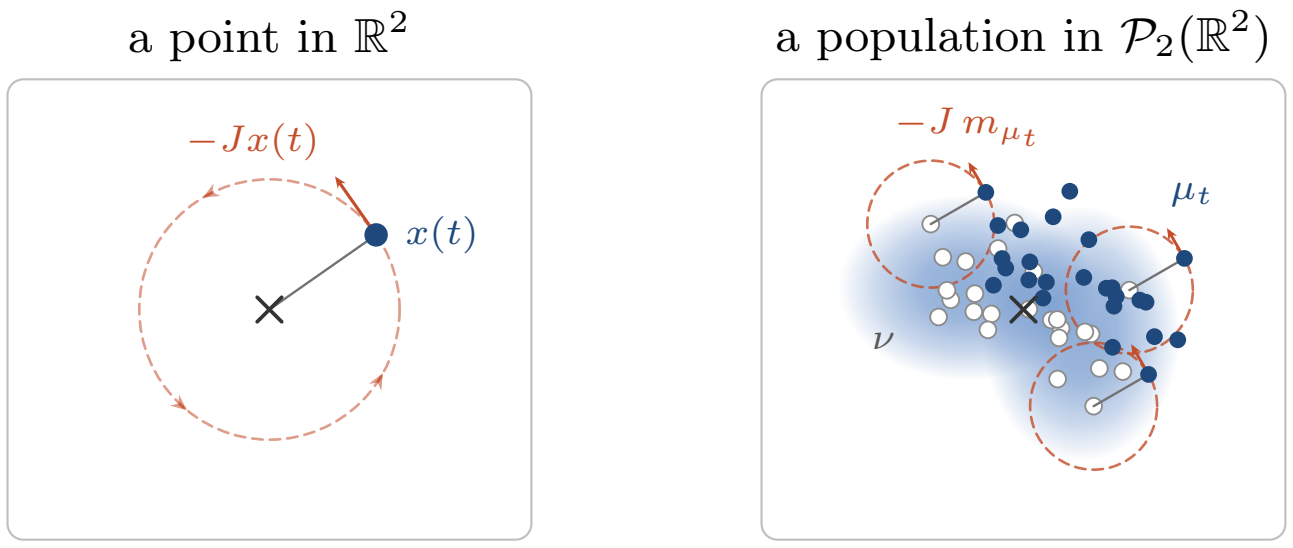
\begin{figure}[H]
\centering
\definecolor{agentcol}{RGB}{31,73,125}
\definecolor{velcol}{RGB}{196,78,42}
\definecolor{denscol}{RGB}{70,120,190}
\pgfdeclareradialshading{failbumpshade}{\pgfpointorigin}{%
  color(0bp)=(pgftransparent!0); color(9bp)=(pgftransparent!12);
  color(17bp)=(pgftransparent!50); color(25bp)=(pgftransparent!100);
  color(50bp)=(pgftransparent!100)}%
\pgfdeclarefading{failbump}{\pgfuseshading{failbumpshade}}
\begin{tikzpicture}[
  scale=1.05, line cap=round,
  frame/.style={draw=black!25, rounded corners=3pt, line width=0.4pt},
  paneltitle/.style={font=\small},
  panelsub/.style={font=\scriptsize, text=black!75},
  agent/.style={fill=agentcol},
  vel/.style={-{Stealth[length=2.2pt, width=1.8pt]}, draw=velcol, line width=0.45pt},
  dens/.style={fill=denscol},
  orbit/.style={draw=velcol, opacity=0.55, line width=0.5pt, dash pattern=on 2pt off 1.3pt},
  link0/.style={draw=black!55, line width=0.45pt},
  zero/.style={draw=black!80, line width=0.7pt},
  link/.style={-{Stealth[length=4pt]}, line width=0.6pt, draw=black!60}
]
  \draw[frame] (-1.650,-1.450) rectangle (1.650,1.450);
  \node[paneltitle] at (0.000,1.720) {a point in $\mathbb R^2$};
  \draw[orbit] (0,0) circle (0.82);
  \draw[orbit, -{Stealth[length=3pt, width=2.6pt]}] (-0.401,0.715) -- (-0.419,0.705);
  \draw[orbit, -{Stealth[length=3pt, width=2.6pt]}] (-0.535,-0.622) -- (-0.519,-0.635);
  \draw[orbit, -{Stealth[length=3pt, width=2.6pt]}] (0.705,-0.419) -- (0.715,-0.401);
  \draw[link0] (0,0) -- (0.672,0.470);
  \draw[vel, line width=0.6pt] (0.672,0.470) -- ++(-0.259,0.369) node[pos=1, above left=-1pt, font=\scriptsize, text=velcol] {$-Jx(t)$};
  \fill[agent] (0.672,0.470) circle (2.1pt) node[right=2pt, font=\scriptsize, text=agentcol] {$x(t)$};
  \draw[zero] (0.000,0.000) +(-2.2pt,-2.2pt) -- +(2.2pt,2.2pt) +(-2.2pt,2.2pt) -- +(2.2pt,-2.2pt);
  \draw[frame] (3.100,-1.450) rectangle (6.400,1.450);
  \node[paneltitle] at (4.750,1.720) {a population in $\mathcal P_2(\mathbb R^2)$};
  \begin{scope}
  \clip (3.100,-1.450) rectangle (6.400,1.450);
  \fill[dens, path fading=failbump, fill opacity=0.6] (4.390,0.14) ellipse (0.8 and 0.58);
  \fill[dens, path fading=failbump, fill opacity=0.6] (5.110,-0.15) ellipse (0.62 and 0.8);
  \end{scope}
  \draw[orbit, opacity=0.8] (4.166,0.538) circle (0.4);
  \draw[orbit, opacity=0.8] (5.192,-0.609) circle (0.4);
  \draw[orbit, opacity=0.8] (5.417,0.122) circle (0.4);
  \draw[black!45, line width=0.35pt, fill=white] (4.215,-0.047) circle (1.5pt);
  \draw[black!45, line width=0.35pt, fill=white] (4.386,0.302) circle (1.5pt);
  \draw[black!45, line width=0.35pt, fill=white] (4.815,0.240) circle (1.5pt);
  \draw[black!45, line width=0.35pt, fill=white] (4.292,0.061) circle (1.5pt);
  \draw[black!45, line width=0.35pt, fill=white] (4.695,0.545) circle (1.5pt);
  \draw[black!45, line width=0.35pt, fill=white] (4.548,-0.028) circle (1.5pt);
  \draw[black!45, line width=0.35pt, fill=white] (4.166,0.538) circle (1.5pt);
  \draw[black!45, line width=0.35pt, fill=white] (4.526,-0.128) circle (1.5pt);
  \draw[black!45, line width=0.35pt, fill=white] (4.437,-0.014) circle (1.5pt);
  \draw[black!45, line width=0.35pt, fill=white] (4.268,0.121) circle (1.5pt);
  \draw[black!45, line width=0.35pt, fill=white] (4.242,0.329) circle (1.5pt);
  \draw[black!45, line width=0.35pt, fill=white] (4.444,0.100) circle (1.5pt);
  \draw[black!45, line width=0.35pt, fill=white] (4.590,0.384) circle (1.5pt);
  \draw[black!45, line width=0.35pt, fill=white] (4.987,-0.118) circle (1.5pt);
  \draw[black!45, line width=0.35pt, fill=white] (4.925,-0.066) circle (1.5pt);
  \draw[black!45, line width=0.35pt, fill=white] (5.176,-0.154) circle (1.5pt);
  \draw[black!45, line width=0.35pt, fill=white] (4.974,-0.179) circle (1.5pt);
  \draw[black!45, line width=0.35pt, fill=white] (4.965,-0.438) circle (1.5pt);
  \draw[black!45, line width=0.35pt, fill=white] (5.228,-0.371) circle (1.5pt);
  \draw[black!45, line width=0.35pt, fill=white] (5.417,0.122) circle (1.5pt);
  \draw[black!45, line width=0.35pt, fill=white] (4.784,0.002) circle (1.5pt);
  \draw[black!45, line width=0.35pt, fill=white] (4.963,-0.063) circle (1.5pt);
  \draw[black!45, line width=0.35pt, fill=white] (5.192,-0.609) circle (1.5pt);
  \draw[black!45, line width=0.35pt, fill=white] (5.137,-0.141) circle (1.5pt);
  \draw[black!45, line width=0.35pt, fill=white] (5.376,-0.391) circle (1.5pt);
  \draw[link0] (4.166,0.538) -- (4.513,0.738);
  \draw[link0] (5.192,-0.609) -- (5.538,-0.409);
  \draw[link0] (5.417,0.122) -- (5.764,0.322);
  \draw[vel, line width=0.55pt] (4.513,0.738) -- ++(-0.110,0.191);
  \draw[vel, line width=0.55pt] (5.538,-0.409) -- ++(-0.110,0.191);
  \draw[vel, line width=0.55pt] (5.764,0.322) -- ++(-0.110,0.191);
  \fill[agent] (4.561,0.153) circle (1.5pt);
  \fill[agent] (4.733,0.502) circle (1.5pt);
  \fill[agent] (5.162,0.440) circle (1.5pt);
  \fill[agent] (4.638,0.261) circle (1.5pt);
  \fill[agent] (5.042,0.745) circle (1.5pt);
  \fill[agent] (4.894,0.172) circle (1.5pt);
  \fill[agent] (4.513,0.738) circle (1.5pt);
  \fill[agent] (4.873,0.072) circle (1.5pt);
  \fill[agent] (4.784,0.186) circle (1.5pt);
  \fill[agent] (4.615,0.321) circle (1.5pt);
  \fill[agent] (4.589,0.529) circle (1.5pt);
  \fill[agent] (4.790,0.300) circle (1.5pt);
  \fill[agent] (4.937,0.584) circle (1.5pt);
  \fill[agent] (5.334,0.082) circle (1.5pt);
  \fill[agent] (5.272,0.134) circle (1.5pt);
  \fill[agent] (5.522,0.046) circle (1.5pt);
  \fill[agent] (5.320,0.021) circle (1.5pt);
  \fill[agent] (5.311,-0.238) circle (1.5pt);
  \fill[agent] (5.575,-0.171) circle (1.5pt);
  \fill[agent] (5.764,0.322) circle (1.5pt);
  \fill[agent] (5.130,0.202) circle (1.5pt);
  \fill[agent] (5.309,0.137) circle (1.5pt);
  \fill[agent] (5.538,-0.409) circle (1.5pt);
  \fill[agent] (5.483,0.059) circle (1.5pt);
  \fill[agent] (5.722,-0.191) circle (1.5pt);
  \draw[zero] (4.750,0.000) +(-2.2pt,-2.2pt) -- +(2.2pt,2.2pt) +(-2.2pt,2.2pt) -- +(2.2pt,-2.2pt);
  \node[font=\scriptsize, text=black!65] at (3.870,-0.200) {$\nu$};
  \node[font=\scriptsize, text=agentcol] at (5.814,0.725) {$\mu_t$};
  \node[font=\scriptsize, text=velcol, above=-0.5pt] at (4.403,0.929) {$-J\,m_{\mu_t}$};
\end{tikzpicture}
\caption{On the left, the trajectory of $\dot x=-Jx$ circles the zero of $A(x)=Jx$ at constant distance. On the right, all agents move with the same velocity $-J\,m_{\mu_t}$, and each of them circles its counterpart in the zero $\nu:=(\id-m_{\mu_{t_0}})_\#\mu_{t_0}$ of $A_\mu=J\,m_\mu$. 
}
\label{fig:failure}
\end{figure}
Since translating the population so that its barycenter is at the origin produces the closest zero, a direct computation shows that the distance from $\mu_t$ to $\mathrm{zer}(A)$, the residual, and the metric speed all coincide with the distance from the barycenter to the origin, i.e.,
\[
  \mathcal W_2\bigl(\mu_t,\mathrm{zer}(A)\bigr)=\|A_{\mu_t}\|_{L^2(\mu_t)^d}=|\mu'|(t)=|m_{\mu_t}|=|m_{\mu_{t_0}}|.
\]
Consequently, unless the initial population has its barycenter at the origin, none of these quantities decays, so the population neither approaches the set of zeroes nor slows down. This is the mean-field counterpart of the cycling we saw for a single particle; Figure~\ref{fig:failure} illustrates this.
\subsubsection*{Cocoercive operators on the Wasserstein space}
For the operator \eqref{rotation} above, whose flow makes the barycenter rotate, the monotonicity condition holds with equality. The difference of the two vector fields is, on average, orthogonal to the displacement, so nothing is dissipated. As in Hilbert spaces, this is ruled out by requiring the left-hand side of the monotonicity condition to control the size of this difference, which leads to cocoercivity in the Wasserstein space. An operator $A$ is called \emph{$\beta$-cocoercive} if
\[
  \int\langle A_\mu(x)-A_\nu(y),x-y\rangle\,d\gamma(x,y)\geq\beta\int|A_\mu(x)-A_\nu(y)|^2\,d\gamma(x,y)
\]
for all $\mu,\nu\in\mathcal P_2(\mathbb R^d)$ and $\gamma\in\Gamma_o(\mu,\nu)$. Along the solution of (\ref{mf_flow}), cocoercivity turns the Fej\'er monotonicity of the trajectory into the dissipation estimate
\[
  \frac12\mathcal W_2^2(\mu_t,\nu)+\beta\int_{t_0}^t\|A_{\mu_s}\|_{L^2(\mu_s)^d}^2\,ds\leq\frac12\mathcal W_2^2(\mu_{t_0},\nu)
  \quad\text{for every } \nu\in\mathrm{zer}(A)\text{ and } t\geq t_0,
\]
so that the residual is square integrable in time, which is precisely what fails in the rotation example. This estimate is proved in Theorem~\ref{thm1}, which also shows that, as long as the operator has at least one zero, cocoercivity is enough to guarantee that the flow converges in the Wasserstein space to a zero of the operator, and that, moreover, both the residual and the metric speed decay as $o(t^{-1/2})$, so the agents eventually come to rest. 
This is the Wasserstein counterpart of Bruck's theorem \cite[Theorem~6.3]{PS_2010} and of the rates in \cite[Theorem~A.1]{BN_2026}.
We point out that the limit need not be a Dirac mass. When the flow converges to a zero that is a Dirac mass, we say that the agents have reached consensus.
\begin{figure}[t]
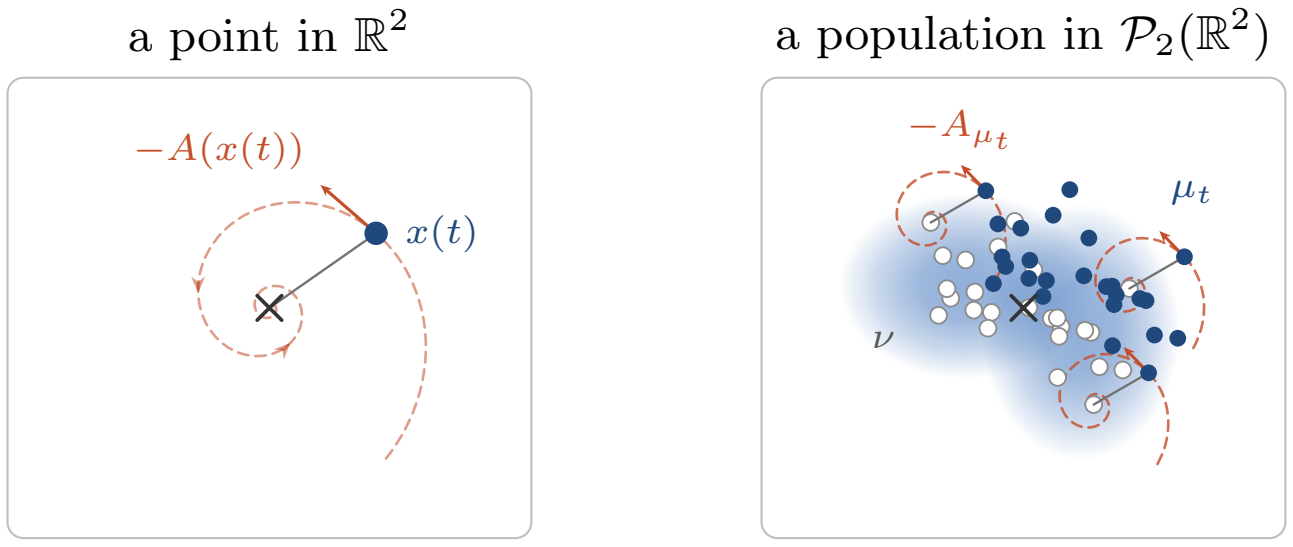

\centering
\definecolor{agentcol}{RGB}{31,73,125}
\definecolor{velcol}{RGB}{196,78,42}
\definecolor{denscol}{RGB}{70,120,190}
\pgfdeclareradialshading{cocobumpshade}{\pgfpointorigin}{%
  color(0bp)=(pgftransparent!0); color(9bp)=(pgftransparent!12);
  color(17bp)=(pgftransparent!50); color(25bp)=(pgftransparent!100);
  color(50bp)=(pgftransparent!100)}%
\pgfdeclarefading{cocobump}{\pgfuseshading{cocobumpshade}}

\caption{On the left, for the cocoercive operator $A(x)=\alpha x+Jx$ with $\alpha>0$, the trajectory of $\dot x=-A(x)$ spirals into the zero of $A$. On the right, for the cocoercive operator $A_\mu=\alpha\,m_\mu+J\,m_\mu$, each agent spirals into its counterpart in the zero $\nu:=(\id-m_{\mu_{t_0}})_\#\mu_{t_0}$ (dashed), so that $\mu_t$ converges to $\nu$, which is not a Dirac mass.}
\label{fig:cocoercive}
\end{figure}
\subsection{Tikhonov regularization and decay rates}

\subsubsection*{On Tikhonov regularization}
Cocoercivity is considerably stronger than monotonicity, and it excludes important examples, among them the saddle operator \eqref{rotation}. One way to recover it is to regularize. Penalizing the second moments of both mixed strategies replaces $A_\mu$ by $A_\mu+\kappa \id$ for some $\kappa>0$, and the resulting operator is cocoercive. This, however, changes the operator and, in general, its zeroes. For \eqref{rotation}, the regularized operator has $\delta_0$ as its only zero, whereas every measure with barycenter at the origin is a zero of the original one.
\smallbreak\noindent
The idea behind Tikhonov regularization is to let the weight $\kappa$ depend on time and vanish in the long term, so that the original operator is recovered asymptotically. To this end, we replace $\kappa$ by a nonincreasing function $\varepsilon:[t_0,+\infty)\to(0,+\infty)$ that vanishes as $t\longrightarrow+\infty$. This leads to the Tikhonov regularization of \eqref{mf_flow}, i.e., to the nonautonomous continuity equation
\begin{equation}\label{tik_flow}
  \partial_t\mu_t-\mathrm{div}\Bigl(\bigl(A_{\mu_t}+\varepsilon(t)\,\id\bigr)\mu_t\Bigr)=0.
\end{equation}
From the point of view of the agents, each of them still follows the negative of the operator, evaluated against the whole population, but it now also has a velocity component pulling it towards the origin, proportional to its distance from the origin with a factor $\varepsilon(t)$ that is the same for all agents and fades away as time goes on. At every time, the regularized operator $A+\varepsilon(t)\id$ is strongly monotone, i.e., its monotonicity condition holds with $\varepsilon(t) \mathlarger{\int}|x-y|^2\,d\gamma(x,y)$ in place of zero on the right-hand side. Thus, some dissipation is present even if the operator itself does not dissipate. The parameter $\varepsilon$ has to vanish, so that the population eventually settles at a zero of $A$ rather than at a zero of a perturbation of it, but slowly enough that $\mathlarger{\int}_{t_0}^{+\infty}\varepsilon(t)\,dt=+\infty$, so that the accumulated contraction makes the trajectory forget its initial datum.
\smallbreak\noindent
In the mean-field regime, the Tikhonov term has a further effect, which concerns the shape of the population and not merely its position, and which, as such, has no counterpart in the single-particle case. 
Since the factor $\varepsilon(t)$ is the same for all agents, the pull towards the origin brings the agents not only closer to the origin but also closer to each other, so that the population shrinks around its barycenter. By monotonicity, the operator does not move the agents away from their barycenter on average, i.e., it does not increase the variance $\mathlarger{\int}|x-m_{\mu_t}|^2\,d\mu_t(x)$ of the population, whereas the Tikhonov term contracts it, so that the variance of $\mu_t$ is at most $\exp\Bigl(-2 \mathlarger{\int}_{t_0}^t\varepsilon(s)\,ds\Bigr)$ times that of $\mu_{t_0}$. 
Thus, the divergence condition $\mathlarger{\int}_{t_0}^{+\infty}\varepsilon(t)\,dt=+\infty$ forces the variance of the population to vanish, i.e., the agents eventually gather around a single point. Since the population also comes to rest, one expects this point to be some $x\in\mathbb R^d$ with $\delta_x\in\mathrm{zer}(A)$, and the pull towards the origin to select, among these points, the one closest to the origin. We show that this is indeed the case. In Lemma~\ref{lemsecmoment} we prove that $A$ has a unique zero of minimal second moment, i.e., a unique zero closest to $\delta_0$, which plays the role of the zero of minimal norm, and that this zero is a Dirac mass. Theorem~\ref{thm2} then shows that the regularized flow converges to it. In contrast with Theorem~\ref{thm1}, the limit is therefore always a Dirac mass, and it does not depend on the initial datum, i.e., the regularized flow brings the population to consensus.
\smallbreak\noindent
The regularization also yields explicit decay rates; see Theorem~\ref{thm2} and Corollary~\ref{mccr}. For the parameter $\varepsilon(t)=\alpha/t$ with $\alpha\geq1$, Corollary~\ref{mccr} shows that the residual and the metric speed decay as $O(1/t)$, which is the rate known in the Hilbert space setting \cite[Corollary~2.10]{BN_2026}. Compared with Theorem~\ref{thm1}, this rate is faster than $o(t^{-1/2})$, and it holds for merely monotone demicontinuous operators, such as the saddle operator \eqref{rotation}, for which, without regularization, the population need not come to rest. The trade-off is that the dynamics are perturbed by a vanishing term. The rotation example also shows that the condition $\alpha\geq1$ cannot be dropped, since for $\alpha<1$ the residual decays more slowly than $1/t$. The other choices of $\varepsilon$ in Corollary~\ref{mccr} show that the choice $\alpha/t$ is a balanced one. If $\varepsilon$ vanishes more slowly, the metric speed still decays as $O(1/t)$, but the residual only decays as fast as $\varepsilon$ does. If $\varepsilon$ vanishes faster, as for $\varepsilon(t)=1/(t\log t)$, the regularization is weaker, and the residual and the metric speed decay only logarithmically.
\begin{figure}[t]
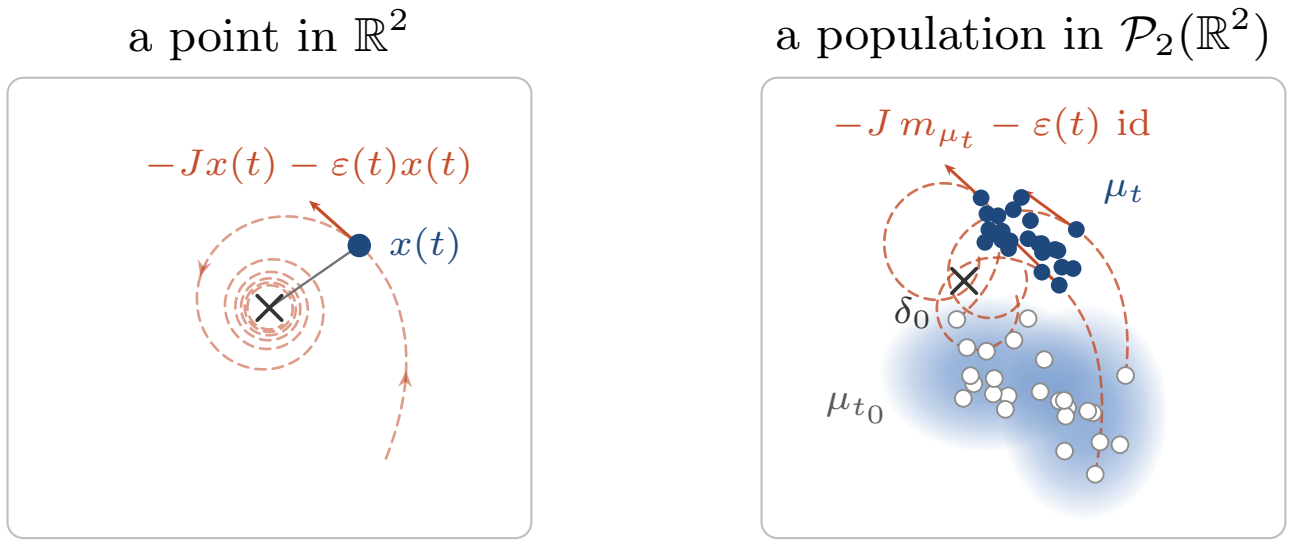

\centering
\definecolor{agentcol}{RGB}{31,73,125}
\definecolor{velcol}{RGB}{196,78,42}
\definecolor{denscol}{RGB}{70,120,190}
\pgfdeclareradialshading{tikbumpshade}{\pgfpointorigin}{%
  color(0bp)=(pgftransparent!0); color(9bp)=(pgftransparent!12);
  color(17bp)=(pgftransparent!50); color(25bp)=(pgftransparent!100);
  color(50bp)=(pgftransparent!100)}%
\pgfdeclarefading{tikbump}{\pgfuseshading{tikbumpshade}}

\caption{On the left, for $A(x)=Jx$ and $\varepsilon(t)=\alpha/t$ with $\alpha>0$, the trajectory of $\dot x=-Jx-\varepsilon(t)x$ spirals into the zero of $A$. On the right, for the Tikhonov regularization of \eqref{rotation} with the same $\varepsilon$, the barycenter of $\mu_t$ follows a trajectory of the system on the left, and the population shrinks around it without changing its shape, so that each agent spirals into the origin (dashed) and $\mu_t$ converges to $\delta_0$, the zero of \eqref{rotation} of minimal second moment.}
\label{fig:tikhonov}
\end{figure}
\subsubsection*{Nonexpansive operators and a Halpern-type algorithm}
The Halpern iteration for nonexpansive mappings \cite{H_1967} can be thought of as a discrete-time analogue of the Tikhonov regularized flow. It is obtained by an explicit discretization of the flow in which the regularization term is centered at a reference point instead of the origin; see \cite[Section~4.2]{BN_2026}. An operator $B$ on $\mathcal P_2(\mathbb R^d)$ is called \emph{nonexpansive} if
\[
  \int|B_\mu(x)-B_\nu(y)|^2\,d\gamma(x,y)\leq\int|x-y|^2\,d\gamma(x,y)
\]
for all $\mu,\nu\in\mathcal P_2(\mathbb R^d)$ and $\gamma\in\Gamma_o(\mu,\nu)$, and its fixed points are the measures $\mu$ with $B_\mu=\id$ $\mu$-almost everywhere. In Proposition~\ref{prop:wass_nonexpansive} we show that, as in Hilbert spaces, $B$ is nonexpansive if and only if $\id-B$ is $\frac12$-cocoercive, and that the fixed points of $B$ are exactly the zeroes of $\id-B$. 
\smallbreak\noindent
If $A$ is $\beta$-cocoercive and $\tau\in(0,2\beta]$, then the explicit Euler step $\id-\tau A$ discussed above is nonexpansive, and its fixed points are the zeroes of $A$. Given a reference point $v\in\mathbb R^d$, parameters $\alpha_k\in(0,1)$ and an initial measure $\mu_0\in\mathcal P_2(\mathbb R^d)$, we consider the iteration
\[
  \mu_{k+1}:=\bigl(\alpha_kv+(1-\alpha_k)B_{\mu_k}\bigr)_\#\mu_k.
\]
This iteration has a natural interpretation in terms of a population of interacting agents. At each step, every agent adjusts its state according to $B_{\mu_k}$, which depends on the current distribution of the population, and has an additional displacement component directed towards the common reference point $v$ and proportional to $\alpha_k$. As in continuous time, it is the displacement towards the common reference point that brings the agents together. Since $B$ is nonexpansive, it does not increase the variance of the population, while the displacement towards $v$ multiplies it by $(1-\alpha_k)^2$ at the $k$-th step. Hence, as long as $\sum_{k=0}^{+\infty}\alpha_k=+\infty$, the variance vanishes and the agents eventually gather around a single point. 
\smallbreak\noindent
In Lemma~\ref{lemfixmoment} we prove that, as long as $B$ has at least one fixed point, it has a unique fixed point closest to $\delta_v$, and that this fixed point is a Dirac mass. In Theorem~\ref{thm:wass_halpern} we show that, under suitable assumptions on $(\alpha_k)$, the sequence $(\mu_k)$ converges in $\mathcal W_2$ to this fixed point. Thus, the population reaches consensus at a fixed point of $B$, the reference point $v$ plays the role of the origin in the regularized flow, and the limit is determined by $v$ rather than by the initial datum. Finally, in Corollary~\ref{cor:halpern_cocoercive} we apply this result to the explicit Euler step of a cocoercive operator and obtain that the iterates converge to the zero of the operator closest to $\delta_v$, which, in contrast with the limit of the flow in Theorem~\ref{thm1}, is always a Dirac mass.

%%%%%%%%%%%%%%%%%%%%%%%%%%%%%%%%%%%
%%%%%%%%%%%%%%%%%%%%%%%%%%%%%%%%%%% 

%%%%%%%%%%%%%%%%%%%%%%%%%%%%%%%%%%%
%%%%%%%%%%%%%%%%%%%%%%%%%%%%%%%%%%% 

\section{Overview of the main results}\label{s2}

\subsection{Operators on the Wasserstein space}
Let $\big(\mathcal P_2(\mathbb R^d),\mathcal W_2\big)$ denote the space of probability measures on $\mathbb R^d$ with finite second moment endowed with the Wasserstein distance. 

\medskip\noindent
We refer to a mapping $A:\mathcal P_2(\mathbb R^d)\to\bigsqcup_{\mu\in\mathcal P_2(\mathbb R^d)} L^2(\mu)^d$ as an \textit{operator on $\mathcal P_2(\mathbb R^d)$} if $A(\mu)\in L^2(\mu)^d$ for all $\mu\in \mathcal P_2(\mathbb R^d)$. Alternatively, with a little abuse of notation, we can see an operator on $\mathcal P_2(\mathbb R^d)$ as a family $A=(A_\mu)_{\mu\in\mathcal P_2(\mathbb R^d)}$ satisfying $A_\mu\in L^2(\mu)^d$ for every
$\mu\in\mathcal P_2(\mathbb R^d)$. 
\smallbreak\noindent
Let us begin by recalling the definition of monotonicity for such operators. 
\begin{definition}
    We say that an operator $A$ on $\mathcal P_2(\mathbb R^d)$ is monotone if
    \begin{align}\label{monotone_prop}
        \int_{\mathbb R^{2d}} \left\langle A_\mu(x)-A_\nu(y),x-y\right\rangle     \,d\gamma(x,y)\ge 0\,\quad \forall \gamma\in\Gamma_{o}(\mu,\nu),
    \end{align}
    for any $\mu,\nu\in\mathcal P_2(\mathbb R^d)$, 
    where $\Gamma_{o}(\mu,\nu)$ denotes the set of optimal couplings between $\mu$ and $\nu$.
\end{definition}\noindent
Several examples of monotone operators will be given in Section~\ref{sec:example}. Let us make a couple of remarks on the definition of monotonicity.
\begin{remark}
Condition \eqref{monotone_prop} is tested only along optimal couplings. 
Up to a change of sign, this corresponds to the metric dissipativity of a deterministic probability vector field in the sense of \cite[Definition~4.3 and Remark~4.4]{CSS_2023}, and requiring the same inequality along arbitrary couplings gives total dissipativity \cite{CSS_2026}.
The two notions coincide under more regularity, see Proposition \ref{p1} below.
The stronger, all-couplings form is the one used
in the mean field games literature, where it is known as \textit{displacement monotonicity};
it first appeared as weak monotonicity in \cite[condition~(A4)]{A_2016}, in dimension one,
and its multidimensional formulation was given in \cite[Definition~2.2]{GMMZ_2022}. A
corresponding condition for nonseparable Hamiltonians was introduced in
\cite[Definition~3.4]{GMMZ_2022} and later weakened in
\cite[Lemma~2.7 and Remark~2.8]{MM_2024}. Displacement monotonicity and the Lasry--Lions
condition \cite{LL_2007} are in general independent; see \cite[Remark~2.5]{GMMZ_2022},
\cite[Remark~2.1]{MM_2024} and \cite[Theorem~4.1]{GM_2023}.
\end{remark}

\begin{remark}
Condition \eqref{monotone_prop} may be viewed as a Wasserstein counterpart of
classical monotonicity in Hilbert spaces. If $A$ satisfies
\eqref{monotone_prop}, then testing it on Dirac masses shows that
$x\longmapsto A_{\delta_x}(x)$ is a monotone vector field on $\mathbb R^d$.
Conversely, every monotone vector field $a:\mathbb R^d\to\mathbb R^d$ satisfying $a \in L^2(\mu)^d$ for every $\mu \in \mathcal{P}_2(\R^d)$ induces
an operator satisfying \eqref{monotone_prop} by defining the operator $A$ on $\mathcal P_2(\mathbb R^d)$ as $A_\mu:=a$. 
\end{remark}\noindent
As in the classical theory, monotonicity  alone is of little use, unless it is combined with a regularity assumption, say 
a continuity property. We use the notion of demicontinuity of
\cite[Definition~3.22]{CSS_2026}, based on the strong--weak topology employed in
\cite{CSS_2023} and also used in the multispecies setting in
\cite[Definition~C.4]{CHMR_2025}.
\begin{definition}
    We say that an operator $A$ on $\mathcal P_2(\mathbb R^d)$ is demicontinuous at $\mu$ if, for any sequence $ (\mu_n)$ such that
    $\mu_n\xrightarrow{n\to+\infty}\mu$ in $\mathcal W_2$, there holds
    \[
        \sup_{n\in\mathbb N}\|A_{\mu_n}\|_{L^2(\mu_n)^d}<+\infty   \qquad\text{and}\qquad         (\mathrm{id},A_{\mu_n})_\#\mu_n \xrightharpoonup{n\to+\infty} (\mathrm{id},A_\mu)_\#\mu        \,\,\,\text{narrowly}.
    \]
    We say that $A$ is demicontinuous if it is demicontinuous at every
    $\mu\in\mathcal P_2(\mathbb R^d)$.
\end{definition}
\noindent
Demicontinuity plays the same role as hemicontinuity in the classical theory of monotone operators. The following proposition records some consequences when combined with monotonicity.
\begin{proposition}\label{p1}
    Let $A$ be a demicontinuous monotone operator on
    $\mathcal P_2(\mathbb R^d)$. The following statements hold.
    \begin{itemize}[label=\raisebox{0.2ex}{\scriptsize$\bullet$}]
        \item The monotonicity condition holds for all couplings, i.e.,
        \[
            \int_{\mathbb R^{2d}} \left\langle A_\mu(x)-A_\nu(y),x-y \right\rangle \,d\gamma(x,y)\geq 0\, \quad \forall \gamma\in\Gamma(\mu,\nu),
        \]
        for every $\mu,\nu\in\mathcal P_2(\mathbb R^d)$.

        \smallbreak
        \item The operator $A$ is maximally monotone, in the sense that, if $\nu\in\mathcal P_2(\mathbb R^d)$ and $\xi\in L^2(\nu)^d$ satisfy
        \[
            \int_{\mathbb R^{2d}} \left\langle A_\mu(x)-\xi(y),x-y \right\rangle \,d\gamma(x,y)\geq 0\,\quad \forall \gamma\in\Gamma(\mu,\nu),
        \]
        for every $\mu\in\mathcal P_2(\mathbb R^d)$, then $\xi(x)=A_\nu(x)$ for $\nu\text{-a.e. }x\in\mathbb R^d$.

        \smallbreak
        \item The set $\mathrm{zer}(A)$ is convex along the
        interpolations induced by arbitrary couplings, i.e.,
        \[
            \bigl((1-t)\mathrm{proj}_1+t\mathrm{proj}_2\bigr)_\#\gamma \in\mathrm{zer}(A)\,\quad  \forall \gamma\in\Gamma(\mu,\nu),
        \]
        for every $t\in[0,1]$ and every $\mu,\nu\in\mathrm{zer}(A)$, where       $\mathrm{proj}_1,\mathrm{proj}_2:\mathbb R^d\times\mathbb R^d\to\mathbb R^d$ denote the canonical projections. 
        Here, $\mathrm{zer}(A) := \left\{ \mu\in\mathcal P_2(\mathbb R^d) :\, A_\mu(x)= 0 \,\,\, \text{for }\mu\text{-a.e. $x\in\mathbb R^d$} \right\}$.
    \end{itemize}
\end{proposition}\noindent
We make two remarks on this proposition.

\begin{remark}
The first two statements follow from~\cite[Theorem~7.2]{CSS_2026}, where it is  shown that an everywhere-defined deterministic demicontinuous probability vector field that is dissipative
along optimal couplings is \textit{maximally totally dissipative}.
The result applies here after reversing the sign of the operator. 
\end{remark}

\begin{remark}
    The second and third assertions of Proposition~\ref{p1} are the transposition of familiar facts. 
The second one is the analogue of the theorem of Minty \cite[Theorem~4]{M_1962} asserting that a
monotone hemicontinuous operator defined on the whole of a Hilbert space is maximally monotone,
and its proof follows the corresponding Minty argument. 
The third one extends the classical result that the set of zeroes of a maximal monotone operator in a Hilbert space is convex \cite[Proposition~23.39]{BC_2017}.
\end{remark}

\subsection{Operator flows and long-time behavior}
\noindent
Given an operator $A$ on $\mathcal P_2(\mathbb R^d)$, we are interested in the evolution of probability measures under the velocity fields induced by $-A$. More precisely, we consider curves of measures whose mass is transported according to these velocity fields, without creation or loss of mass. Such evolutions are described by the associated continuity equation.
\smallbreak\noindent
Fix an initial time $t_0\in\mathbb R$. We begin by recalling what it means that a curve be driven by a monotone flow. 
\begin{definition}
    Let $A$ be an operator on $\mathcal P_2(\mathbb R^d)$ and let $\mu_0\in\mathcal P_2(\mathbb R^d)$.
    We say that a locally absolutely continuous curve $\mu:[t_0,+\infty)\to\mathcal P_2(\mathbb R^d)$ satisfying $\mu(t_0) = \mu_0$ is a solution of
    \begin{equation}\label{eqtflw}
        \begin{cases}
            \partial_t\mu_t -\mathrm{div}\bigl(A_{\mu_t}\mu_t\bigr)= 0\\[0.3em]
            \mu_{t_0}=\mu_0
        \end{cases}
    \end{equation}
    if  one has $\displaystyle\int_{t_0}^\tau\|A_{\mu_t}\|_{L^2(\mu_t)^d}\,dt<+\infty$ for any $\tau\in [t_0,+\infty)$ and 
    \[
        \dfrac{d}{dt}\int_{\mathbb R^d}\varphi(x)\,d\mu_t(x) = -\int_{\mathbb R^d}        \left\langle\nabla\varphi(x),A_{\mu_t}(x)\right\rangle
        \,d\mu_t(x)\quad \text{ for every } \varphi\in C_c^\infty(\mathbb R^d)
    \]
    in the sense of distributions on $(t_0,+\infty)$.
\end{definition}\noindent
The following lemma establishes the well-posedness of solutions under monotonicity and demicontinuity.
\begin{lemma}\label{l1}
Let $A$ be a demicontinuous monotone operator on $\mathcal P_2(\mathbb R^d)$. 
Then, for every $\mu_0\in\mathcal P_2(\mathbb R^d)$, there exists a unique solution $\mu:[t_0,+\infty)\to\mathcal P_2(\mathbb R^d)$ of~\eqref{eqtflw}.
\end{lemma}\noindent
The preceding lemma is the Wasserstein analogue of the well-posedness and contraction theory for evolution equations governed by maximal monotone operators in Hilbert spaces \cite[Chapter~3]{B_1973}. Its conclusions can also be deduced from \cite[Theorems~4.2, 4.4 and~7.2]{CSS_2026}, building on the theory developed in \cite[Section~5]{CSS_2023}. We give an alternative proof based on Theorem~\ref{general_existence}.
\smallbreak\noindent
We now turn to the long-time behavior of the flow. As in the classical theory, the notion of \textit{cocoercivity} plays a central role in the convergence of the trajectory to a zero of the operator. Let us introduce the corresponding notion for the space of probability measures.
\begin{definition}
    Let $\beta>0$. We say that an operator $A$ on $\mathcal P_2(\mathbb R^d)$ is $\beta$-cocoercive if 
    \begin{equation}\label{cocoercive_prop}
        \begin{aligned}
            &\int_{\mathbb R^{2d}}
            \left\langle A_\mu(x)-A_\nu(y),x-y \right\rangle \,d\gamma(x,y)\geq \beta             \int_{\mathbb R^{2d}} \left|A_\mu(x)-A_\nu(y)\right|^2\,d\gamma(x,y)\,\quad \forall \gamma\in\Gamma_o(\mu,\nu),
        \end{aligned}
    \end{equation}
    for every $\mu,\nu \in \mathcal P_2(\mathbb R^d)$.
    %We say that $A$ is cocoercive if it is $\beta$-cocoercive for some $\beta>0$.
\end{definition}\noindent
Examples of cocoercive operators will be given in Section~\ref{sec:example1}. Notice that cocoercivity is stronger than monotonicity and also provides the continuity needed for the well-posedness of the associated flow.
 \smallbreak\noindent
The classical convergence result for cocoercive operators carries over to the Wasserstein setting as follows.
\begin{theorem}\label{thm1}
Let $\mu_0\in\mathcal P_2(\mathbb R^d)$ and $\beta>0$ be given. Suppose that $A$ is a $\beta$-cocoercive operator satisfying $\mathrm{zer}(A)\neq\varnothing$. Let $\mu:[t_0,+\infty)\to\mathcal P_2(\mathbb R^d)$ be the unique solution of $(\ref{eqtflw})$. 
The following statements hold.
\begin{itemize}[label=\raisebox{0.2ex}{\scriptsize$\bullet$}]
    \item For every $t \in [t_0,+\infty)$ and  every $\nu\in\mathrm{zer}(A)$,
    \begin{equation} \label{eq.cocoercive-energy}
\frac12\mathcal W_2^2(\mu_t,\nu) + \beta\int_{t_0}^t \|A_{\mu_s}\|_{L^2(\mu_s)^d}^2\,ds \leq \frac12\mathcal W_2^2(\mu_0,\nu).
\end{equation}
 \smallbreak
    \item  $\| A_{\mu_t} \|_{L^2(\mu_t)^d} = o(t^{-1/2})$ \,and\,\,   $|\mu'|(t) = o(t^{-1/2})$\, as\,  $t \longrightarrow+\infty$.
\end{itemize}
Moreover, there exists $\mu^\circ\in\mathrm{zer}(A)$ such that $\mathcal W_2(\mu_t,\mu^\circ)\longrightarrow0$\,  as\,  $t\longrightarrow+\infty$. 
\end{theorem}
We make now some remarks. 

\begin{remark}\label{non-dirac-mass}
The limit $\mu^\circ\in\mathrm{zer}(A)$ in Theorem~\ref{thm1} need not be a Dirac mass. For instance, consider the operator $A$ on $\mathcal P_2(\mathbb R^d)$ given by $A_\mu:=m_\mu$, where $m_\mu\in\mathbb R^d$ denotes the barycenter of $\mu$. This operator is $1$-cocoercive, and the corresponding flow converges to $(\id-m_{\mu_0})_\#\mu_0$. In particular, the limit is a Dirac mass if and only if the initial measure $\mu_0$ is a Dirac mass.
\end{remark}

\begin{remark}
For the flow of a cocoercive operator, one can check that
both mappings $t\longmapsto\|A_{\mu_t}\|_{L^2(\mu_t)^d}$
and $t\longmapsto|\mu'|(t)$ are continuous.
Thus, although the metric derivative of a locally
absolutely continuous curve is in general defined only
almost everywhere, the decay estimates in
Theorem~\ref{thm1} hold pointwise as $t\longrightarrow+\infty$.
\end{remark}

\begin{remark}
Theorem~\ref{thm1} is the Wasserstein counterpart of Bruck's weak convergence theorem \cite[Theorem~1]{B_1975}, with Opial's lemma \cite[Lemma~1]{O_1967} replaced by its Wasserstein analogue \cite[Corollary~5.3]{NS_2021}. In the Wasserstein setting, however, convergence is strong (in $\mathcal W_2$); this is due to the relative compactness of the trajectories. In Hilbert spaces, the corresponding residual decay estimates are given in \cite[Theorem~A.1]{BN_2026}; related convergence results for monotone flows of measures can be found in \cite[Theorem~3.2 and Proposition~3.3]{CHMR_2025}.
\end{remark}

\subsection{Tikhonov regularization}
The convergence and decay results of the previous subsection rely on cocoercivity. For merely demicontinuous monotone operators, however, the flow need not converge to a zero, and the residual may fail to vanish. A classical approach to recovering asymptotic stationarity is Tikhonov regularization, which consists in adding to the dynamics a time-dependent vanishing multiple of the identity.

\medskip\noindent
We fix an initial time $t_0\in\mathbb R$ and a nonincreasing, locally integrable function $\varepsilon:[t_0,+\infty)\to(0,+\infty)$, which we refer to as the Tikhonov regularization parameter throughout this subsection.
\begin{definition}
    Let $A$ be an operator on $\mathcal P_2(\mathbb R^d)$ and let
    $\mu_0\in\mathcal P_2(\mathbb R^d)$. We say that a locally absolutely
    continuous curve
    $\mu:[t_0,+\infty)\to\mathcal P_2(\mathbb R^d)$ satisfying $\mu(t_0)=\mu_0$ is a solution of
    \begin{equation}\label{eqtfl}
        \begin{cases}
            \partial_t\mu_t
            -\mathrm{div}\bigl((A_{\mu_t}+\varepsilon(t)\id)\mu_t\bigr)= 0\\[0.3em]
            \mu_{t_0}=\mu_0
        \end{cases}
    \end{equation}
    if  one has $\displaystyle\int_{t_0}^\tau\|A_{\mu_t}\|_{L^2(\mu_t)^d}\,dt<+\infty$ for any $\tau\in [t_0,+\infty)$ and 
    \[
        \dfrac{d}{dt}\int_{\mathbb R^d}\varphi(x)\,d\mu_t(x)
        =  -\int_{\mathbb R^d} \left\langle \nabla\varphi(x), A_{\mu_t}(x)+\varepsilon(t)x      \right\rangle \,d\mu_t(x) \quad \text{ for every } \varphi\in C_c^\infty(\mathbb R^d)
    \]
    in the sense of distributions on $(t_0,+\infty)$.
\end{definition}\noindent
The Tikhonov term can be thought of as making the operator strongly monotone at each time. The following lemma establishes well-posedness of the regularized flow. 
\begin{lemma}\label{l2}
    Suppose that 
    $\varepsilon:[t_0,+\infty)\to(0,+\infty)$ is Borel measurable and locally integrable. 
    Let $A$ be a demicontinuous monotone operator on $\mathcal P_2(\mathbb R^d)$. 
    Then, for every $\mu_0\in\mathcal P_2(\mathbb R^d)$, there exists a unique solution $\mu:[t_0,+\infty)\to\mathcal P_2(\mathbb R^d)$ of \eqref{eqtfl}.
\end{lemma}
With well-posedness established, we now turn to the long-time behavior of the regularized flow, but before we give another relevant lemma. 
\begin{lemma}\label{lemsecmoment}
Let $A$ be a demicontinuous monotone operator on $\mathcal P_2(\mathbb R^d)$ such that
$\mathrm{zer}(A)\neq\varnothing$.
Then the set $K:=\{x\in\mathbb R^d:\,\delta_x\in\mathrm{zer}(A)\}$
is nonempty, closed and convex. 
Moreover,
\[
    \big\{ \delta_{\mathrm{proj}_K(0)} \big\} = \displaystyle\argmin_{\nu\in\mathrm{zer}(A)}\mathcal M_2(\nu),
\]
where $\mathcal M_2(\nu)$ denotes the second moment of a measure $\nu\in\mathcal P_2(\mathbb R^d)$.
\end{lemma}
\iffalse 
Given $\mu\in\mathcal P_2(\mathbb R^d)$, we denote by . 
\smallbreak\noindent
\fi 
The following theorem provides decay estimates for the metric speed and the operator residual under suitable assumptions on the Tikhonov parameter.
\begin{theorem}\label{thm2}
    Let $A$ be a demicontinuous monotone operator on $\mathcal P_2(\mathbb R^d)$ such that    $\mathrm{zer}(A)\neq\varnothing$. 
    Suppose that $\varepsilon:[t_0,+\infty)\to(0,+\infty)$ is locally absolutely continuous and satisfies
    \[
        \lim_{t\longrightarrow+\infty}\varepsilon(t)= 0, \qquad       \int_{t_0}^{+\infty}\varepsilon(t)\,dt=+\infty, \qquad         \int_{t_0}^{+\infty}|\dot\varepsilon(t)|\,dt<+\infty,
    \]
    and define $\gamma_\varepsilon:[t_0,+\infty)\to (0,+\infty)$ by $\gamma_\varepsilon(t) :=         \exp\left( \mathlarger{\int}_{t_0}^t\varepsilon(s)\,ds\right).$ 
    Let $\mu_0\in\mathcal P_2(\mathbb R^d)$ and let $\mu:[t_0,+\infty)\to\mathcal P_2(\mathbb R^d)$ be the unique solution of  \eqref{eqtfl}.
    Then
    \[
        |\mu'|(t)\longrightarrow 0
        \quad\text{and}\qquad \|A_{\mu_t}\|_{L^2(\mu_t)^d}\longrightarrow 0 \qquad\text{as}\quad t\longrightarrow +\infty.
    \]
% More precisely, there holds
%    \[
%        |\mu'|(t) =  O\left(\frac{1}{\gamma_\varepsilon(t)}\right) +  O\left(
%        \frac{1}{\gamma_\varepsilon(t)} \int_{t_0}^t       |\dot\varepsilon(\tau)|\gamma_\varepsilon(\tau)\,d\tau
%        \right),
%    \]
%    \[
%        \|A_{\mu_t}\|_{L^2(\mu_t)^d} = O\left(\frac{1}{\gamma_\varepsilon(t)}\right) +       O\left( \frac{1}{\gamma_\varepsilon(t)} \int_{t_0}^t        |\dot\varepsilon(\tau)|\gamma_\varepsilon(\tau)\,d\tau   \right)  +       O(\varepsilon(t)).
%    \]
More precisely, there holds
\[
    |\mu'|(t) = O(\min \{ R_1(t), R_2(t) \}) \quad \text{ and } \quad \| A_{\mu_t} \|_{L^2(\mu_t)^d} = O \big( \min \{ R_1(t), R_2(t) \} + \varepsilon(t) \big)
\]
as $t\longrightarrow+\infty$, where $R_1,R_2:[t_0,+\infty)\to \mathbb R$ are given by 
\begin{align*}
    R_1(t) := & \, \dfrac{1}{\gamma_\varepsilon(t) } \left( 1 + \int_{t_0}^t |\dot{\varepsilon}(\tau)| \gamma_\varepsilon(\tau) \, d\tau \right) \quad\text{and}\quad
    R_2(t) :=  \,\dfrac{1}{\gamma_\varepsilon(t) } \left( 1 + \int_{t_0}^t |\dot{\varepsilon}(\tau) + \varepsilon^2(\tau)| \gamma_\varepsilon(\tau) \, d\tau \right),
\end{align*}
respectively. Moreover, $\mathcal W_2\big(\mu_t, \delta_{\mathrm{proj}_K(0)}\big)\longrightarrow 0$ as $t\longrightarrow+\infty$, where $K\subseteq\mathbb R^d$ is as in Lemma \ref{lemsecmoment}.
\end{theorem}\noindent

\noindent
Since the preceding theorem applies to a broad class of Tikhonov regularization parameters, some specializations are in order. The following corollary provides explicit decay rates for several particular choices.
\begin{corollary}\label{mccr}
    Let $A$ be a demicontinuous monotone operator on $\mathcal P_2(\R^d)$ such that $\mathrm{zer}(A)\neq\varnothing$.
    Let $\mu_0\in\mathcal P_2(\R^d)$, and $\varepsilon:[t_0,+\infty)\to(0,+\infty)$ be given, and let $\mu:[t_0,+\infty)\to\mathcal P_2(\R^d)$ be the unique solution of \eqref{eqtfl} with initial datum $\mu_0$. 
    Then, the following statements hold.
    \begin{itemize}[label=\raisebox{0.2ex}{\scriptsize$\bullet$}]
        \item If $t_0>0$ and $\varepsilon(t)=\alpha/t$ for some $\alpha>0$, then
        \[
            |\mu'|(t)=O\bigl(t^{-\min\{\alpha,1\}}\bigr)
            \qquad\text{and}\qquad
            \|A_{\mu_t}\|_{L^2(\mu_t)^d}=O\bigl(t^{-\min\{\alpha,1\}}\bigr) \quad \text{as $t\longrightarrow+\infty$.}
        \]
        \smallbreak
        \item If $t_0>0$ and $\varepsilon(t)=\alpha/t^q$ for some $\alpha>0$ and $q\in(0,1)$, then
        \[
            |\mu'|(t)=O\left(\frac1t\right)
            \qquad\text{and}\qquad
            \|A_{\mu_t}\|_{L^2(\mu_t)^d}=O\left(\frac1{t^q}\right) \quad \text{as $t\longrightarrow+\infty$.}
        \]
         \smallbreak
        \item If $t_0>1$ and $\varepsilon(t)=1/(t\log t)$, then
        \[
            |\mu'|(t)=O\left(\frac1{\log t}\right)
            \qquad\text{and}\qquad
            \|A_{\mu_t}\|_{L^2(\mu_t)^d}=O\left(\frac1{\log t}\right) \quad \text{as $t\longrightarrow+\infty$.}
        \]
        
    \end{itemize}
\end{corollary}
\smallbreak\noindent
We close the subsection with a couple of remarks.

\begin{remark}
The decay estimates for the metric derivative in Theorem~\ref{thm2}
are understood along sequences of times at which it exists.
This means, precisely, that there exists a set
$E\subseteq[t_0,+\infty)$ of full Lebesgue measure
such that $|\mu'|(t)$ exists for every $t\in E$ and,
for every sequence $(t_k)_{k\in\mathbb N}\subseteq E$
with $t_k\longrightarrow+\infty$, one has
$|\mu'|(t_k)=O(\min\{R_1(t_k),R_2(t_k)\})$
as $k\longrightarrow+\infty$.
\end{remark}

\begin{remark}
Theorem~\ref{thm2} and Corollary~\ref{mccr} extend to the Wasserstein space the decay estimates for Tikhonov regularized monotone operator flows obtained in \cite[Proposition~2.5, Theorem~2.8, and Corollaries~2.9--2.11]{BN_2026} for Hilbert spaces. In particular, the rate $O(1/t)$ for the metric speed and the operator residual is preserved when $\varepsilon(t)=\alpha/t$ with $\alpha\geq 1$.
\end{remark}

\subsection{A Halpern--type algorithm}
We now turn to a discrete-time scheme. The so-called Halpern iteration for nonexpansive mappings is a discrete analogue of the Tikhonov regularized flow. It is obtained by an explicit discretization of the flow in which the regularization term is centered at a fixed reference point instead of the origin; see \cite[Section~4.2]{BN_2026}. 
\smallbreak\noindent
We introduce the following notion of nonexpansiveness on the Wasserstein space.
\begin{definition}\label{def:wass_nonexpansive}
Let $B$ be an operator on $\mathcal P_2(\mathbb R^d)$. We say that $B$ is nonexpansive if
\begin{equation}\label{eq:wass_nonexpansive}
\int_{\mathbb R^{2d}}|B_\mu(x)-B_\nu(y)|^2\,d\gamma(x,y)\leq\int_{\mathbb R^{2d}}|x-y|^2\,d\gamma(x,y) \quad \forall \gamma \in \Gamma_o(\mu, \nu),
\end{equation}
for every $\mu, \nu \in \mathcal{P}_2(\R^d)$.
\end{definition}\noindent
Examples of nonexpansive operators will be given in Section~\ref{sec:example2}. In analogy with the set of zeroes of an operator, we denote the set of fixed points of an operator $B$ on $\mathcal P_2(\mathbb R^d)$ by 
\[
    \mathrm{fix}(B):=\left\{\mu\in\mathcal P_2(\mathbb R^d):\, B_\mu(x)=x  \,\,\, \text{for }\mu\text{-a.e. $x\in\mathbb R^d$}\right\}.
\]
Notice that $\mu\in\mathrm{fix}(B)$ is a stronger requirement than the invariance $(B_\mu)_\#\mu=\mu$, since it asks that $\mu$-almost every point be left fixed by $B_\mu$ and not only that the mass be rearranged onto itself. As in Hilbert spaces, nonexpansive and cocoercive operators are related through the passage from $B$ to $\id-B$, and fixed points of the former correspond to zeroes of the latter. The following proposition makes this precise.
\begin{proposition}\label{prop:wass_nonexpansive}
Let $B$ be an operator on $\mathcal P_2(\mathbb R^d)$, and consider the operator $A$ on $\mathcal P_2(\mathbb R^d)$ given by $A_\mu(x):=x-B_\mu(x)$. Then, $B$ is nonexpansive if and only if $A$ is $\frac{1}{2}$-cocoercive. 
    Moreover, $\mathrm{fix}(B) = \mathrm{zer}(A)$.
\end{proposition}\noindent
We remark a couple of immediate  consequences of the previous proposition. 
\begin{remark}\label{rmk.nonexp}
    It follows easily from Proposition~\ref{prop:wass_nonexpansive} that if $B$ is nonexpansive, then the mapping $\mathcal B:\mathcal P_2(\mathbb R^d)\to\mathcal P_2(\mathbb R^d)$ defined by $\mathcal B(\mu):=(B_\mu)_\#\mu$ is nonexpansive on the metric space $(\mathcal P_2(\mathbb R^d),\mathcal W_2)$. However, simple examples show that the converse does not hold. 
\end{remark}

\begin{remark}
It follows from Proposition~\ref{prop:wass_nonexpansive} that if $B$ is nonexpansive, then $\delta_{m_\mu}\in\mathrm{fix}(B)$ for every $\mu\in\mathrm{fix}(B)$, where $m_\mu$ denotes the barycenter of $\mu$.
\end{remark}

\noindent
Before turning to the convergence of the discrete-time scheme, we give a lemma describing some properties of the fixed points of nonexpansive operators. 
\begin{lemma}\label{lemfixmoment}
Let $B$ be a nonexpansive operator on
$\mathcal P_2(\mathbb R^d)$ such that
$\mathrm{fix}(B)\neq\varnothing$.
Then the set
$Q:=\{x\in\mathbb R^d:\,\delta_x\in\mathrm{fix}(B)\}$
is nonempty, closed and convex. Moreover,
\[
    \delta_{\mathrm{proj}_Q(v)} = \mathrm{proj}_{\mathrm{fix}(B)} (\delta_v)
\]
for any $v\in\mathbb R^d$.
\end{lemma}
We are now ready to state our convergence result for the aforementioned scheme.
\begin{theorem}\label{thm:wass_halpern}
Let $B$ be a nonexpansive operator on $\mathcal P_2(\mathbb R^d)$ satisfying $\mathrm{fix}(B)\neq\varnothing$.
Let $v\in\mathbb R^d$, $\mu_0\in\mathcal P_2(\mathbb R^d)$ and  $(\alpha_k)_{k\in\mathbb N\cup\{0\}}\subseteq(0,1)$ satisfy
\begin{equation}\label{eq:halpern_parameters}
\sum_{k = 0}^{+\infty} \alpha_k = +\infty\qquad\text{and}\qquad\alpha_k \xrightarrow{k \longrightarrow +\infty} 0.
\end{equation}
Assume, in addition, that either
\begin{equation}\label{eq:halpern_wittmann}
\sum_{k = 0}^{+\infty}|\alpha_{k+1}-\alpha_k| < +\infty\qquad\text{or}\qquad \frac{|\alpha_{k+1}-\alpha_k|}{\alpha_{k+1}}\xrightarrow{k \longrightarrow +\infty}0.
\end{equation}
Starting from $\mu_0$, define recursively
\begin{equation}\label{eq:wass_halpern}
\mu_{k+1} := \bigl(\alpha_kv+(1-\alpha_k)B_{\mu_k}\bigr)_\#\mu_k \quad \text{for }k\in\mathbb N\cup\{0\}.
\end{equation}
Then, $\mathcal W_2\big(\mu_k,\delta_{\mathrm{proj}_Q(v)}\big)\xrightarrow{k \longrightarrow +\infty}0$, where $Q$ is the set
defined in Lemma~\ref{lemfixmoment}.
\end{theorem}
Before passing to some of the consequences of this result, let us make some remarks.

\begin{remark}
The condition $\alpha_k\xrightarrow{k \longrightarrow +\infty}0$ in Theorem~\ref{thm:wass_halpern}
cannot in general be omitted. If
$\alpha_k=\alpha\in(0,1)$ for $k\in\mathbb N\cup\{0\}$, then the mapping
$\mathcal H_\alpha:\mathcal P_2(\mathbb R^d)
\to\mathcal P_2(\mathbb R^d)$ defined by
\[
\mathcal H_\alpha(\mu)
:=(\alpha v+(1-\alpha)B_\mu)_\#\mu
\]
is a contraction with constant $1-\alpha$.
Consequently, the iteration converges geometrically to
a unique Dirac mass, which need not belong to
$\mathrm{fix}(B)$. %For instance, if $B\equiv 0$ and $v\neq 0$, then the iteration converges to $\delta_{\alpha v}$, whereas $\mathrm{fix}(B)=\{\delta_0\}$.
\end{remark}

\begin{remark}
Let $B$ be a nonexpansive operator and define $\mathcal B:\mathcal P_2(\mathbb R^d)\to \mathcal P_2(\mathbb R^d)$ by 
$\mathcal B(\mu):=(B_\mu)_\#\mu$. Then
\[
Q=\{m_\mu:\mu\in\mathrm{fix}(B)\}
=\{m_\mu:\mathcal B(\mu)=\mu\},
\]
where $Q$ is the set introduced in Lemma~\ref{lemfixmoment}.
In particular, every invariant measure of $\mathcal B$ produces a fixed point of $B$, i.e., $\delta_{m_\mu}\in\mathrm{fix}(B)$ whenever $\mathcal B(\mu)=\mu$.
\end{remark}

\begin{remark}\label{rem:halpern_rates}
The estimates obtained in the proof of Theorem~\ref{thm:wass_halpern}, together with the arguments of \cite[Theorem~4.6]{BN_2026}, can yield quantitative convergence rates for the choice $\alpha_k=(\alpha-1)/(k+\alpha)$ with $\alpha>1$. More precisely, one has
\[
    \mathrm{Var}(\mu_k)^{1/2}=O(k^{-(\alpha-1)})\qquad \text{and}\qquad \mathcal W_2(\mu_k,\mu_{k-1}) + \mathcal W_2(\mu_k,(B_{\mu_k})_\#\mu_k) = O(k^{-\min\{\alpha-1,1\}})
\]
as $k\longrightarrow+\infty$, where $\mathrm{Var}(\nu)$ denotes the variance of a measure $\nu\in\mathcal P_2(\mathbb R^d)$. 
\end{remark}

\begin{remark}
Averaged operators arise naturally in splitting schemes in Hilbert spaces. In the Wasserstein space, we may define an  operator $B$ on $\mathcal P_2(\mathbb R^d)$  to be $\theta$-averaged, for some $\theta\in(0,1)$, if $B=(1-\theta)\id+\theta N$ for some nonexpansive operator $N$. It is not difficult to see that averaged operators are nonexpansive, so Theorem~\ref{thm:wass_halpern} applies to this class of operators whenever $\mathrm{fix}(B)\neq\varnothing$.
\end{remark}

\noindent
We close the subsection with a corollary on a Halpern--type algorithm for finding zeros of cocoercive operators.

\begin{corollary}\label{cor:halpern_cocoercive}
Let $\mu_0\in\mathcal P_2(\mathbb R^d)$ and $\beta>0$ be given. Suppose that $A$ is a $\beta$-cocoercive operator on $\mathcal P_2(\mathbb R^d)$ satisfying $\mathrm{zer}(A)\neq\varnothing$, and let $\tau\in(0,2\beta]$, $v\in\mathbb R^d$ and $(\alpha_k)_{k\in\mathbb N\cup\{0\}}\subseteq(0,1)$ satisfy \eqref{eq:halpern_parameters}.
Assume, in addition, that either of the conditions in \eqref{eq:halpern_wittmann} holds.
Starting from $\mu_0$, define recursively
\begin{equation}\label{eq:halpern_cocoercive}
\mu_{k+1}:=\bigl(\alpha_kv+(1-\alpha_k)(\id-\tau A_{\mu_k})\bigr)_\#\mu_k
\quad\text{for }k\in\mathbb N\cup\{0\}.
\end{equation}
Then, $\mathcal W_2\big(\mu_k,\delta_{\mathrm{proj}_K(v)}\big)\xrightarrow{k\longrightarrow+\infty}0$, where $K$ is the set defined in Lemma~\ref{lemsecmoment}.
\end{corollary}

\iffalse 
\noindent
Let us finally record a quantitative consequence for Theorem~\ref{thm:wass_halpern}.
For the classical Halpern iteration in Hilbert spaces, Bo{\c t} and Nguyen showed in \cite[Theorem~4.6]{BN_2026} that, for $\alpha_k = (\alpha-1)/(k+\alpha)$ with $\alpha>1$, both the successive difference and the fixed-point residual decay as $O(k^{-\min\{\alpha-1,1\}})$.

For $\mu\in\mathcal P_2(\mathbb R^d)$, we respectively denote its barycenter and variance by 
\[ 
m_\mu:=\int_{\mathbb R^d}x\,d\mu(x) \quad \text{ and } \quad  \mathrm{Var}(\mu):=\int_{\mathbb R^d}|x-m_\mu|^2\,d\mu(x).
\]

\begin{corollary}\label{cor:halpern_rates}
In the framework of Theorem~\ref{thm:wass_halpern}, let $\alpha_k := (\alpha-1)/(k+\alpha)$, for every $k \in \mathbb N$, for some $\alpha > 1$.
Then, one has, as $k \longrightarrow +\infty$,
\begin{align*}
    (i) & \qquad \mathrm{Var}(\mu_k)^{1/2} = O \big(k^{-(\alpha-1)}\big);\\
    (ii) & \qquad \mathcal W_2(\mu_k,\mu_{k-1}) = O\big(k^{-\min\{\alpha-1,1\}} \big);\\
    (iii) & \qquad \mathcal W_2(\mu_k,(B_{\mu_k})_\#\mu_k) = O\big(k^{-\min\{\alpha-1,1\}}\big).
\end{align*}
\end{corollary}
\fi 

\section{Notation and preliminaries}

Throughout the paper, $d\in\mathbb N$ is a fixed dimension, $\langle\cdot,\cdot\rangle$ and $|\cdot|$ denote the Euclidean scalar product and norm on $\mathbb R^d$, respectively, and $t_0\in\mathbb R$ is a fixed initial time. 
We write $C_b(\mathbb R^m)$ for the space of bounded continuous functions on $\mathbb R^m$, $C_c^\infty(\mathbb R^m)$ for the space of smooth functions with compact support and $\id$ for the identity mapping of $\mathbb R^d$. 
Given a Banach space $E$, a mapping $u:I\to E$ defined on an interval $I\subseteq\mathbb R$ is said to be locally absolutely continuous if it is absolutely continuous on every compact subinterval of $I$ and we denote by $\mathrm{AC}_{\mathrm{loc}}(I;E)$ the set of such mappings.

\subsection{Measures, push-forwards and narrow convergence}
Denote by $\mathcal P(\mathbb R^m)$ the set of Borel probability measures on $\mathbb R^m$.
Given $\mu\in\mathcal P(\mathbb R^m)$ and a Borel mapping $T:\mathbb R^m\to\mathbb R^k$, the push-forward of $\mu$ through $T$ is the measure $T_\#\mu\in\mathcal P(\mathbb R^k)$ defined by $T_\#\mu(B):=\mu(T^{-1}(B))$ for every Borel set $B\subseteq\mathbb R^k$; equivalently,
\[
    \int\psi\,d(T_\#\mu)=\int\psi\circ T\,d\mu \quad \text{ for every bounded Borel function $\psi$}.
\]
The same notation is used when the underlying spaces are Polish rather than Euclidean, as is the case for the evaluation mappings of Section~\ref{sec:conti_eq}.

\medskip\noindent
A sequence $(\lambda_n)_{n\in\mathbb N}\subseteq\mathcal P(\mathbb R^m)$ converges narrowly to $\lambda\in\mathcal P(\mathbb R^m)$, written $\lambda_n\rightharpoonup\lambda$, if $\mathlarger{\int} \psi\,d\lambda_n\to \mathlarger{\int} \psi\,d\lambda$ for every $\psi\in C_b(\mathbb R^m)$. 
The narrow topology on $\mathcal P(\mathbb R^m)$ is Polish and by Prokhorov's theorem a subset of $\mathcal P(\mathbb R^m)$ is relatively compact for it if and only if it is tight.
In particular, by Markov's inequality, any family of measures with uniformly bounded second moments is tight and hence relatively compact in the narrow topology. 
We shall also use repeatedly the existence of a countable family $(\psi_j)_{j\in\mathbb N}\subseteq C_b(\mathbb R^m)$ determining probability measures, i.e., such that two measures $\lambda,\sigma\in\mathcal P(\mathbb R^m)$ coincide as soon as  $\mathlarger{\int} \psi_j\,d\lambda = \mathlarger{\int} \psi_j\,d\sigma$ for every $j\in\mathbb N$; such a family exists because the narrow topology is separable and metrizable.

\subsection{Optimal transport and the Wasserstein distance}
Given a measure $\mu \in \mathcal{P}(\mathbb{R}^d)$, we define its second moment by $\mathcal M_2(\mu) := \mathlarger{\int}_{\mathbb{R}^d} |x|^2 \mathrm{d}\mu(x)$ and let $\mathcal P_2(\mathbb{R}^d)$ be the set of all Borel probability measures whose second moment is finite. 
Given $\mu,\nu\in\mathcal P_2(\mathbb R^d)$, we denote by $\Gamma(\mu,\nu)$ the set of couplings, or transport plans, between $\mu$ and $\nu$, i.e., the set of measures $\gamma\in\mathcal P(\mathbb R^{2d})$ whose marginals are $(\mathrm{proj}_1)_\#\gamma=\mu$ and $(\mathrm{proj}_2)_\#\gamma=\nu$, where $\mathrm{proj}_1,\mathrm{proj}_2:\mathbb R^d\times\mathbb R^d\to\mathbb R^d$ are the canonical projections.
The $2$-Wasserstein distance between two measures $\mu,\nu\in\mathcal P_2(\mathbb R^d)$ is given by
\begin{equation*}
\mathcal W_2(\mu, \nu) := \bigg( \inf \bigg\{ \int_{\mathbb{R}^{2d}} |x - y|^2 \mathrm{d}\gamma(x, y) ~\, : ~\, \gamma \in \Gamma(\mu, \nu) \bigg\} \bigg)^{1/2}, 
\end{equation*}
and in the sequel, we denote by 
\begin{equation*}
\Gamma_o(\mu,\nu)
:=
\argmin_{\gamma \in \Gamma(\mu,\nu)}
\int_{\mathbb R^{2d}} |x-y|^2 \,\mathrm d\gamma(x,y)
\end{equation*}
the set of so-called \textit{optimal transport plans}. 
The set $\Gamma(\mu,\nu)$ is nonempty, since it contains the product measure $\mu\otimes\nu$ and it is narrowly compact.
Since the cost is lower semicontinuous, the direct method of the calculus of variations shows that $\Gamma_o(\mu,\nu)$ is nonempty and narrowly compact.
It is a basic fact in optimal transport theory that $(\mathcal P_2(\mathbb{R}^d), \mathcal W_2(\cdot,\cdot))$ is a complete separable metric space and that $\mathcal W_2$ is lower semicontinuous with respect to narrow convergence in both of its arguments; moreover, $\mathcal W_2(\mu_n,\mu)\to0$ if and only if $\mu_n\rightharpoonup\mu$ narrowly and $\mathcal M_2(\mu_n)\to\mathcal M_2(\mu)$. 
We refer to \cite{AGS_2008} for these facts.

\medskip\noindent
For every $\mu\in\mathcal P_2(\mathbb R^d)$ and every
$\xi\in L^2(\mu)^d$, we set
\[
\|\xi\|_{L^2(\mu)^d}^2
:=
\int_{\mathbb R^d}|\xi(x)|^2\,d\mu(x),
\]
so that $L^2(\mu)^d$ is a separable Hilbert space with the associated scalar product, which we denote by $\langle\cdot,\cdot\rangle_{L^2(\mu)^d}$. 
The same notation is used for the space $L^2(\eta)^d$ attached to a probability measure $\eta$ on a Polish space. 
Two estimates will be used constantly. First, since every coupling is admissible in the definition of $\mathcal W_2$,
\begin{equation}\label{plan_estimate}
\mathcal W_2^2(\mu,\nu)\leq\int_{\mathbb R^{2d}}|x-y|^2\,d\gamma(x,y) \quad \text{ for every } \gamma\in\Gamma(\mu,\nu).
\end{equation}
Second, if $\eta$ is a probability measure on a Polish space and $X,Y\in L^2(\eta)^d$, then $(X,Y)_\#\eta\in\Gamma(X_\#\eta,Y_\#\eta)$, so that \eqref{plan_estimate} gives
\begin{equation}\label{lift_estimate} 
\mathcal W_2(X_\#\eta,Y_\#\eta)\leq\|X-Y\|_{L^2(\eta)^d}.
\end{equation}
Inequality \eqref{lift_estimate} is what allows us to transfer estimates obtained in a Hilbert
space to the Wasserstein space and it is the mechanism underlying most of the proofs of
Section~\ref{s2}.

\subsection{Barycenters and variance}
For $\mu\in\mathcal P_2(\mathbb R^d)$, we respectively denote its barycenter and variance by
\[
m_\mu:=\int_{\mathbb R^d}x\,d\mu(x) \quad \text{ and } \quad \mathrm{Var}(\mu):=\int_{\mathbb R^d}|x-m_\mu|^2\,d\mu(x),
\]
both of which are finite because $\mu$ has finite second moment. 
Expanding the square and using the definition of $m_\mu$, we obtain
\begin{equation}\label{variance_identity}
\int_{\mathbb R^d}|x-v|^2\,d\mu(x)=\mathrm{Var}(\mu)+|m_\mu-v|^2 \quad \text{ for every } v\in\mathbb R^d;
\end{equation}
in particular, $\mathrm{Var}(\mu)=\mathcal M_2(\mu)-|m_\mu|^2$ and $m_\mu$ is the unique minimizer of $v\longmapsto \mathlarger{\int} |x-v|^2\,d\mu(x)$. 
Since $\mu\otimes\delta_v$ is the only coupling between $\mu$ and $\delta_v$, identity \eqref{variance_identity} also reads
\begin{equation}\label{dirac_distance}
\mathcal W_2^2(\mu,\delta_v)=\mathrm{Var}(\mu)+|m_\mu-v|^2 \quad \text{ for every } v\in\mathbb R^d.
\end{equation}
Hence $\mathcal W_2(\delta_x,\delta_y)=|x-y|$ for all $x,y\in\mathbb R^d$, the variance $\mathrm{Var}(\mu)=\mathcal W_2^2(\mu,\delta_{m_\mu})$ is the squared distance from $\mu$ to the set of Dirac masses and $\mathrm{Var}(\mu)= 0$ if and only if $\mu=\delta_{m_\mu}$.

\medskip\noindent
If $\nu\in\mathcal P_2(\mathbb R^d)$ and $\gamma\in\Gamma(\mu,\nu)$, then Jensen's inequality gives 
\[
    |m_\mu-m_\nu|^2 = \left|\int(x-y)\,d\gamma(x,y) \right|^2\leq\int|x-y|^2\,d\gamma(x,y).
\]
Choosing $\gamma\in\Gamma_o(\mu,\nu)$, we obtain
\begin{equation}\label{barycenter_lipschitz}
|m_\mu-m_\nu|\leq\mathcal W_2(\mu,\nu) \quad \text{ for every } \mu,\nu\in\mathcal P_2(\mathbb R^d).
\end{equation}
Finally, if $T:\mathbb R^d\to\mathbb R^d$ is given by $T(x):=a+sx$ for some $a\in\mathbb R^d$ and $s\in\mathbb R$, then the change of variables formula for push-forwards gives $m_{T_\#\mu}=a+s\,m_\mu$ and $\mathrm{Var}(T_\#\mu)=s^2\,\mathrm{Var}(\mu)$.

\subsection{Interpolations and geodesic convexity}
Given $\mu_0,\mu_1\in\mathcal P_2(\mathbb R^d)$ and $\gamma\in\Gamma(\mu_0,\mu_1)$, we write
\[
\mu_t^\gamma:=\bigl((1-t)\mathrm{proj}_1+t\mathrm{proj}_2\bigr)_\#\gamma \quad \text{ for every } t\in[0,1]
\]
for the interpolation induced by $\gamma$. By \eqref{plan_estimate}, this curve satisfies
\[
    \mathcal W_2(\mu^\gamma_s,\mu^\gamma_t)\leq|t-s| \left( \mathlarger{\int} |x-y|^2\,d\gamma \right)^{1/2},
\]
for all $s,t\in[0,1]$; when $\gamma\in\Gamma_o(\mu_0,\mu_1)$, all these inequalities are equalities and
$(\mu^\gamma_t)_{t\in[0,1]}$ is a constant speed geodesic joining $\mu_0$ to $\mu_1$, every
geodesic being of this form. 
A set $S\subseteq\mathcal P_2(\mathbb R^d)$ is geodesically convex
if $\mu^\gamma_t\in S$ for every $\mu_0,\mu_1\in S$, every $\gamma\in\Gamma_o(\mu_0,\mu_1)$
and every $t\in[0,1]$.

\subsection{Metric projections onto closed geodesically convex sets}
Given a nonempty, closed and convex set $K\subseteq\mathbb R^d$ and $v\in\mathbb R^d$, we denote by $\mathrm{proj}_K(v)$ the unique point of $K$ closest to $v$.
Recall that $\mathrm{proj}_K(v) \in K$ is characterized by 
\[
\left\langle x-\mathrm{proj}_K(v),v-\mathrm{proj}_K(v)\right\rangle\leq 0 \quad \text{ for every } x\in K;
\]
see \cite[Theorem~3.16]{BC_2017}. 
Analogously, given $S\subseteq\mathcal P_2(\mathbb R^d)$ and $\mu\in\mathcal P_2(\mathbb R^d)$, we set
\[
\mathrm{proj}_S(\mu):=\argmin_{\nu\in S}\mathcal W_2(\mu,\nu).
\]
In contrast with the Hilbert space setting, the function $\mathcal W_2^2(\cdot,\mu)$ is in general not convex along geodesics \cite[Example~9.1.5]{AGS_2008}. Therefore, the classical argument for the existence and uniqueness of the projection onto a nonempty closed geodesically convex set does not apply to an arbitrary $\mu$. 
It does apply, however, when $\mu$ is a Dirac mass. Thus $\mathrm{proj}_S(\delta_v)$ consists of exactly one element for all $v\in\mathbb R^d$

\subsection{Absolutely continuous curves and the continuity equation}\label{sec:conti_eq}
A curve $\mu:[t_0,+\infty)\to\mathcal P_2(\mathbb R^d)$ is locally absolutely continuous if there
exists $g\in L^1_{\mathrm{loc}}([t_0,+\infty))$ such that
\[
\mathcal W_2(\mu_s,\mu_t)\leq\int_s^t g(r)\,dr
\quad \text{ for every } t_0\leq s\leq t.
\]
For such a curve the metric derivative
\[
|\mu'|(t):=\lim_{h\to0}\frac{\mathcal W_2(\mu_{t+h},\mu_t)}{|h|}
\]
exists for almost every $t\in[t_0,+\infty)$, belongs to
$L^1_{\mathrm{loc}}([t_0,+\infty))$ and is the smallest admissible $g$ in the displayed
inequality; see \cite[Theorem~1.1.2]{AGS_2008}. Since the metric derivative of a locally absolutely continuous curve need not exist everywhere, all asymptotic statements involving $|\mu'|(t)$ are understood along the set of times at which the metric derivative exists.

\smallbreak\noindent
Given a Borel vector field $b:[t_0,+\infty)\times\mathbb R^d\to\mathbb R^d$, we say that a
narrowly continuous curve $\mu:[t_0,+\infty)\to\mathcal P(\mathbb R^d)$ solves the continuity
equation $\partial_t\mu_t+\mathrm{div}(b(t,\cdot)\mu_t)= 0$ if
$\displaystyle\int_{t_0}^\tau\int_{\mathbb R^d}|b(t,x)|\,d\mu_t(x)\,dt<+\infty$ for every $\tau>t_0$ and
\[
\frac{d}{dt}\int_{\mathbb R^d}\varphi(x)\,d\mu_t(x)
=
\int_{\mathbb R^d}\left\langle\nabla\varphi(x),b(t,x)\right\rangle\,d\mu_t(x)
\quad \text{ for every } \varphi\in C_c^\infty(\mathbb R^d)
\]
in the sense of distributions on $(t_0,+\infty)$. Testing against tensor products
$(t,x)\longmapsto\alpha(t)\varphi(x)$, with $\alpha\in C_c^\infty((t_0,+\infty))$ and
$\varphi\in C_c^\infty(\mathbb R^d)$, whose finite linear combinations are dense in
$C_c^\infty((t_0,+\infty)\times\mathbb R^d)$, shows that this formulation is equivalent to the
usual one. The flows studied in this paper are continuity equations whose velocity field depends
on the solution itself and their analysis rests on the superposition principle
\cite[Theorem~3.4]{AC_2014}, which represents any such solution as the law of a random absolutely
continuous curve solving the associated characteristic equation; the precise statement we need is given in Lemma~\ref{superposition}.

\smallbreak\noindent

We denote by $C([t_0,+\infty);\mathbb R^d)$ the space of continuous curves endowed with the
topology of uniform convergence on compact intervals, which makes it a Polish space and, for $t\in[t_0,+\infty)$, by $e_t:C([t_0,+\infty);\mathbb R^d)\to\mathbb R^d$ the evaluation mapping $e_t(\omega):=\omega(t)$, which is continuous. 
Finally, we shall use the disintegration theorem
in the form of \cite[Theorem~5.3.1]{AGS_2008}: if $\lambda$ is a probability measure on a Polish
space $Z$, $\pi:Z\to Y$ is Borel with values in a Polish space and $\sigma:=\pi_\#\lambda$,
then there exists a $\sigma$-a.e. uniquely determined Borel family
$(\lambda_y)_{y\in Y}\subseteq\mathcal P(Z)$ such that $\lambda_y$ is concentrated on
$\pi^{-1}(y)$ for $\sigma$-a.e. $y$ and
$\lambda= \mathlarger{\int}_Y\lambda_y\,d\sigma(y)$.

% \newpage 
\section{A superposition principle and a global existence theorem}
The classical long-time analysis of monotone flows usually relies on differentiating squared distances between trajectories and then using monotonicity. In the Wasserstein space, such arguments are not directly available, since the velocity field at each time belongs to a Lebesgue space attached to the current measure, so the underlying Hilbert space changes along the flow.
\smallbreak\noindent
The superposition principle provides a way to overcome this difficulty by representing a solution as a probability measure on path space whose time marginals are the evolving measures. The evaluation mappings then form a curve in a fixed Hilbert space, where the usual differentiation arguments can be carried out, and the resulting estimates pass to the Wasserstein distance through \eqref{lift_estimate}. The standard superposition principle, say \cite[Theorem~8.2.1]{AGS_2008}, does not apply directly under our assumptions, since the velocity fields are given as equivalence classes in varying Lebesgue spaces and satisfy only a weaker integrability condition in time. Moreover, our arguments require an evolution equation in a fixed Hilbert space, rather than just a pathwise identity.
\smallbreak\noindent
In the first subsection, we establish a version of the superposition principle suitable for these purposes. We then use it to prove a general well-posedness result for flows driven by time-dependent monotone operators.
\subsection{A superposition principle}

Let $\eta\in\mathcal P\bigl(C([t_0,+\infty);\mathbb R^d)\bigr)$ and let
$\mu:[t_0,+\infty)\to\mathcal P_2(\mathbb R^d)$ be such that
$(e_t)_\#\eta=\mu_t$ for every $t\in[t_0,+\infty)$. Given
$\xi\in L^2(\mu_t)^d$, we denote by $\xi(e_t)$ the pull-back of
$\xi$ through $e_t$, i.e., the class of $\widehat\xi\circ e_t$, where
$\widehat\xi$ is any Borel representative of $\xi$. This is a well-defined
element of $L^2(\eta)^d$. If $\widehat\xi_1=\widehat\xi_2$
$\mu_t$-a.e., then the identity $(e_t)_\#\eta=\mu_t$ gives
\[
\eta\bigl(\{\omega:\widehat\xi_1(\omega(t))\neq\widehat\xi_2(\omega(t))\}\bigr)
=
\mu_t\bigl(\{x:\widehat\xi_1(x)\neq\widehat\xi_2(x)\}\bigr)= 0,
\]
and the same identity yields $\|\xi(e_t)\|_{L^2(\eta)^d}=\|\xi\|_{L^2(\mu_t)^d}$.
Expressions such as $f_t[\mu_t](e_t)$ are always understood in this sense, as
pull-backs of $L^2$-equivalence classes and not as pointwise evaluations.

\medskip\noindent
The superposition principle in this subsection requires a velocity field which is jointly Borel, whereas the operators we consider produce, at each time, an element of an equivalence class with respect to a measure which itself varies with time.
The following lemma removes this discrepancy.

\begin{lemma}\label{borel_representative}
Let $f=(f_t[\nu])$ be a family indexed by
$(t,\nu)\in[t_0,+\infty)\times\mathcal P_2(\mathbb R^d)$ such that
$f_t[\nu]\in L^2(\nu)^d$ for every $t\in[t_0,+\infty)$ and every
$\nu\in\mathcal P_2(\mathbb R^d)$. Assume that
\begin{itemize}[label=\raisebox{0.2ex}{\scriptsize$\bullet$}]
    \item the mapping $t\longmapsto f_t[\nu]$ is Borel measurable from
    $[t_0,+\infty)$ into $L^2(\nu)^d$, endowed with its norm
    topology, for every $\nu\in\mathcal P_2(\mathbb R^d)$;

    \smallbreak
    \item the operator $\nu\longmapsto f_t[\nu]$ is demicontinuous for every
    $t\in[t_0,+\infty)$.
\end{itemize}
Let $\mu:[t_0,+\infty)\to\mathcal P_2(\mathbb R^d)$ be a continuous curve.
Then the mapping $[t_0,+\infty)\ni t\longmapsto\|f_t[\mu_t]\|_{L^2(\mu_t)^d}\in\mathbb R$ is Borel measurable and there exists a Borel mapping $\widehat f:[t_0,+\infty)\times\mathbb R^d\to\mathbb R^d$ such that $\widehat f(t,\cdot)\in L^2(\mu_t)^d$ for every $t\in[t_0,+\infty)$ and $\widehat f(t,\cdot)=f_t[\mu_t]$ in $L^2(\mu_t)^d$ for almost every $t\in[t_0,+\infty)$.
\end{lemma}
\begin{proof}
We first prove that the mapping
\begin{align*}
\mathfrak g: [t_0,+\infty)\times\mathcal P_2(\mathbb R^d) & \longrightarrow\mathcal P(\mathbb R^{2d}) \\
(t, \nu) &  \longmapsto \mathfrak{g}(t,\nu):= \big( \id,f_t[\nu] \big)_\#\nu
\end{align*}
is Borel measurable, where $\mathcal P(\mathbb R^{2d})$ is endowed with the
narrow topology.

\medskip\noindent
Fix $\nu\in\mathcal P_2(\mathbb R^d)$ and consider the mapping $G_\nu:L^2(\nu)^d\to\mathcal P(\mathbb R^{2d})$ given by $G_\nu(\xi):=(\id,\xi)_\#\nu$.
We claim that $G_\nu$ is continuous from the
norm topology of $L^2(\nu)^d$ to the narrow topology. 
Indeed, let $\xi_n\to\xi$ in $L^2(\nu)^d$ and let $\psi\in C_b(\mathbb R^{2d})$. 
Then $\xi_n\to\xi$ in $\nu$-measure and therefore $\psi(\cdot,\xi_n(\cdot))\to\psi(\cdot,\xi(\cdot))$ in $\nu$-measure as well, since every subsequence of $(\xi_n)_n$ admits a further subsequence converging $\nu$-a.e., along which the continuity of $\psi$ gives $\nu$-a.e. convergence. 
Consequently, we have, for every $\delta>0$,
\[
\int_{\mathbb R^d}\left|\psi(x,\xi_n(x))-\psi(x,\xi(x))\right|\,d\nu(x)
\leq \delta+2\|\psi\|_\infty \nu\bigl(\{x:\left|\psi(x,\xi_n(x))-\psi(x,\xi(x))\right|>\delta\}\bigr).
\]
Taking the limit superior as $n\to+\infty$ and then letting $\delta \searrow 0$, we obtain
\begin{align*}
    \lim_{n \to \infty} \int_{\R^{2d} } \psi \, dG_\nu(\xi_n) = \lim_{n \to \infty} \int_{\R^d} \psi(x, \xi_n(x)) \, d\nu(x) = \int_{\R^d} \psi(x, \xi(x)) \, d\nu(x) = \int_{\R^{2d}} \psi \, dG_\nu(\xi), 
\end{align*}
which proves the claim.
Since $t\longmapsto f_t[\nu]$ is Borel measurable by assumption and
$\mathfrak g(t,\nu)=G_\nu(f_t[\nu])$, the mapping $t\longmapsto\mathfrak g(t,\nu)$ is Borel measurable for every fixed $\nu$.

\medskip\noindent
Fix now $t\in[t_0,+\infty)$ and let $\nu_n\to\nu$ in $\mathcal W_2$. 
The demicontinuity of $\nu\longmapsto f_t[\nu]$ implies that $\mathfrak g(t,\nu_n)\rightharpoonup\mathfrak g(t,\nu)$ narrowly. 
Since $\mathcal P_2(\mathbb R^d)$, endowed with $\mathcal W_2$, and $\mathcal P(\mathbb R^{2d})$, endowed with the narrow topology, are metrizable,
sequential continuity is equivalent to continuity and hence $\nu\longmapsto\mathfrak g(t,\nu)$ is continuous. 
Thus $\mathfrak g$ is Borel measurable in its first variable and continuous in its second one and, the
space $\mathcal P_2(\mathbb R^d)$ being separable and metrizable, it is jointly
Borel measurable; see, e.g., \cite[Lemma~4.51]{AB_2006}.

\medskip\noindent
Since $[t_0,+\infty)\ni t\longmapsto(t,\mu_t)$ is continuous, we conclude that
\[
[t_0,+\infty)\ni t\longmapsto
\lambda_t:=\mathfrak g(t,\mu_t)=(\id,f_t[\mu_t])_\#\mu_t
\in\mathcal P(\mathbb R^{2d})
\]
is Borel measurable. 
The first assertion of the lemma follows at once since one has
\[
\|f_t[\mu_t]\|_{L^2(\mu_t)^d}^2=\int_{\mathbb R^{2d}}|z|^2\,d\lambda_t(x,z)
\] 
and the mapping $\lambda\longmapsto   \displaystyle\int|z|^2\,d\lambda$ is the increasing limit of the mappings $\lambda\longmapsto \displaystyle\int\min\{|z|^2,n\}\,d\lambda$, each of which is narrowly continuous because the corresponding integrand is bounded and continuous.
Hence $t\longmapsto\|f_t[\mu_t]\|_{L^2(\mu_t)^d}$ is Borel measurable.

\medskip\noindent
To show that $t\longmapsto\lambda_t(B)$ is Borel measurable
for every Borel set $B\subseteq\mathbb R^{2d}$, consider
the class
\[
\mathcal{D} := \{ B \in \mathcal{B}(\R^{2d}):\, t \longmapsto  \lambda_t(B)\, \text{ is Borel measurable} \}.
\]
This class contains all open sets, since $\lambda\longmapsto\lambda(O)$ is lower semicontinuous for every open $O\subseteq\mathbb R^{2d}$ and $t \longmapsto \lambda_t$ is Borel measurable.
Observe that $\mathcal D$ is a Dynkin system.
Since the open sets are stable under finite intersections and generate $\mathcal{B}(\R^{2d})$, Dynkin's $\pi$-$\lambda$ theorem yields $\mathcal D  = \mathcal{B}(\R^{2d})$, which is the desired assertion.
Thus $(t, B) \longmapsto\lambda_t(B)$ is a Borel probability kernel and consequently
\[
\Lambda(dt,dx,dz):=e^{t_0-t}\,dt\,\lambda_t(dx,dz)
\]
defines a Borel probability measure on
$[t_0,+\infty)\times\mathbb R^d\times\mathbb R^d$,
whose projection onto the first two variables is the probability measure
$q(dt,dx):=e^{t_0-t}\,dt\,\mu_t(dx)$, since the first marginal of $\lambda_t$ is
$\mu_t$. By the disintegration theorem, in the form recalled in Section~\ref{sec:conti_eq}, there exists a Borel family $(\nu_{t,x})$ of probability measures on $\mathbb R^d$ such that
$\Lambda(dt,dx,dz) = q(dt,dx)\,\nu_{t,x}(dz)$.

\medskip\noindent
Let $(\psi_j)_{j\in\mathbb N}\subseteq C_b(\mathbb R^{2d})$ be a countable
family determining probability measures on $\mathbb R^{2d}$ and set
\[
a_j(t):=\int_{\mathbb R^{2d}}\psi_j\,d\lambda_t
\qquad\text{and}\qquad
b_j(t):=\int_{\mathbb R^d}\int_{\mathbb R^d}
\psi_j(x,z)\,d\nu_{t,x}(z)\,d\mu_t(x).
\]
Both mappings are bounded and Borel measurable, the first one by the Borel measurability of $t\longmapsto\lambda_t$ and the second one because
$(t,x)\longmapsto \mathlarger{\int} \psi_j(x,z)\,d\nu_{t,x}(z)$ is bounded and Borel and $(t, B) \longmapsto\mu_t(B)$ is a Borel probability kernel. 
By the definition of $\Lambda$ and the disintegration identity, for every bounded Borel mapping
$\alpha:[t_0,+\infty)\to\mathbb R$ with compact support, one has
\[
\int_{t_0}^{+\infty}\alpha(t)e^{t_0-t}a_j(t)\,dt
=
\int_{[t_0,+\infty)\times\mathbb R^{2d}}\alpha(t)\psi_j(x,z)\,d\Lambda(t,x,z)
=
\int_{t_0}^{+\infty}\alpha(t)e^{t_0-t}b_j(t)\,dt.
\]
Since the weight $e^{t_0-t}$  is positive, $a_j = b_j$ almost everywhere on $[t_0,+\infty)$. 
Since the family $(\psi_j)_j$ is countable, there exists a Borel set $J\subseteq[t_0,+\infty)$ whose complement is Lebesgue negligible and such that $a_j(t)=b_j(t)$ for every $j\in\mathbb N$ and every $t\in J$. 
The fact that $(\psi_j)_j$ determines probability measures implies that $\lambda_t(dx,dz)=\mu_t(dx)\,\nu_{t,x}(dz)$ for every $t\in J$.

\medskip\noindent
Fix $t\in J$ and let $w_t:\mathbb R^d\to\mathbb R^d$ be a Borel representative of $f_t[\mu_t]$. 
Since $\lambda_t=(\id,w_t)_\#\mu_t$, the disintegration of $\lambda_t$ with respect to its first marginal is $x\longmapsto\delta_{w_t(x)}$ and the $\mu_t$-a.e. uniqueness of disintegrations yields $\nu_{t,x}=\delta_{w_t(x)}$ for $\mu_t$-a.e. $x\in\mathbb R^d$.

\medskip\noindent
It remains to construct a Borel representative from the kernel $(\nu_{t,x})$. 
Set 
\[
r(t,x):=\int_{\mathbb R^d}|z|\,d\nu_{t,x}(z)\in[0,+\infty],
\]
which is Borel as the increasing limit of the Borel mappings 
\[
(t,x)\longmapsto \displaystyle\int\min\{|z|,n\}\,d\nu_{t,x}(z).
\]
Set $\mathcal A:=\{(t,x):r(t,x)<+\infty\}$. 
For every $i\in\{1,\cdots,d\}$ and every $n\in\mathbb N$, the mapping
\[
(t,x)\longmapsto\int\max\{-n,\min\{z_i,n\}\}\,d\nu_{t,x}(z)
\]
is Borel and on $\mathcal A$ these mappings converge, by dominated convergence, to 
\[
(t,x)\longmapsto\int z_i\,d\nu_{t,x}(z).
\]
Hence the mapping $\widehat f:[t_0,+\infty)\times\mathbb R^d\to\mathbb R^d$ given by
\[
\widehat f(t,x):= 
\begin{cases}
\displaystyle\int_{\mathbb R^d}z\,d\nu_{t,x}(z)
&\quad\text{if}\quad t\in J\quad\text{and}\quad r(t,x)<+\infty,\\
0&\quad\text{otherwise},
\end{cases}
\]
is Borel, because both $J$ and $\mathcal A$ are Borel.
If $t\in J$, then, for $\mu_t$-a.e. $x$, we have $r(t,x)=|w_t(x)|<+\infty$ and $\widehat f(t,x)=w_t(x)$ so that $\widehat f(t,\cdot)=f_t[\mu_t]$ in $L^2(\mu_t)^d$, while $\widehat f(t,\cdot)= 0$ for $t\notin J$.
In particular, $\widehat f(t,\cdot)\in L^2(\mu_t)^d$ for every $t\in[t_0,+\infty)$ and $\widehat f(t,\cdot)=f_t[\mu_t]$ in $L^2(\mu_t)^d$ for almost every $t\in[t_0,+\infty)$.
\end{proof}

We are now ready to state the aforementioned superposition principle.
\begin{lemma}\label{superposition}
Let $f=(f_t[\nu])$ indexed by $(t,\nu)\in[t_0,+\infty)\times\mathcal P_2(\mathbb R^d)$ be such that $f_t[\nu]\in L^2(\nu)^d$ for every $t\in[t_0,+\infty)$ and every $\nu\in\mathcal P_2(\mathbb R^d)$. 
Assume that
\begin{itemize}[label=\raisebox{0.2ex}{\scriptsize$\bullet$}]
    \item the mapping $t\longmapsto f_t[\nu]$ is Borel measurable from $[t_0,+\infty)$ into $L^2(\nu)^d$, endowed with its norm topology, for every $\nu\in\mathcal P_2(\mathbb R^d)$;

    \smallbreak
    \item the operator $\nu\longmapsto f_t[\nu]$ is demicontinuous for every $t\in[t_0,+\infty)$.
\end{itemize}
Let $\mu_0 \in \mathcal{P}_2(\R^d)$ and let $\mu:[t_0,+\infty)\to\mathcal P_2(\mathbb R^d)$ be a locally absolutely continuous curve such that
\begin{equation}\label{eq.superposition_hypothesis}
    \begin{cases}
        \partial_t\mu_t -\mathrm{div}\bigl(f_t[\mu_t]\mu_t\bigr)= 0,\\[0.3em]
        \mu_{t_0}=\mu_0,
    \end{cases}
    \qquad\text{and}\qquad
    \int_{t_0}^\tau\|f_t[\mu_t]\|_{L^2(\mu_t)^d}\,dt<+\infty \quad \text{for every } \tau > t_0.
\end{equation}
Then there exists a probability measure $\eta\in\mathcal P\bigl(C([t_0,+\infty);\mathbb R^d)\bigr)$ concentrated on $\mathrm{AC}_{\mathrm{loc}}([t_0,+\infty);\mathbb R^d)$ and satisfying the following properties.
\begin{itemize}[label=\raisebox{0.2ex}{\scriptsize$\bullet$}]
    \item $(e_t)_\#\eta=\mu_t$ for all $t\in[t_0,+\infty)$;
    \smallbreak
    \item if $\widehat f$ is a Borel mapping as in Lemma~\ref{borel_representative}, then
    \[
        \dot\omega(t)+\widehat f(t,\omega(t))= 0
    \]
    for $\eta$-a.e. $\omega$ and almost every $t\in[t_0,+\infty)$; further, any two such mappings agree along $\eta$-a.e. curve for almost every $t$, so that  we may write $\dot\omega(t)+f_t[\mu_t](\omega(t))= 0$;
    \smallbreak
    \item the mapping $[t_0,+\infty)\ni t\longmapsto e_t\in L^2(\eta)^d$ is     locally absolutely continuous and
    \[
        \dot e_t+f_t[\mu_t](e_t)= 0
    \]
    in $L^2(\eta)^d$ for almost every $t\in[t_0,+\infty)$.
\end{itemize}
\end{lemma}

\begin{proof}
Thanks to Lemma~\ref{borel_representative}, the mapping $t\longmapsto\|f_t[\mu_t]\|_{L^2(\mu_t)^d}$ is Borel measurable and so we notice that the integrability requirement above is meaningful.
For each $n\in\mathbb N\cup\{0\}$, set $s_n:=t_0+n$ and $I_n:=[s_n,s_{n+1}]$. 
Let $\widehat f$ be the Borel mapping provided by Lemma~\ref{borel_representative} and set $b(t,x):=-\widehat f(t,x)$. 
Since $\widehat f(t,\cdot) = f_t[\mu_t]$ in $L^2(\mu_t)^d$ for almost every $t$ and  each $\mu_t$ is a probability measure, the H\"older inequality then yields
\begin{equation}\label{borel_velocity}
\int_{t_0}^\tau\int_{\mathbb R^d}|b(t,x)|\,d\mu_t(x)\,dt \leq \int_{t_0}^\tau\|b(t,\cdot)\|_{L^2(\mu_t)^d}\,dt = \int_{t_0}^\tau\|f_t[\mu_t]\|_{L^2(\mu_t)^d}\,dt<+\infty \quad \text{for every } \tau > t_0.
\end{equation}
Moreover, as observed in Section~\ref{sec:conti_eq}, the distributional identity in \eqref{eq.superposition_hypothesis} amounts to $\partial_t\mu_t+\mathrm{div}(b(t,\cdot)\mu_t)= 0$ in the sense of distributions on $(t_0,+\infty)\times\mathbb R^d$.
We divide the remainder of the proof into three steps.

\medskip\noindent
\textbf{Step 1:} \emph{Superposition on compact intervals.}
Fix $n\in\mathbb N\cup\{0\}$. The curve $t\longmapsto\mu_t$ is narrowly continuous on $I_n$ and solves the continuity equation with the Borel velocity field $b$, which satisfies \eqref{borel_velocity}. 
Thanks to the integrability condition~\eqref{borel_velocity}, the superposition principle for the continuity equation \cite[Theorem~3.4]{AC_2014} yields a measure $\eta^n\in\mathcal P\bigl(C(I_n;\mathbb R^d)\bigr)$ concentrated on curves $\omega\in\mathrm{AC}(I_n;\mathbb R^d)$ satisfying $\dot\omega(t)=b(t,\omega(t))$ for almost every $t\in I_n$ and $(e_t)_\#\eta^n=\mu_t$ for every $t\in I_n$.

\medskip\noindent
\textbf{Step 2:} \emph{Concatenation and the pathwise equation.}
We construct by induction $\vartheta^n\in\mathcal P\bigl(C([t_0,s_{n+1}];\mathbb R^d)\bigr)$, concentrated on curves $\omega\in\mathrm{AC}([t_0,s_{n+1}];\mathbb R^d)$ with $\dot\omega(t)=b(t,\omega(t))$ for almost every $t$, such that
\begin{itemize}
    \item $(e_t)_\#\vartheta^n=\mu_t$ for all $t\in[t_0,s_{n+1}]$;

    \item the restriction of $\vartheta^{n+1}$ to $[t_0,s_{n+1}]$ equals $\vartheta^n$.
\end{itemize}
Set $\vartheta^0:=\eta^0$. 
Given $\vartheta^n$, disintegrate $\eta^{n+1}$ with respect to the continuous evaluation mapping $e_{s_{n+1}}$.
Since $(e_{s_{n+1}})_\#\eta^{n+1}=\mu_{s_{n+1}}$, the disintegration theorem provides a Borel probability kernel $x\longmapsto\eta^{n+1}_x$ such that 
\[
\eta^{n+1}=\int_{\mathbb R^d}\eta^{n+1}_x\,d\mu_{s_{n+1}}(x)
\]
and $\eta^{n+1}_x$ is concentrated on the curves with $\omega(s_{n+1})=x$ for $\mu_{s_{n+1}}$-a.e. $x \in \R^d$.
Let $\rho^{n+1}$ be the probability measure on $C([t_0,s_{n+1}];\mathbb R^d)\times C(I_{n+1};\mathbb R^d)$ defined by
\[
\rho^{n+1}(d\omega_1,d\omega_2) :=\vartheta^n(d\omega_1)\,\eta^{n+1}_{\omega_1(s_{n+1})}(d\omega_2),
\]
which is well defined because $\omega_1\longmapsto\omega_1(s_{n+1})$ is continuous and $x\longmapsto\eta^{n+1}_x$ is a Borel kernel.
Since the law of $\omega_1(s_{n+1})$ under $\vartheta^n$ is $\mu_{s_{n+1}}$, the closed set $D:=\{(\omega_1,\omega_2):\omega_1(s_{n+1})=\omega_2(s_{n+1})\}$ satisfies $\rho^{n+1}(D)=1$.
On $D$, let $\mathcal C = \mathcal C(\omega_1, \omega_2)$ denote the concatenation mapping which sends $(\omega_1, \omega_2)$ to the curve equal to $\omega_1$ on $[t_0, s_{n +1}]$ and to $\omega_2$ on $I_{n + 1}$.
This defines a continuous mapping from $D$ into $C([t_0,s_{n+2}];\mathbb R^d)$ since the uniform distance between two concatenated curves is the maximum of the uniform distances between their pieces.
We extend $\mathcal C$ to a Borel mapping on the whole product space by assigning a fixed curve outside $D$, which is immaterial since $\rho^{n + 1}(D) = 1$.
Set $\vartheta^{n+1}:=\mathcal C_\#\rho^{n+1}$.
For $t\in[t_0,s_{n+1}]$, we have $(e_t)_\#\vartheta^{n+1}=(e_t)_\#\vartheta^n=\mu_t$.
For $t\in I_{n+1}$, the identity $(e_{s_{n+1}})_\#\vartheta^n=\mu_{s_{n+1}}$ gives
\[
(e_t)_\#\vartheta^{n+1} =\int_{\mathbb R^d}(e_t)_\#\eta^{n+1}_x\,d\mu_{s_{n+1}}(x) =(e_t)_\#\eta^{n+1} =\mu_t.
\]
By the same argument, the restriction of $\vartheta^{n+1}$ to $[t_0,s_{n+1}]$ is $\vartheta^n$. Finally, for $\rho^{n+1}$-a.e. $(\omega_1,\omega_2)$, both pieces are absolutely continuous and their values agree at $s_{n+1}$; joining their derivatives produces an $L^1$ function on $[t_0,s_{n+2}]$ whose integral recovers the concatenated curve, so that $\mathcal C(\omega_1,\omega_2)$ belongs to $\mathrm{AC}([t_0,s_{n+2}];\mathbb R^d)$ and 
\[
\frac{d}{dt}\mathcal C(\omega_1,\omega_2)(t) = b(t,\mathcal C(\omega_1,\omega_2)(t)) \quad \text{ for almost every $t\in[t_0,s_{n+2}]$.}
\]
This completes the induction.

\medskip\noindent
The spaces $C([t_0,s_{n+1}];\mathbb R^d)$ are Polish and the restriction mappings between them are continuous and surjective.
Moreover, $C([t_0,+\infty);\mathbb R^d)$, endowed with the topology of uniform convergence on compact intervals, is canonically identified with their projective limit, since every compatible family of curves determines a unique continuous curve on $[t_0,+\infty)$ and conversely.
Hence, by the projective limit theorem for probability measures on Polish spaces (see, e.g., \cite[Chapter~V]{P_1967}), the consistent family
$(\vartheta^n)_{n}$ admits a unique projective limit
$\eta\in\mathcal P\bigl(C([t_0,+\infty);\mathbb R^d)\bigr)$ whose restriction to $C([t_0,s_{n+1}];\mathbb R^d)$ is $\vartheta^n$ for every $n\geq 0$. 
Consequently $(e_t)_\#\eta=\mu_t$ for all $t\in[t_0,+\infty)$, which proves the first assertion.

\medskip

\noindent
Further, for each $n$, the set of curves whose restriction to $[t_0,s_{n+1}]$ is absolutely continuous and satisfies $\dot\omega(t)=b(t,\omega(t))$ for almost every $t\in[t_0,s_{n+1}]$ has full $\eta$-measure.
This set is Borel, being described by the conditions
\[
    \int_{t_0}^{s_{n+1}}|b(r,\omega(r))|\,dr<+\infty \quad \text{ and } \quad \omega(t)=\omega(t_0)+\int_{t_0}^tb(r,\omega(r))\,dr \quad \text{for every rational $t\in[t_0,s_{n+1}]$},
\]
which are Borel because $b$ is Borel and $(t,\omega)\longmapsto\omega(t)$ is continuous. 
Intersecting these countably many sets, we conclude that $\eta$ is concentrated on $\mathrm{AC}_{\mathrm{loc}}([t_0,+\infty);\mathbb R^d)$ and that $\dot\omega(t)+\widehat f(t,\omega(t))= 0$ for $\eta$-a.e. $\omega$ and almost every $t\in[t_0,+\infty)$. This identity does not depend on the chosen Borel representative because if $\widetilde f$ is another Borel mapping with $\widetilde f(t,\cdot)=\widehat f(t,\cdot)$ $\mu_t$-a.e. for almost every $t$, then, using $(e_t)_\#\eta=\mu_t$,
\[
\int_{t_0}^\tau\int \mathbbm 1_{\{\widehat f(t,\omega(t))\neq\widetilde f(t,\omega(t))\}}\,d\eta(\omega)\,dt
= \int_{t_0}^\tau \mu_t\bigl(\{x:\widehat f(t,x)\neq\widetilde f(t,x)\}\bigr)\,dt= 0
\]
for every $\tau>t_0$ and Fubini's theorem shows that, for $\eta$-a.e. $\omega$, the two mappings agree along $\omega$ for almost every $t$. 
This proves the second assertion.

\medskip\noindent
\textbf{Step 3:} \emph{The $L^2(\eta)^d$-valued equation.}
For every $t\in[t_0,+\infty)$, the identity $(e_t)_\#\eta=\mu_t$ gives 
\[
\int|e_t|^2\,d\eta=\mathcal M_2(\mu_t)<+\infty \quad \text{ and so } \quad e_t\in L^2(\eta)^d.
\]
Set $V_t:=b(t,e_t)$, i.e., $V_t(\omega)=b(t,\omega(t))$.
Consequently, $\|V_t\|_{L^2(\eta)^d}=\|b(t,\cdot)\|_{L^2(\mu_t)^d}$ and $V_t \in L^2(\eta)^d$ for every $t\in[t_0,+\infty)$, since $\widehat f(t,\cdot)\in L^2(\mu_t)^d$ for every $t\in[t_0,+\infty)$.
Also, thanks to~\eqref{borel_velocity}, the mapping $t\longmapsto\|V_t\|_{L^2(\eta)^d}$ is locally integrable.
Since $b$ is Borel and the evaluation $(t,\omega)\longmapsto\omega(t)$ is  continuous on $[t_0,+\infty)\times C([t_0,+\infty);\mathbb R^d)$, the mapping $(t,\omega)\longmapsto V_t(\omega)$ is jointly Borel and hence $t\longmapsto V_t\in L^2(\eta)^d$ is weakly measurable, it is therefore strongly measurable and Bochner integrable on compact subintervals of $[t_0,+\infty)$. 
By Step 2, for $\eta$-a.e. $\omega$ and all $t_0\leq s\leq t$,
\[
\omega(t)-\omega(s)=\int_s^t b(r,\omega(r))\,dr.
\]
Given $g\in L^2(\eta)^d$, integrating this identity against $g\,d\eta$ and using Fubini's theorem, we obtain
\[
\langle e_t-e_s,g\rangle_{L^2(\eta)^d} =\int_s^t\langle V_r,g\rangle_{L^2(\eta)^d}\,dr \quad 
\text{ and hence } \quad
e_t-e_s=\int_s^t V_r\,dr\qquad\text{in }L^2(\eta)^d,
\]
the integral being a Bochner integral.
It follows that $t\longmapsto e_t\in L^2(\eta)^d$ is locally absolutely continuous and that $\dot e_t=V_t$ for almost every $t\in[t_0,+\infty)$.
Finally, since $\widehat f(t,\cdot)=f_t[\mu_t]$ in $L^2(\mu_t)^d$ for almost every $t$, pulling this identity back through $e_t$ gives $V_t=-f_t[\mu_t](e_t)$ in
$L^2(\eta)^d$ for almost every $t$, which proves the last assertion.
\end{proof}
Before closing this subsection, we make a couple of remarks. 

\begin{remark}
Under the hypotheses of Lemma~\ref{superposition} alone, the measure $\eta$ need not be induced by a flow map, i.e., it need not be of the form $X_\#\mu_0$ for a Borel mapping $X:\mathbb R^d\to C([t_0,+\infty);\mathbb R^d)$, as  characteristics issuing from the same initial point may branch. However, when the monotonicity condition \eqref{general_lambda_monotonicity} holds, such branching is excluded and $\eta$ is induced by a flow map.
\end{remark}

\begin{remark}
Combining \eqref{lift_estimate} with the last assertion of Lemma~\ref{superposition} gives
\[
\mathcal W_2(\mu_s,\mu_t)\leq\|e_t-e_s\|_{L^2(\eta)^d}\leq\int_s^t\|f_r[\mu_r]\|_{L^2(\mu_r)^d}\,dr \quad \text{for every } t_0\leq s\leq t,
\]
so that $|\mu'|(t)\leq\|f_t[\mu_t]\|_{L^2(\mu_t)^d}$ for almost every $t\in[t_0,+\infty)$. Equality need not hold, since the driving field $-f_t[\mu_t]$ need not belong to the tangent space of $\mathcal P_2(\mathbb R^d)$ at $\mu_t$, i.e., it need not be the velocity field of minimal norm among those driving $\mu$; see \cite[Theorem~8.3.1 and Proposition~8.4.5]{AGS_2008}.
\end{remark}

\subsection{A global existence theorem}
We now establish an existence and uniqueness result for a general time-dependent
continuity equation whose velocity field is monotone. It can directly be applied both to the monotone operator
flow \eqref{eqtflw} and to its Tikhonov regularization \eqref{eqtfl} in
Section~\ref{s2}.

\begin{theorem}\label{general_existence}
Let $f=(f_t[\nu])$ be a family indexed by
$(t,\nu)\in[t_0,+\infty)\times\mathcal P_2(\mathbb R^d)$ such that
$f_t[\nu]\in L^2(\nu)^d$ for every $t\in[t_0,+\infty)$ and every
$\nu\in\mathcal P_2(\mathbb R^d)$. Assume that
\begin{itemize}[label=\raisebox{0.2ex}{\scriptsize$\bullet$}]
    \item the mapping $t\longmapsto f_t[\nu]$ is Borel measurable from
    $[t_0,+\infty)$ into $L^2(\nu)^d$, endowed with its norm
    topology, for every $\nu\in\mathcal P_2(\mathbb R^d)$;

    \smallbreak
    \item the operator $\nu\longmapsto f_t[\nu]$ is demicontinuous for every
    $t\in[t_0,+\infty)$;

    \smallbreak
    \item there exist $\lambda\in L^1_{\mathrm{loc}}([t_0,+\infty))$ and a
    Lebesgue negligible set $N \subseteq[t_0,+\infty)$ such that, for every
    $t\in[t_0,+\infty)\setminus N$,
    \begin{equation}\label{general_lambda_monotonicity}
        \int_{\mathbb R^{2d}} \left\langle f_t[\mu](x)-f_t[\nu](y),x-y \right\rangle        \,d\gamma(x,y) \geq \lambda(t) \int_{\mathbb R^{2d}}|x-y|^2\,d\gamma(x,y)
    \end{equation}
    for every $\mu,\nu\in\mathcal P_2(\mathbb R^d)$ and every
    $\gamma\in\Gamma(\mu,\nu)$;

    \smallbreak
    \item for every $\tau>t_0$ and every compact set $K\subseteq\mathcal P_2(\mathbb R^d)$, there exist $g_{\tau,K}\in L^1(t_0,\tau)$ and a Lebesgue negligible set   $N_{\tau,K}\subseteq(t_0,\tau)$ such that
    \[
        \|f_t[\nu]\|_{L^2(\nu)^d}\leq g_{\tau,K}(t) 
        \quad \text{for every } \nu\in K \text{ and every $t\in(t_0,\tau)\setminus N_{\tau,K}$.}
    \]
\end{itemize}
Then, for every $\mu_0\in\mathcal P_2(\mathbb R^d)$, there exists a unique
locally absolutely continuous curve
$\mu:[t_0,+\infty)\to\mathcal P_2(\mathbb R^d)$ satisfying
\begin{equation}\label{general_flow}
    \begin{cases} 
        \partial_t\mu_t -\mathrm{div}\bigl(f_t[\mu_t]\mu_t\bigr)= 0,\\[0.3em]
        \mu_{t_0}=\mu_0,
    \end{cases}
    \qquad\text{and}\qquad
    \int_{t_0}^\tau\|f_t[\mu_t]\|_{L^2(\mu_t)^d}\,dt<+\infty
    \qquad\text{for every } \tau>t_0.
\end{equation}
Moreover, if $\mu$ and $\nu$ are the solutions corresponding to the initial
data $\mu_0$ and $\nu_0$, respectively, then
\begin{equation}\label{general_contraction}
    \mathcal W_2(\mu_t,\nu_t)
    \leq
    \exp\left(-\int_{t_0}^t\lambda(s)\,ds\right) \mathcal W_2(\mu_0,\nu_0) \qquad\text{for every } t\geq t_0.
\end{equation}
\end{theorem}
\begin{proof}
We divide the argument into three steps. 
In the first one we prove the contraction estimate \eqref{general_contraction}, which contains in particular the uniqueness assertion. 
In the second one we solve \eqref{general_flow} when
the initial datum is finitely supported and in the third one we let the initial datum vary and pass to the limit.

\medskip\noindent
\textbf{Step 1:} \emph{Contraction and uniqueness.}
Let $\mu_0,\nu_0\in\mathcal P_2(\mathbb R^d)$ and let $\mu$ and $\nu$ be solutions of \eqref{general_flow} with initial data $\mu_0$ and $\nu_0$, respectively.
Both curves satisfy the hypotheses of Lemma~\ref{superposition}, so
there exist
$\eta,\widetilde\eta\in\mathcal P\bigl(C([t_0,+\infty);\mathbb R^d)\bigr)$,
concentrated on $\mathrm{AC}_{\mathrm{loc}}([t_0,+\infty);\mathbb R^d)$,
such that, writing $e^\eta_t$ and $e^{\widetilde\eta}_t$ for the evaluation
mappings at time $t$ on the two path spaces,
$(e^\eta_t)_\#\eta=\mu_t$ and $(e^{\widetilde\eta}_t)_\#\widetilde\eta=\nu_t$
for every $t\in[t_0,+\infty)$.
Moreover, the mappings $t\longmapsto e^\eta_t$ and
$t\longmapsto e^{\widetilde\eta}_t$ are locally absolutely
continuous with values in $L^2(\eta)^d$ and
$L^2(\widetilde\eta)^d$ and satisfy
\[
\dot e^\eta_t+f_t[\mu_t](e^\eta_t)= 0 \quad\text{in }L^2(\eta)^d \quad \text{ and } \quad \dot e^{\widetilde\eta}_t+f_t[\nu_t](e^{\widetilde\eta}_t)= 0 \quad\text{in }L^2(\widetilde\eta)^d
\]
for almost every $t\in[t_0,+\infty)$, respectively.

\medskip\noindent
We couple the two measures in such a way that the initial values are optimally
coupled. Since $(e^\eta_{t_0})_\#\eta=\mu_0$, the disintegration theorem provides a
Borel family $(\eta_x)_{x\in\mathbb R^d}$ of probability measures on $C([t_0,+\infty);\mathbb R^d)$ such that $\eta = \mathlarger{\int}_{\mathbb R^d}\eta_x\,d\mu_0(x)$ and $\eta_x$ is concentrated on the curves issuing from $x$ for $\mu_0$-a.e. $x$.
Let $(\widetilde\eta_y)_{y\in\mathbb R^d}$ be the analogous family for $\widetilde\eta$.
Fix $\gamma_0\in\Gamma_o(\mu_0,\nu_0)$ and set
\[
\Theta:=\int_{\mathbb R^{2d}}\eta_x\otimes\widetilde\eta_y\,d\gamma_0(x,y) \in\mathcal P\bigl(C([t_0,+\infty);\mathbb R^d)\times C([t_0,+\infty);\mathbb R^d)\bigr),
\]
and denote by $E_t$ and $\widetilde E_t$ the evaluations of the first and of
the second curve, i.e., $E_t(\omega,\widetilde\omega):=\omega(t)$ and
$\widetilde E_t(\omega,\widetilde\omega):=\widetilde\omega(t)$. Since the
marginals of $\gamma_0$ are $\mu_0$ and $\nu_0$, the marginals of
$\Theta$ are $\eta$ and $\widetilde\eta$, so that $(E_t)_\#\Theta=\mu_t$
and $(\widetilde E_t)_\#\Theta=\nu_t$ for every $t\in[t_0,+\infty)$.
Moreover, since the measures $\eta_x\otimes\widetilde\eta_y$ are concentrated on pairs of curves issuing from $x$ and $y$, we have
$(E_{t_0},\widetilde E_{t_0})_\#\Theta=\gamma_0$.

\medskip\noindent
The pull-back through the projection onto the first factor is a linear isometry
from $L^2(\eta)^d$ into $L^2(\Theta)^d$ carrying $e^\eta_t$ to
$E_t$; as bounded linear mappings commute with Bochner integration,
$t\longmapsto E_t\in L^2(\Theta)^d$ is locally absolutely continuous
and $\dot E_t=-f_t[\mu_t](E_t)$ for almost every $t$ and analogously for
$\widetilde E$. Consequently, the function
$t\longmapsto\dfrac{1}{2}\|E_t-\widetilde E_t\|^2_{L^2(\Theta)^d}$ is
locally absolutely continuous and, for a.e. $t \in [t_0,+\infty)$,
\[
\frac{d}{dt}\,\frac12\|E_t-\widetilde E_t\|^2_{L^2(\Theta)^d}
=
-\int \left\langle f_t[\mu_t](E_t)-f_t[\nu_t](\widetilde E_t),E_t-\widetilde E_t \right\rangle\,d\Theta.
\]
The measure $\gamma_t:=(E_t,\widetilde E_t)_\#\Theta$ belongs to
$\Gamma(\mu_t,\nu_t)$ and the right-hand side above equals
\[
-\int_{\mathbb R^{2d}} \left\langle f_t[\mu_t](x)-f_t[\nu_t](y),x-y\right\rangle\,d\gamma_t(x,y)
\leq 
-\lambda(t)\int_{\mathbb R^{2d}}|x-y|^2\,d\gamma_t(x,y)
=
-\lambda(t)\|E_t-\widetilde E_t\|^2_{L^2(\Theta)^d}
\]
by \eqref{general_lambda_monotonicity}. It follows that the locally absolutely continuous function
\[
t\longmapsto\exp\left( 2\int_{t_0}^t\lambda(s)\,ds \right) \|E_t-\widetilde E_t\|^2_{L^2(\Theta)^d}
\]
has nonpositive derivative
almost everywhere, hence is nonincreasing and therefore
\[
\|E_t-\widetilde E_t\|^2_{L^2(\Theta)^d}
\leq
\exp\left(-2\int_{t_0}^t\lambda(s)\,ds\right) \int_{\mathbb R^{2d}}|x-y|^2\,d\gamma_0(x,y)
=
\exp\left(-2\int_{t_0}^t\lambda(s)\,ds\right)\mathcal W_2^2(\mu_0,\nu_0).
\]
Since $\gamma_t\in\Gamma(\mu_t,\nu_t)$, estimate \eqref{lift_estimate} gives $\mathcal W_2(\mu_t,\nu_t)\leq\|E_t-\widetilde E_t\|_{L^2(\Theta)^d}$ and \eqref{general_contraction} follows.
In particular, if $\mu_0=\nu_0$, then $\mu_t=\nu_t$ for every $t\in[t_0,+\infty)$, which is the uniqueness assertion.

\medskip\noindent
\textbf{Step 2:} \emph{Finitely supported initial data.}
Let $N\in\mathbb N$ and let $m_1,\cdots,m_N>0$ satisfy
$\sum_{i=1}^Nm_i=1$. 
For $\mathbf x=(x_1,\cdots,x_N)\in(\mathbb R^d)^N$, we set $\nu_{\mathbf x} := \sum_{i=1}^N m_i\delta_{x_i}$ and we endow $(\mathbb R^d)^N$ with the scalar product $\langle\mathbf x,\mathbf y\rangle_{\mathbf m} := \sum_{i=1}^N m_i\langle x_i,y_i\rangle$ and the associated norm. Since
$\sum_{i=1}^N m_i\delta_{(x_i,y_i)}\in\Gamma(\nu_{\mathbf x},\nu_{\mathbf y})$,
inequality \eqref{plan_estimate} gives
$\mathcal W_2(\nu_{\mathbf x},\nu_{\mathbf y})
\leq\|\mathbf x-\mathbf y\|_{\mathbf m}$, so that
$\mathbf x\longmapsto\nu_{\mathbf x}$ is continuous from $(\mathbb R^d)^N$
into $\mathcal P_2(\mathbb R^d)$. Every atom of $\nu_{\mathbf x}$ has
positive mass, whence any two Borel representatives of $f_t[\nu_{\mathbf x}]$
coincide at each $x_i$ and the mapping
\[
\mathbf F_t(\mathbf x)
:=
\bigl(f_t[\nu_{\mathbf x}](x_1),\cdots,f_t[\nu_{\mathbf x}](x_N)\bigr)
\]
is well defined; grouping the indices according to the coincidences among the
$x_i$, one checks that $\|\mathbf F_t(\mathbf x)\|_{\mathbf m} =\|f_t[\nu_{\mathbf x}]\|_{L^2(\nu_{\mathbf x})^d}$. We consider the system
\begin{equation}\label{particle_system}
\begin{dcases}
\dot{\mathbf x}(t)+\mathbf F_t(\mathbf x(t))= 0,  \\ \mathbf x(t_0)=\mathbf x^0.
\end{dcases}
\end{equation}
\smallbreak\noindent
For any fixed $\mathbf x$, the mapping $t\longmapsto\mathbf F_t(\mathbf x)$ is Borel measurable, because $\xi\longmapsto(\xi(x_1),\cdots,\xi(x_N))$ is a linear isometry from $L^2(\nu_{\mathbf x})^d$ into $\bigl((\mathbb R^d)^N,\|\cdot\|_{\mathbf m}\bigr)$. 
For any fixed $t$, the mapping $\mathbf{x} \longmapsto \mathbf F_t(\mathbf{x})$ is continuous.
Indeed, let $\mathbf x^n\to\mathbf x$ in $(\mathbb R^d)^N$ and write $F^n_i:=f_t[\nu_{\mathbf x^n}](x^n_i)$ and $F_i:=f_t[\nu_{\mathbf x}](x_i)$. Since $\nu_{\mathbf x^n}\to\nu_{\mathbf x}$
in $\mathcal W_2$, demicontinuity gives
\begin{equation}\label{eq.discrete_demicontinuity}
\sum_{i=1}^Nm_i\delta_{(x^n_i,F^n_i)} =(\id,f_t[\nu_{\mathbf x^n}])_\#\nu_{\mathbf x^n}
\rightharpoonup (\id,f_t[\nu_{\mathbf x}])_\#\nu_{\mathbf x}
=\sum_{i=1}^Nm_i\delta_{(x_i,F_i)}
\end{equation}
narrowly, together with
$\sup_n\|f_t[\nu_{\mathbf x^n}]\|_{L^2(\nu_{\mathbf x^n})^d}<+\infty$. Since
$m_i|F^n_i|^2\leq\|f_t[\nu_{\mathbf x^n}]\|^2_{L^2(\nu_{\mathbf x^n})^d}$ and
$m_i>0$, the sequence $(\mathbf F_t(\mathbf x^n))_n$ is bounded in
$(\mathbb R^d)^N$. Let $\overline{\mathbf F}$ be the limit of a convergent
subsequence. Along it, the left-hand side of
\eqref{eq.discrete_demicontinuity} converges narrowly to
$\sum_im_i\delta_{(x_i,\overline F_i)}$, whence
\[
\sum_{i=1}^Nm_i\delta_{(x_i,\overline F_i)}
=
\sum_{i=1}^Nm_i\delta_{(x_i,F_i)}.
\]
Fix $i$ and set $I:=\{j:x_j=x_i\}$ and $m_I:=\sum_{j\in I}m_j$.
Disintegrating both sides with respect to their common first marginal
$\nu_{\mathbf x}$, evaluating at the atom $x_i$ and using that $F_j=F_i$
for $j\in I$, we obtain $m_I^{-1}\sum_{j\in I}m_j\delta_{\overline F_j}=\delta_{F_i}$.
Since a convex combination of Dirac masses with positive weights is a Dirac mass only if all of them are concentrated at the same point,
this forces $\overline F_j=F_i$ for every $j\in I$ and in particular $\overline F_i = F_i$.
Thus every subsequence of $(\mathbf F_t(\mathbf x^n))_n$ admits a further subsequence converging to $\mathbf F_t(\mathbf x)$, which proves the continuity of $\mathbf F_t$.

\medskip\noindent
Finally, given $\tau>t_0$ and $R>0$, the set $K_R:=\{\nu_{\mathbf x}:\|\mathbf x\|_{\mathbf m}\leq R\}$ is a compact subset of $\mathcal P_2(\mathbb R^d)$, being the image of a compact set under a continuous mapping and therefore, for almost every $t\in(t_0,\tau)$,
\[
\|\mathbf F_t(\mathbf x)\|_{\mathbf m} =\|f_t[\nu_{\mathbf x}]\|_{L^2(\nu_{\mathbf x})^d} 
\leq g_{\tau,K_R}(t) \quad \text{ for every } \mathbf x\text{ with }\|\mathbf x\|_{\mathbf m}\leq R.
\]
Thus $\mathbf F$ satisfies the Carath\'eodory conditions and the Carath\'eodory existence and  continuation theorems (see, e.g., \cite[Chapter~I, §5]{H_1980}) provide a locally absolutely continuous solution of \eqref{particle_system} on a maximal interval $[t_0,\tau_{\max})$, along which $\|\mathbf x(t)\|_{\mathbf m}\to+\infty$ as $t\nearrow\tau_{\max}$ if $\tau_{\max}<+\infty$.

\medskip\noindent
We claim that $\tau_{\max}=+\infty$. Let $\mathbf 0$ denote the configuration all of whose components vanish, so that $\nu_{\mathbf 0}=\delta_0$ and $\|\mathbf F_t(\mathbf 0)\|_{\mathbf m}=\|f_t[\delta_0]\|_{L^2(\delta_0)^d}=:h(t)$, which belongs to $L^1(t_0,\tau)$ for every $\tau>t_0$ by the fourth hypothesis applied with $K=\{\delta_0\}$.
Since $\sum_{i=1}^Nm_i\delta_{(x_i,0)}\in\Gamma(\nu_{\mathbf x},\delta_0)$, inequality \eqref{general_lambda_monotonicity} reads $\langle\mathbf F_t(\mathbf x)-\mathbf F_t(\mathbf 0),\mathbf x\rangle_{\mathbf m} \geq\lambda(t)\|\mathbf x\|^2_{\mathbf m}$. 
The function $u(t):=\|\mathbf x(t)\|_{\mathbf m}$ is locally absolutely continuous and at almost every $t$ at which $u(t)>0$ the Cauchy--Schwarz inequality gives
\[
\dot u(t) = -\frac{\langle\mathbf F_t(\mathbf x(t)),\mathbf x(t)\rangle_{\mathbf m}}{u(t)} \leq -\lambda(t)u(t)+h(t);
\]
the same inequality holds at almost every $t$ at which $u(t)= 0$, because $\dot u$ vanishes almost everywhere on the level set $\{u= 0\}$ and $h\geq 0$.
Gronwall's lemma (see, e.g., \cite[Chapter~1]{D_2003}) then gives
\[
u(t) \leq \exp\left(-\int_{t_0}^t\lambda(s)\,ds\right) \left(u(t_0)+\int_{t_0}^t\exp\left(\int_{t_0}^s\lambda(r)\,dr\right)h(s)\,ds\right) \quad \text{ for every } t\in[t_0,\tau_{\max}),
\]
so that $u$ is bounded on bounded intervals by a quantity depending only on $\lambda$, $h$ and $\mathcal M_2(\nu_{\mathbf x^0})$. 
Hence the solution cannot blow up in finite time and $\tau_{\max}=+\infty$.

\medskip\noindent

Set $\mu_t:=\nu_{\mathbf x(t)}$.
Given $\tau>t_0$, the image of $[t_0,\tau]$ under the continuous mapping $t\longmapsto\mu_t$ is a compact set $K\subseteq\mathcal P_2(\mathbb R^d)$ and therefore $\|f_t[\mu_t]\|_{L^2(\mu_t)^d}=\|\dot{\mathbf x}(t)\|_{\mathbf m}\leq g_{\tau,K}(t)$ for almost every $t\in(t_0,\tau)$, which gives the required integrability. 
Moreover, $\mathcal W_2(\mu_s,\mu_t) \leq \|\mathbf x(t)-\mathbf x(s)\|_{\mathbf m} \leq \mathlarger{\int}_s^t\|\dot{\mathbf x}(r)\|_{\mathbf m}\,dr$ for all $t_0\leq s\leq t$, so $\mu$ is locally absolutely continuous.
Finally, for every $\varphi\in C_c^\infty(\mathbb R^d)$ the function $t\longmapsto \mathlarger{\int} \varphi\,d\mu_t = \sum_{i=1}^Nm_i\varphi(x_i(t))$ is locally absolutely continuous and
\[
\frac{d}{dt}\int_{\mathbb R^d}\varphi(x)\,d\mu_t(x)
=
\sum_{i=1}^Nm_i\left\langle\nabla\varphi(x_i(t)),\dot x_i(t)\right\rangle
=
-\int_{\mathbb R^d}
\left\langle\nabla\varphi(x),f_t[\mu_t](x)\right\rangle\,d\mu_t(x)
\]
for almost every $t$, hence in the sense of distributions on $(t_0,+\infty)$. 
Thus $\mu$ is a solution of \eqref{general_flow} with initial datum $\nu_{\mathbf x^0}$.

\medskip\noindent
\textbf{Step 3:} \emph{Passage to the limit.}
Let $\mu_0\in\mathcal P_2(\mathbb R^d)$. We first approximate it by finitely supported measures. 
For $n\in\mathbb N$, cover the ball of radius $n$ centered
at the origin by finitely many pairwise disjoint nonempty Borel sets of diameter at most $1/n$, choose a point in each of them and let $T_n:\mathbb R^d\to\mathbb R^d$ be the Borel mapping sending each of these sets to the corresponding point and the complement of the ball to the origin. 
Then $\mu^n_0:=(T_n)_\#\mu_0$ is finitely supported and, since
$(\id,T_n)_\#\mu_0\in\Gamma(\mu_0,\mu^n_0)$,
\[
\mathcal W_2^2(\mu^n_0,\mu_0) \leq \int_{\mathbb R^d}|T_n(x)-x|^2\,d\mu_0(x) \leq \frac{1}{n^2}+\int_{\{|x|\geq n\}}|x|^2\,d\mu_0(x) \xrightarrow{n \to + \infty} 0 .
\] 
Discarding the chosen points carrying zero mass and merging those that coincide, we may write $\mu^n_0=\nu_{\mathbf x^{0, n}}$ for a configuration $\mathbf x^{0, n}$ with pairwise distinct components and with strictly positive masses summing to one, so that Step 2 applies.
By Step 2, for every $n$ there is a solution $\mu^n$ of \eqref{general_flow} with initial datum $\mu^n_0$ and by Step 1,
\[
\mathcal W_2(\mu^n_t,\mu^k_t) \leq \exp\left(-\int_{t_0}^t\lambda(s)\,ds\right) \mathcal W_2(\mu^n_0,\mu^k_0)
\quad \text{ for every } t\geq t_0.
\]
Since $\lambda$ is locally integrable, the exponential factor is bounded on $[t_0,\tau]$ for every $\tau>t_0$, so $(\mu^n)_n$ is a Cauchy sequence for the uniform distance on compact intervals; as
$(\mathcal P_2(\mathbb R^d),\mathcal W_2)$ is complete, there exists a continuous curve $\mu:[t_0,+\infty)\to\mathcal P_2(\mathbb R^d)$ such that $\mu^n_t\to\mu_t$ in $\mathcal W_2$, uniformly on compact intervals.
In particular $\mu_{t_0}=\mu_0$.

\medskip\noindent
Fix $\tau>t_0$ and let $K_\tau$ be the closure of $\{\mu^n_t:n\in\mathbb N,\ t\in[t_0,\tau]\}$. 
This set is compact: given $\varepsilon>0$, choose $n_0$ such that
$\sup_{t\in[t_0,\tau]}\mathcal W_2(\mu^n_t,\mu_t)<\varepsilon$ for every $n\geq n_0$; the images of $[t_0,\tau]$ under the finitely many continuous curves $\mu^1,\cdots,\mu^{n_0-1},\mu$ form a compact, hence totally bounded, set whose $\varepsilon$-neighbourhood contains $\{\mu^n_t:n\in\mathbb N,\ t\in[t_0,\tau]\}$, so that the latter is totally bounded and its closure is compact.
Note that $\mu_t\in K_\tau$ for every $t\in[t_0,\tau]$.
Writing $g_\tau:=g_{\tau,K_\tau}$, the fourth hypothesis gives
\[
\|f_t[\mu^n_t]\|_{L^2(\mu^n_t)^d}\leq g_\tau(t) \qquad\text{and}\qquad \|f_t[\mu_t]\|_{L^2(\mu_t)^d}\leq g_\tau(t)
\]
for almost every $t\in(t_0,\tau)$ and every $n\in\mathbb N$; in particular $\mathlarger{\int}_{t_0}^\tau\|f_t[\mu_t]\|_{L^2(\mu_t)^d}\,dt<+\infty$, the integrand being Borel measurable by  Lemma~\ref{borel_representative}.
Moreover, by Step 2 we have $\mathcal W_2(\mu^n_s,\mu^n_t)\leq \mathlarger{\int}_s^tg_\tau(r)\,dr$ for all $t_0\leq s\leq t\leq\tau$ and letting $n\to+\infty$ we conclude that $\mu$ is locally absolutely continuous.

\medskip\noindent
It remains to pass to the limit in the weak formulation. 
Let $t\in(t_0,\tau)$ be such that the preceding bounds hold and set $\Lambda^n_t:=(\id,f_t[\mu^n_t])_\#\mu^n_t$ and $\Lambda_t:=(\id,f_t[\mu_t])_\#\mu_t$; since $\mu^n_t\to\mu_t$ in $\mathcal W_2$, demicontinuity gives $\Lambda^n_t \rightharpoonup \Lambda_t$ narrowly. 
Fix $\varphi\in C_c^\infty(\mathbb R^d)$ and, for $R>0$, let $\theta_R:\mathbb R^d\to\mathbb R^d$ be given by $\theta_R(z) := z$ if $|z|\leq R$ and $\theta_R(z):=Rz/|z|$ otherwise.
The function  $(x,z)\longmapsto\langle\nabla\varphi(x),\theta_R(z)\rangle$ is bounded and continuous, while, for $\sigma\in\{\Lambda^n_t, \Lambda_t\}$,
\[
\int_{\mathbb R^{2d}} \left|\left\langle\nabla\varphi(x),z-\theta_R(z)\right\rangle\right|\,d\sigma(x,z)
\leq
\|\nabla\varphi\|_\infty\int_{\{|z|>R\}}|z|\,d\sigma(x,z)
\leq
\frac{\|\nabla\varphi\|_\infty}{R}\,g_\tau^2(t),
\]
because  $\mathlarger{\int} |z|^2\,d\Lambda^n_t=\|f_t[\mu^n_t]\|^2_{L^2(\mu^n_t)^d}\leq g_\tau^2(t)$ and similarly for $\Lambda_t$.
Letting first $n \to +\infty$ and then $R \to +\infty$, we obtain
\[
\int_{\mathbb R^d} \left\langle\nabla\varphi,f_t[\mu^n_t]\right\rangle\,d\mu^n_t \longrightarrow \int_{\mathbb R^d} \left\langle\nabla\varphi,f_t[\mu_t]\right\rangle\,d\mu_t \quad \text{ as } n \longrightarrow + \infty.
\]
These integrals are bounded in modulus by $\|\nabla\varphi\|_\infty g_\tau(t)$, which is integrable on $(t_0,\tau)$ and $\left| \mathlarger{\int} \varphi\,d\mu^n_t \right|\leq\|\varphi\|_\infty$. 
Hence, given $\alpha\in C_c^\infty((t_0,+\infty))$ and assuming, as we may, that $\tau>\sup\mathrm{supp}\alpha$ and writing the identity satisfied by $\mu^n$ as
\[
-\int_{t_0}^{+\infty}\dot\alpha(t)\int_{\mathbb R^d}\varphi\,d\mu^n_t\,dt
=
-\int_{t_0}^{+\infty}\alpha(t)
\int_{\mathbb R^d}\left\langle\nabla\varphi,f_t[\mu^n_t]\right\rangle\,d\mu^n_t\,dt,
\]
and letting $n\to+\infty$ with the aid of the dominated convergence theorem, we conclude that $\mu$ satisfies the same identity. 
Thus $\mu$ is a solution of \eqref{general_flow} with initial datum $\mu_0$ and the proof is complete.
\end{proof}

\noindent

We close this subsection with three remarks concerning Theorem~\ref{general_existence}.

\begin{remark}
The proof of Theorem~\ref{general_existence} is reminiscent of Dobrushin's approach to mean-field limits \cite[Section~6]{D_1979}, where stability in a transport distance yields convergence of empirical measures. Unlike Dobrushin, who establishes existence by a fixed-point argument, we construct solutions directly from particle approximations, using monotonicity rather than Lipschitz continuity to obtain stability.
\end{remark}

\begin{remark}
Although $\mathcal P_2(\mathbb R^d)$ has nonnegative curvature in the sense of Alexandrov, the stability estimate \eqref{general_contraction} retains the form of its classical Hilbert-space counterpart. The proof recovers the underlying Hilbertian structure by lifting the dynamics to a space of coupled trajectories.
\end{remark}

\begin{remark}
The estimate \eqref{general_contraction} is sharp even for genuinely measure-dependent operators. For instance, a skew-symmetric rotation of the barycenter yields, when $\lambda\equiv0$, a group of Wasserstein isometries that are not induced by Euclidean isometries.
\end{remark}

%%%%%%%%%%%%%%%%%%%%%%%%%%%%
%%%%%%%%%%%%%%%%%%%%%%%%%%%%

%%%%%%%%%%%%%%%%%%%%%%%%%%%%
%%%%%%%%%%%%%%%%%%%%%%%%%%%%

\section{Proofs of Section \ref{s2}}\label{pruebas}
\subsection{Proof of Proposition \ref{p1}}
We associate with $A$ the deterministic probability vector field
\[
\boldsymbol F[\mu]:=(\id,-A_\mu)_\#\mu \quad \text{ for }
\mu\in\mathcal P_2(\mathbb R^d),
\]
and we claim that it is maximal totally dissipative in the sense of \cite[Definition~3.6]{CSS_2026}.
The monotonicity condition \eqref{monotone_prop} implies the dissipativity of $\boldsymbol F$.
Moreover, $\boldsymbol F$ is demicontinuous in the sense of
\cite[Definition~3.22]{CSS_2026}.
Indeed, if $\mu_n\to\mu$ in $\mathcal W_2$, then the measures $(\id,-A_{\mu_n})_\#\mu_n$ converge narrowly to $(\id,-A_\mu)_\#\mu$, their second moments in the first variable converge, since $\mathcal M_2(\mu_n)\to\mathcal M_2(\mu)$ and their second moments in the second variable are uniformly bounded, being equal to $\|A_{\mu_n}\|^2_{L^2(\mu_n)^d}$; see
\cite[Proposition~2.5]{CSS_2026}. 
Hence \cite[Theorem~7.2]{CSS_2026} applies and yields
the claim. 
For a deterministic probability vector field, total dissipativity amounts to dissipativity along every coupling; see \cite[Remark~3.7]{CSS_2026}.
Since $\boldsymbol F$ is obtained from $A$ by a change of sign, this gives the first assertion. 
Moreover, by \cite[Theorem~3.12]{CSS_2026}, the Lagrangian representation of $\boldsymbol F$ in the sense of~\cite[Definition 3.8]{CSS_2026} is maximally dissipative and the usual maximality criterion in the Lagrangian space yields the second assertion.

\medskip\noindent

It remains to prove the third one. 
Let $\mu_0,\mu_1\in\mathrm{zer}(A)$, let $\gamma\in\Gamma(\mu_0,\mu_1)$, let $t\in[0,1]$ and write $\mu^\gamma_t=\bigl((1-t)\mathrm{proj}_1+t\mathrm{proj}_2\bigr)_\#\gamma$. By the second assertion, applied with $\nu:=\mu^\gamma_t$ and $\xi:= 0$, it suffices to show that 
\begin{equation}\label{zer_test}
\int_{\mathbb R^{2d}}\left\langle A_\mu(x),x-z\right\rangle\,d\sigma(x,z)\geq 0
\end{equation}
for every $\mu\in\mathcal P_2(\mathbb R^d)$ and every
$\sigma\in\Gamma(\mu,\mu^\gamma_t)$. 
Fix such a $\mu$ and such a $\sigma$ and disintegrate $\sigma$ and $\gamma$ with respect to $\mu^\gamma_t$, i.e., write $\sigma  = \mathlarger{\int}_{\mathbb R^d}\sigma_z\otimes\delta_z\,d\mu^\gamma_t(z)$ and $\gamma = \mathlarger{\int}_{\mathbb R^d}\gamma_z\,d\mu^\gamma_t(z)$, where $\sigma_z$ is a probability measure on $\mathbb R^d$ and $\gamma_z$ is concentrated on the set of pairs $(y_0,y_1)$ with $(1-t)y_0+ty_1=z$, for $\mu^\gamma_t$-a.e. $z$.
Then
\[
\alpha:=\int_{\mathbb R^d}\sigma_z\otimes\gamma_z\,d\mu^\gamma_t(z) \in\mathcal P(\mathbb R^d\times\mathbb R^{2d})
\]
is a probability measure with finite second moment, whose marginal in the variable
$x$ is $\mu$, whose marginal in the pair $(y_0,y_1)$ is $\gamma$ and whose push-forward through $(x,y_0,y_1)\longmapsto(x,(1-t)y_0+ty_1)$ is $\sigma$.
Denoting by $\sigma^i$ the marginal of $\alpha$ in the pair $(x,y_i)$, we have
$\sigma^i\in\Gamma(\mu,\mu_i)$ for $i= 0,1$ and, since $A_{\mu_i}= 0$
$\mu_i$-a.e., the first assertion gives
\[
\int\left\langle A_\mu(x),x-y_i\right\rangle\,d\alpha(x,y_0,y_1)
=
\int_{\mathbb R^{2d}} \left\langle A_\mu(x)-A_{\mu_i}(y),x-y\right\rangle\,d\sigma^i(x,y) \geq 0, \quad \text{ for  } i \in \{0, 1\}.
\]
Multiplying the inequality for $i= 0$ by $1-t$, the one for $i=1$ by $t$ and adding, we obtain
\[
\int\left\langle A_\mu(x),x-\bigl((1-t)y_0+ty_1\bigr)\right\rangle \,d\alpha(x,y_0,y_1)\geq 0,
\]
which is \eqref{zer_test}. Hence $A_{\mu^\gamma_t}= 0$ $\mu^\gamma_t$-a.e., so $\mu^\gamma_t\in\mathrm{zer}(A)$.

\subsection{Proof of Lemma \ref{l1}}
We apply Theorem~\ref{general_existence} to the time-independent family $f_t[\nu]:=A_\nu$, for which \eqref{general_flow} is precisely \eqref{eqtflw}.
The first two hypotheses are immediate and \eqref{general_lambda_monotonicity} holds with
$\lambda\equiv0$ by the first assertion of Proposition~\ref{p1}. As for the last one, demicontinuity implies that $\nu\longmapsto\|A_\nu\|_{L^2(\nu)^d}$ is bounded on
compact subsets of $\mathcal P_2(\mathbb R^d)$, so that $g_{\tau,K}$ may be taken constant. 
The conclusion, together with \eqref{general_contraction} for $\lambda\equiv0$, is the statement of the lemma.

%%%%%%%%%%%%%%%%%%%%%%%%%%%%%%%%%%%%%%%%%%%%%%%%%%%%%%%%%%%%%%%%%%%%%%%%%%%%%%%%%%%%%%%%%%%%%%%%%%%%%%%%%%%%

\subsection{Proof of Theorem~\ref{thm1}}
We divide the proof into four steps.

\medskip\noindent
\textbf{Step 1:} \emph{Demicontinuity and the energy inequality.}
We first observe that cocoercivity implies demicontinuity. 
Let $\mu,\nu\in\mathcal P_2(\mathbb R^d)$ and let $\gamma\in\Gamma_o(\mu,\nu)$. 
By cocoercivity and the Cauchy--Schwarz inequality,
\begin{align}
\beta \int_{\mathbb R^{2d}} |A_\mu(x)-A_\nu(y)|^2\,d\gamma(x,y)
&\leq 
\int_{\mathbb R^{2d}} \left\langle A_\mu(x)-A_\nu(y),x-y\right\rangle\,d\gamma(x,y)\\
&\leq
\left( \int_{\mathbb R^{2d}} |A_\mu(x)-A_\nu(y)|^2\,d\gamma(x,y) \right)^{1/2} \mathcal W_2(\mu,\nu).
\end{align}
Therefore,
\[
\left( \int_{\mathbb R^{2d}} |A_\mu(x)-A_\nu(y)|^2\,d\gamma \right)^{1/2}
\leq\beta^{-1}\mathcal W_2(\mu,\nu).
\]
Here, the case in which the integral vanishes is trivial.
Since
\[
\bigl( (\mathrm{proj}_1,A_\mu\circ\mathrm{proj}_1), (\mathrm{proj}_2,A_\nu\circ\mathrm{proj}_2) \bigr)_\#\gamma
\in 
\Gamma\bigl( (\id,A_\mu)_\#\mu, (\id,A_\nu)_\#\nu \bigr),
\]
it follows from $\gamma \in \Gamma_o(\mu, \nu)$ and inequality \eqref{plan_estimate} that
\[
\mathcal W_2^2\bigl( (\id,A_\mu)_\#\mu, (\id,A_\nu)_\#\nu \bigr)
\leq
\int_{\mathbb R^{2d}} \left( |x-y|^2+|A_\mu(x)-A_\nu(y)|^2
\right)d\gamma(x,y) 
\leq
\left(1+\frac1{\beta^2}\right) \mathcal W_2^2(\mu,\nu).
\]
In particular, if $\mu_n\to\mu$ in $\mathcal W_2$, then $(\id,A_{\mu_n})_\#\mu_n\to(\id,A_\mu)_\#\mu$ in $\mathcal W_2$, hence narrowly and with convergence of the second moments.
Since $\mathcal M_2\bigl((\id,A_{\mu_n})_\#\mu_n\bigr) =\mathcal M_2(\mu_n)+\|A_{\mu_n}\|_{L^2(\mu_n)^d}^2$ and $\mathcal M_2(\mu_n)\to\mathcal M_2(\mu)$, the sequence  $\big( \|A_{\mu_n}\|_{L^2(\mu_n)^d} \big)$ is bounded.
Thus $A$ is demicontinuous. 
Since cocoercivity implies monotonicity and $A$ is demicontinuous, the first assertion of Proposition~\ref{p1} shows that the monotonicity inequality holds along every coupling. Hence Lemma~\ref{l1} applies and \eqref{eqtflw} has a unique solution $\mu$.

\medskip\noindent
Fix now $\nu\in\mathrm{zer}(A)$ and let $\rho$ be a solution of
\eqref{eqtflw} with arbitrary initial datum $\rho_{t_0} \in \mathcal{P}_2(\R^d)$.
Set $\varphi_\nu(t):=\frac12\mathcal W_2^2(\rho_t,\nu)$, which is locally
absolutely continuous because $\rho$ is. 
Using a Borel representative of the velocity field provided by Lemma~\ref{borel_representative}, the differentiation formula for the squared Wasserstein distance \cite[Theorem~2.13(6a)]{CSS_2026} gives, for almost every $t>t_0$ and every $\gamma_t\in\Gamma_o(\rho_t,\nu)$,
\begin{align*}
\dot\varphi_\nu(t)
&\leq
-\int_{\mathbb R^{2d}} \left\langle A_{\rho_t}(x),x-y\right\rangle\,d\gamma_t(x,y)\\
& = -\int_{\mathbb R^{2d}} \left\langle A_{\rho_t}(x)-A_\nu(y),x-y\right\rangle \,d\gamma_t(x,y)\\
& \leq -\beta \int_{\mathbb R^{2d}} |A_{\rho_t}(x)-A_\nu(y)|^2\,d\gamma_t(x,y) = -\beta\|A_{\rho_t}\|_{L^2(\rho_t)^d}^2,
\end{align*}
where we have used that $A_\nu= 0$ $\nu$-a.e.\ and that the second marginal of $\gamma_t$ is $\nu$.
Integrating over $[t_0,t]$ yields
\begin{equation}\label{eq.energy_general}
\frac12\mathcal W_2^2(\rho_t,\nu) +  \beta\int_{t_0}^t \|A_{\rho_s}\|_{L^2(\rho_s)^d}^2\,ds  \leq \frac12\mathcal W_2^2(\rho_{t_0},\nu)
\end{equation}
for every $t\in[t_0,+\infty)$, every $\nu\in\mathrm{zer}(A)$ and every solution $\rho$ of \eqref{eqtflw}; the computation used only that $\rho$ solves \eqref{eqtflw}. 
Taking $\rho:=\mu$ gives \eqref{eq.cocoercive-energy}. 
In particular, $t\longmapsto\mathcal W_2(\mu_t,\nu)$ is nonincreasing and letting
$t\to+\infty$ in \eqref{eq.cocoercive-energy}, we get
\begin{equation}\label{eq.estimate_Wasserstein}
\int_{t_0}^{+\infty} \|A_{\mu_t}\|_{L^2(\mu_t)^d}^2\,dt \leq \frac1{2\beta}\mathcal W_2^2(\mu_0,\nu)<+\infty.
\end{equation}
Finally, if $\rho$ is a solution of \eqref{eqtflw} and $\eta$ is the superposition measure provided by Lemma~\ref{superposition}, applied to the
time-independent family $f_t[\nu] := A_\nu$, then \eqref{lift_estimate} and
the $L^2(\eta)^d$-valued equation give
\begin{equation}\label{eq.speed_bound}
\mathcal W_2(\rho_s,\rho_t) \leq \|e_t-e_s\|_{L^2(\eta)^d} \leq \int_s^t \|A_{\rho_r}\|_{L^2(\rho_r)^d}\,dr \quad \text{ for every } \,t_0\leq s\leq t,
\end{equation}
again for every solution of \eqref{eqtflw}.

\medskip

We prove that $t \longmapsto \|A_{\mu_t}\|_{L^2(\mu_t)^d}$ is nonincreasing. 
Let $\eta$ be the superposition measure introduced above, so that $(e_t)_\#\eta=\mu_t$ for every $t\ge t_0$ and
\[
\dot e_t+A_{\mu_t}(e_t) = 0 \quad\text{in }L^2(\eta)^d
\]
for almost every $t\in [t_0, +\infty)$. 
For $h>0$, define
\[
q_h(t) := \left\|\frac{e_{t+h}-e_t}{h}\right\|_{L^2(\eta)^d}^2,
\]
which is locally absolutely continuous.
For almost every $t\ge t_0$, we have
\begin{align*}
\dot q_h(t)
&=
\frac2{h^2}
\left\langle e_{t+h}-e_t,\dot e_{t+h}-\dot e_t
\right\rangle_{L^2(\eta)^d}
=
-\frac2{h^2}
\int
\left\langle
A_{\mu_{t+h}}(e_{t+h})-A_{\mu_t}(e_t),
e_{t+h}-e_t
\right\rangle\,d\eta
\le0,
\end{align*}
since $(e_{t+h},e_t)_\#\eta\in\Gamma(\mu_{t+h},\mu_t)$,
and the monotonicity inequality is valid along arbitrary couplings, according to Proposition~\ref{p1}. Thus $q_h$ is nonincreasing.

\medskip

Let $E\subseteq(t_0,+\infty)$ be a full-measure set on which
$t\longmapsto e_t$ is differentiable and
\[
\dot e_t + A_{\mu_t}(e_t) = 0
\]
holds at every $t \in E$. If $s,t\in E$ and $s<t$, then $q_h(t)\le q_h(s)$ for every $h>0$.
Letting $h\searrow 0$ gives
\[
\|\dot e_t\|_{L^2(\eta)^d}
\le
\|\dot e_s\|_{L^2(\eta)^d},
\]
that is,
\[
\|A_{\mu_t}(e_t)\|_{L^2(\eta)^d}
\le
\|A_{\mu_s}(e_s)\|_{L^2(\eta)^d}.
\]
Moreover, the estimate proved at the beginning of Step $1$,
\[
\left(
\int_{\mathbb R^{2d}}
|A_\mu(x)-A_\nu(y)|^2\,d\gamma(x,y)
\right)^{1/2}
\le
\beta^{-1} \mathcal W_2(\mu,\nu),
\quad \text{ for  } \gamma\in\Gamma_o(\mu,\nu),
\]
implies
\[
\big|
\|A_\mu\|_{L^2(\mu)^d}
-
\|A_\nu\|_{L^2(\nu)^d}
\big|
\le
\beta^{-1} \mathcal W_2(\mu,\nu).
\]

Hence, since $t\longmapsto\mu_t$ is continuous in the Wasserstein metric, $t \longmapsto \|A_{\mu_t}(e_t)\|_{L^2(\eta)^d}$ is continuous, thus the inequality $
\|A_{\mu_t}(e_t)\|_{L^2(\eta)^d}
\le
\|A_{\mu_s}(e_s)\|_{L^2(\eta)^d}$ extends
from $s,t\in E$ to every $t_0 \le s<t$, which proves that $t \longmapsto \|A_{\mu_t}\|_{L^2(\mu_t)^d}$ is nonincreasing.

\medskip

Now the energy estimate \eqref{eq.estimate_Wasserstein} gives, for every
$\nu \in \mathrm{zer}(A)$,
\[
\int_{t_0}^{+\infty} \|A_{\mu_s}(e_s)\|_{L^2(\eta)^d}^2\,ds
\le
\frac1{2 \beta} \mathcal W_2^2(\mu_0,\nu)<+\infty.
\]
Thus, for every $t \in [t_0, +\infty)$, we have
\[
\frac{t-t_0}{2} \|A_{\mu_t}(e_t)\|_{L^2(\eta)^d}^2
\le
\int_{\frac{t_0 + t}{2}}^{t}\|A_{\mu_s}(e_s)\|_{L^2(\eta)^d}^2 ds
\le
\int_{\frac{t_0 + t}{2}}^{+\infty}\|A_{\mu_s}(e_s)\|_{L^2(\eta)^d}^2 ds.
\]
Since $\int_{\frac{t_0 + t}{2}}^{+\infty}\|A_{\mu_s}(e_s)\|_{L^2(\eta)^d}^2 ds \longrightarrow 0 \ \mbox{as} \ t \longrightarrow +\infty$,
\[
\|A_{\mu_t}\|_{L^2(\mu_t)^d}
=o(t^{-1/2})
\quad\text{as }t \longrightarrow +\infty.
\]

Finally, estimate \eqref{eq.speed_bound} gives
\[
\mathcal W_2(\mu_s,\mu_t)
\le
\int_s^t \|A_{\mu_\tau}(e_\tau)\|_{L^2(\eta)^d} d\tau,
\]
hence, at every time at which the metric derivative exists,
\[
|\mu'|(t)\le \|A_{\mu_t}(e_t)\|_{L^2(\eta)^d} = o(t^{-1/2}) \quad\text{as }t \longrightarrow +\infty.
\]

\medskip\noindent
\textbf{Step 2:} \emph{Narrow cluster points.}
Fix $\nu\in\mathrm{zer}(A)$.
The Fej\'er monotonicity obtained in Step 1 yields
\[
\sqrt{\mathcal M_2(\mu_t)} = \mathcal W_2(\mu_t,\delta_0) \leq \mathcal W_2(\mu_t,\nu) + \sqrt{\mathcal M_2(\nu)} \leq \mathcal W_2(\mu_0,\nu) +\sqrt{\mathcal M_2(\nu)} \quad \text{ for every } t\geq t_0.
\]
Hence the orbit $\{\mu_t:t\geq t_0\}$ is bounded in $\mathcal P_2(\mathbb R^d)$ and consequently, relatively compact with respect to the narrow topology.

\medskip\noindent
We next record the restart property of the flow. 
Let $\rho$ be a solution of \eqref{eqtflw}, let $h\geq 0$ and set $\sigma_r:=\rho_{r+h}$ for $r\in[t_0,+\infty)$. 
The curve $\sigma$ is locally absolutely continuous, it satisfies 
\[
\int_{t_0}^\tau\|A_{\sigma_r}\|_{L^2(\sigma_r)^d}\,dr =\int_{t_0+h}^{\tau+h}\|A_{\rho_t}\|_{L^2(\rho_t)^d}\,dt<+\infty \quad \text{ for every $\tau>t_0$} 
\]
and, the velocity field of \eqref{eqtflw} being independent of time, the change of variables $t=r+h$ turns the distributional identity satisfied by $\rho$ on $(t_0+h,+\infty)$ into the same identity for $\sigma$ on $(t_0,+\infty)$. 
Hence $\sigma$ solves \eqref{eqtflw} with datum $\sigma_{t_0}=\rho_{t_0+h}$.
For $\rho_0\in\mathcal P_2(\mathbb R^d)$ and $t\geq 0$ we may therefore set $S_t\rho_0:=\sigma_{t_0+t}$, where $\sigma$ is the unique solution of \eqref{eqtflw} with $\sigma_{t_0}=\rho_0$ and the uniqueness assertion of Theorem~\ref{general_existence} gives  $S_{t+s}=S_t\circ S_s$ together with 
\[
\mathcal W_2(S_t\rho,S_t\sigma) \leq  \mathcal W_2(\rho,\sigma)
\quad  \text{ for every } \rho,\sigma\in\mathcal P_2(\mathbb R^d),\quad t\geq 0.
\]
Applying the shift to $\rho:=\mu$ with $h:=t-t_0$, where $t\geq t_0$, we
obtain a solution with datum $\mu_t$, whence $S_h\mu_t=\mu_{t+h}$ for all $t\geq t_0$ and $h\geq 0$; in particular $\mu_t=S_{t-t_0}\mu_0$. Moreover, by \eqref{eq.speed_bound}, the Cauchy--Schwarz inequality and the square integrability established in Step 1,
\[
\mathcal W_2(\mu_{t+1},\mu_t) \leq
\left( \int_t^{t+1} \|A_{\mu_r}\|_{L^2(\mu_r)^d}^2\,dr \right)^{1/2} \longrightarrow 0 \quad\text{ as } \quad t \longrightarrow +\infty.
\]
\smallbreak\noindent
Let now $t_n\to+\infty$ be such that
$\mu_{t_n}\rightharpoonup\bar\mu$ narrowly. Since the orbit is
bounded in $\mathcal P_2(\mathbb R^d)$, the lower semicontinuity of
the second moment with respect to narrow convergence gives
$\bar\mu\in\mathcal P_2(\mathbb R^d)$. The sequence $(\mu_{t_n})_n$ being
bounded in $\mathcal P_2(\mathbb R^d)$ and narrowly convergent to
$\bar\mu$, the Opial inequality \cite[Corollary~5.3]{NS_2021} applies to it
and yields
\begin{equation}\label{eq.opial}
\mathcal W_2^2(\bar\mu,\sigma)
+
\liminf_{n\to+\infty}
\mathcal W_2^2(\mu_{t_n},\bar\mu)
\leq
\liminf_{n\to+\infty}
\mathcal W_2^2(\mu_{t_n},\sigma)
\quad \text{ for every } \sigma\in\mathcal P_2(\mathbb R^d).
\end{equation}
We claim that $\bar\mu\in\mathrm{zer}(A)$. Since $S_1\mu_{t_n}=\mu_{t_n+1}$,
the nonexpansivity of $S_1$ gives
\[
\mathcal W_2(\mu_{t_n},S_1\bar\mu)
\leq
\mathcal W_2(\mu_{t_n},\mu_{t_n+1})
+
\mathcal W_2(\mu_{t_n+1},S_1\bar\mu)
\leq
\mathcal W_2(\mu_{t_n},\mu_{t_n+1})
+
\mathcal W_2(\mu_{t_n},\bar\mu).
\]
The first term tends to zero and the second one is bounded, whence
\[
\liminf_{n\to+\infty}
\mathcal W_2^2(\mu_{t_n},S_1\bar\mu)
\leq
\liminf_{n\to+\infty}
\mathcal W_2^2(\mu_{t_n},\bar\mu)<+\infty,
\]
and \eqref{eq.opial}, applied with $\sigma:=S_1\bar\mu$, gives
$S_1\bar\mu=\bar\mu$.
\smallbreak\noindent
Let $\rho$ be the solution of \eqref{eqtflw} with
$\rho_{t_0}=\bar\mu$. Since $\rho_{t_0+1}=S_1\bar\mu=\bar\mu$,
inequality \eqref{eq.energy_general}, applied to $\rho$ with any
$\nu\in\mathrm{zer}(A)$, gives
\[
\beta\int_{t_0}^{t_0+1} \|A_{\rho_s}\|_{L^2(\rho_s)^d}^2\,ds
\leq
\frac12\mathcal W_2^2(\bar\mu,\nu) - \frac12\mathcal W_2^2(\rho_{t_0+1},\nu) = 0.
\]
Thus $\|A_{\rho_s}\|_{L^2(\rho_s)^d}= 0$ for almost every $s\in(t_0,t_0+1)$ and \eqref{eq.speed_bound}, applied to $\rho$, yields
\[
\mathcal W_2(\bar\mu,\rho_s) \leq\int_{t_0}^s\|A_{\rho_r}\|_{L^2(\rho_r)^d}\,dr= 0
\]
for every $s\in[t_0,t_0+1]$. 
Hence $\rho$ is constant on $[t_0,t_0+1]$,
and choosing $s\in(t_0,t_0+1)$ with $A_{\rho_s}= 0$ $\rho_s$-a.e.\ we
obtain $\bar\mu=\rho_s\in\mathrm{zer}(A)$. Every narrow cluster point of
the trajectory therefore belongs to $\mathrm{zer}(A)$.

\medskip\noindent
\textbf{Step 3:} \emph{Uniqueness of the cluster point.}
For every $\sigma\in\mathrm{zer}(A)$, the Fej\'er monotonicity of Step 1 shows that the limit
\[
L(\sigma) := \lim_{t\to+\infty}\mathcal W_2^2(\mu_t,\sigma)
\]
exists. 
Let $\bar\mu,\widehat\mu$ be two narrow cluster points of the trajectory, which belong to $\mathrm{zer}(A)$ by Step 2. 
Along a sequence $t_n\to+\infty$ such that $\mu_{t_n}\rightharpoonup\bar\mu$, inequality \eqref{eq.opial} with $\sigma:=\widehat\mu$ gives
\[
\mathcal W_2^2(\bar\mu,\widehat\mu) + L(\bar\mu) \leq L(\widehat\mu).
\]
Interchanging $\bar\mu$ and $\widehat\mu$ we obtain $\mathcal W_2^2(\bar\mu,\widehat\mu)+L(\widehat\mu)\leq L(\bar\mu)$.
Adding the two inequalities yields $2\mathcal W_2^2(\bar\mu,\widehat\mu)\leq 0$ and hence $\bar\mu=\widehat\mu$.
Since the orbit is narrowly relatively compact  and the narrow topology is metrizable, it follows that there exists $\mu^\circ\in\mathrm{zer}(A)$ such that $\mu_t\rightharpoonup\mu^\circ$ narrowly as $t\to+\infty$.

\medskip\noindent
\textbf{Step 4:} \emph{Convergence in $\mathcal W_2$.}
Let $(\mu^n_0)_n$ be a sequence of finitely supported measures converging to $\mu_0$ in $\mathcal W_2$ and let $\mu^n$ be the solution of \eqref{eqtflw} with initial datum $\mu^n_0$.
By \eqref{general_contraction} with $\lambda\equiv0$,
\[
\sup_{t\geq t_0}\mathcal W_2(\mu^n_t,\mu_t)
\leq \mathcal W_2(\mu^n_0,\mu_0) \longrightarrow 0 \quad \text{ as } n\longrightarrow+\infty.
\]
Moreover, Step 2 of the proof of Theorem~\ref{general_existence} and uniqueness give
\[
\mu^n_t=\sum_{i=1}^{N_n}m_i^n\delta_{x_i^n(t)},
\qquad m_i^n>0,\qquad \sum_{i=1}^{N_n}m_i^n=1,
\]
where the masses are independent of time.

\medskip\noindent
Fix $\nu\in\mathrm{zer}(A)$.
Since the constant curve equal to $\nu$ is a solution of \eqref{eqtflw}, nonexpansivity and the triangle inequality give
\[
m_i^n|x_i^n(t)|^2 \leq \mathcal M_2(\mu^n_t)
\leq \left(\mathcal W_2(\mu^n_0,\nu)+\sqrt{\mathcal M_2(\nu)}\right)^2
\quad \text{ for every } t\geq t_0,\quad i \in \{ 1,\cdots ,N_n \}.
\]
For fixed $n$, the masses $m_i^n$ being positive and independent of time, the configurations $(x_1^n(t),\cdots,x_{N_n}^n(t))$ are bounded in the finite-dimensional space $(\mathbb R^d)^{N_n}$.
By \eqref{plan_estimate}, the mapping $(x_1,\cdots,x_{N_n})\longmapsto\sum_{i=1}^{N_n}m_i^n\delta_{x_i}$ is continuous in $\mathcal W_2$, so the orbit $\{\mu^n_t:t\geq t_0\}$ is relatively compact in $\mathcal W_2$.

\medskip\noindent
Given $\varepsilon>0$, choose $n$ such that $\sup_{t\geq t_0}\mathcal W_2(\mu^n_t,\mu_t)<\varepsilon$.
A finite $\varepsilon$-net for $\{\mu^n_t:t\geq t_0\}$ is then a finite $2\varepsilon$-net for $\{\mu_t:t\geq t_0\}$.
Thus the latter is totally bounded and, by completeness of $(\mathcal P_2(\mathbb R^d),\mathcal W_2)$, relatively compact in $\mathcal W_2$.
Every sequence $t_n\to+\infty$ therefore admits a subsequence along which $\mu_{t_n}$ converges in $\mathcal W_2$, and the limit is necessarily $\mu^\circ$ by Step 3.
Hence $\mathcal W_2(\mu_t,\mu^\circ)\longrightarrow0$ as $t\longrightarrow+\infty$, which completes the proof.

%%%%%%%%%%%%%%%%%%%%%%%%%%%%%%%%%%%%%%%%%%%%%%%%%%%%%%%%%%%%%%%%%%%%%%%%%%%%%%%%%%%%%%%%%%%%%%%%%%%%%%%%%%%%%%%%%%%%%%%%%%%%%%%%%%%%%%%%%%

%%%%%%%%%%%%%%%%%%%%%%%%%%%%%%%%%%%%%%%%%%%%%%%%%%%%%%%%%%%%%%%%%%%%%%%%%%%%%%%%%%%%%%%%%%%%%%%%%%%%%%%%%%%%%%%%%%%%%%%%%%%%%%%%%%%%%%%%%%%%%%%%%%%%%%%%%

\subsection{Proof of Lemma \ref{l2}}
We apply Theorem~\ref{general_existence} to the family
$f_t[\nu]:=A_\nu+\varepsilon(t)\id$, for which \eqref{general_flow} is precisely
\eqref{eqtfl}. The mapping $t\longmapsto f_t[\nu]$ is Borel measurable because
$\varepsilon$ is and $\nu\longmapsto f_t[\nu]$ is demicontinuous, since
$\nu\longmapsto\id$ is. Adding $\varepsilon(t)\id$ to $A$ adds
$\varepsilon(t) \mathlarger{\int} |x-y|^2\,d\gamma$ to the left-hand side of the first assertion of
Proposition~\ref{p1}, so that \eqref{general_lambda_monotonicity} holds with $\lambda:=\varepsilon$, which is locally integrable. Finally, for every compact set
$K\subseteq\mathcal P_2(\mathbb R^d)$,
\[
\|f_t[\nu]\|_{L^2(\nu)^d}
\leq
\sup_{\sigma\in K}\|A_\sigma\|_{L^2(\sigma)^d}
+\varepsilon(t)\sup_{\sigma\in K}\mathcal M_2(\sigma)^{1/2}
\quad \text{ for every } \nu\in K,
\]
and the right-hand side belongs to $L^1(t_0,\tau)$ for every $\tau>t_0$, the two
suprema being finite by demicontinuity and by the compactness of $K$. 
The conclusion, together with \eqref{general_contraction} for $\lambda=\varepsilon$, is the statement of the lemma.

\subsection{Proof of Lemma \ref{lemsecmoment}}
Let $\nu\in\mathrm{zer}(A)$ and denote its barycenter by $m_\nu$.
For every $\mu\in\mathcal P_2(\mathbb R^d)$, the first assertion of Proposition~\ref{p1}, applied to $\mu\otimes\nu$, gives
\[
\int_{\mathbb R^d}\left\langle A_\mu(x),x-m_\nu\right\rangle\,d\mu(x)
=
\int_{\mathbb R^{2d}}\left\langle A_\mu(x)-A_\nu(y),x-y\right\rangle\,d(\mu\otimes\nu)(x,y)\geq 0,
\]
since $A_\nu= 0$ $\nu$-a.e.
The measure $\mu\otimes\delta_{m_\nu}$ being the only coupling between $\mu$ and $\delta_{m_\nu}$, the second assertion of Proposition~\ref{p1} yields $\delta_{m_\nu}\in\mathrm{zer}(A)$.
Thus $m_\nu\in K$ for every $\nu\in\mathrm{zer}(A)$ and in particular $K\neq\varnothing$.

\medskip\noindent
Let $x,y\in K$ and $t\in[0,1]$.
The third assertion of Proposition~\ref{p1}, applied to $\delta_x\otimes\delta_y$, gives $\delta_{(1-t)x+ty}\in\mathrm{zer}(A)$ and hence $(1-t)x+ty\in K$.
Thus $K$ is convex.
Let now $(x_n)_{n\in\mathbb N}\subseteq K$ with $x_n\to x$.
Since $\delta_{x_n}\to\delta_x$ in $\mathcal W_2$, demicontinuity gives
\[
\delta_{(x_n,0)}=(\id,A_{\delta_{x_n}})_\#\delta_{x_n}
\rightharpoonup
(\id,A_{\delta_x})_\#\delta_x
\]
narrowly.
The left-hand side converges narrowly to $\delta_{(x,0)}$, so that $A_{\delta_x}(x)= 0$.
Hence $x\in K$ and $K$ is closed.

\medskip\noindent
Set $z:=\mathrm{proj}_K(0)$.
For every $\nu\in\mathrm{zer}(A)$, we have $m_\nu\in K$ and therefore
\[
\mathcal M_2(\nu)=\mathrm{Var}(\nu)+|m_\nu|^2\geq|z|^2=\mathcal M_2(\delta_z).
\]
Since $\delta_z\in\mathrm{zer}(A)$, it is a minimizer of the second moment over $\mathrm{zer}(A)$.
If equality holds, then $\mathrm{Var}(\nu)= 0$ and $|m_\nu|=|z|$.
The uniqueness of the projection onto $K$ gives $m_\nu=z$ and hence $\nu=\delta_z$.
This proves the assertion.

\subsection{Proof of Theorem \ref{thm2}}
We divide the proof into five steps.

\medskip\noindent
\textbf{Step 1:} \emph{Superposition and boundedness of the trajectory.}
We apply Lemma~\ref{superposition} to the family $f_t[\nu]:=A_\nu+\varepsilon(t)\id$, whose hypotheses are satisfied, as observed in the proof of Lemma~\ref{l2}. 
Hence there exists $\eta\in\mathcal P\bigl(C([t_0,+\infty);\mathbb R^d)\bigr)$, concentrated on $\mathrm{AC}_{\mathrm{loc}}([t_0,+\infty);\mathbb R^d)$, such that
$(e_t)_\#\eta=\mu_t$ for every $t\in[t_0,+\infty)$, the mapping $[t_0,+\infty)\ni t\longmapsto e_t\in L^2(\eta)^d$ is locally absolutely continuous,
and
\begin{equation}\label{eq.thm2_characteristic}
\dot e_t+A_{\mu_t}(e_t)+\varepsilon(t)e_t= 0
\qquad\text{in }L^2(\eta)^d
\end{equation}
for almost every $t\in[t_0,+\infty)$.

\medskip\noindent
We first prove that the second moments along the trajectory are uniformly bounded. Fix
$\nu\in\mathrm{zer}(A)$ and $\gamma_0\in\Gamma_o(\mu_0,\nu)$. Since
$(e_{t_0})_\#\eta=\mu_0$, the disintegration theorem provides a Borel family
$(\eta_x)_{x\in\mathbb R^d}$ of probability measures on
$C([t_0,+\infty);\mathbb R^d)$ such that $\eta = \mathlarger{\int}_{\mathbb R^d}\eta_x\,d\mu_0(x)$
and $\eta_x$ is concentrated on the curves issuing from $x$ for $\mu_0$-a.e.
$x$; let $(\gamma_x)_{x\in\mathbb R^d}$ be the disintegration of $\gamma_0$ with
respect to its first marginal, so that
$\gamma_0 = \mathlarger{\int}_{\mathbb R^d}\delta_x\otimes\gamma_x\,d\mu_0(x)$.
Set
\[
\Theta:=\int_{\mathbb R^d}\eta_x\otimes\gamma_x\,d\mu_0(x)
\in\mathcal P\bigl(C([t_0,+\infty);\mathbb R^d)\times\mathbb R^d\bigr),
\]
and denote by $E_t$ the evaluation of the curve at time $t$ and by $Y$ the second
variable, i.e., $E_t(\omega,y):=\omega(t)$ and $Y(\omega,y):=y$. The marginals of
$\Theta$ are $\eta$ and $\nu$, so that $(E_t)_\#\Theta=\mu_t$ for every
$t\in[t_0,+\infty)$ and $Y_\#\Theta=\nu$; moreover, the measures $\eta_x$ being
concentrated on the curves issuing from $x$, we have
$(E_{t_0},Y)_\#\Theta=\gamma_0$. The pull-back through the projection onto the first
factor is a linear isometry from $L^2(\eta)^d$ into $L^2(\Theta)^d$ carrying $e_t$
to $E_t$; as bounded linear mappings commute with Bochner integration,
$t\longmapsto E_t\in L^2(\Theta)^d$ is locally absolutely continuous and
\eqref{eq.thm2_characteristic} gives
\[
\dot E_t+A_{\mu_t}(E_t)+\varepsilon(t)E_t= 0
\qquad\text{in }L^2(\Theta)^d,
\]
for almost every $t\in[t_0,+\infty)$.

\medskip\noindent
Set $q(t):=\|E_t-Y\|_{L^2(\Theta)^d}^2$, which is locally absolutely continuous by the
differentiation rule for Hilbert-valued curves.
Since $(E_t,Y)_\#\Theta\in\Gamma(\mu_t,\nu)$ and $A_\nu= 0$ $\nu$-a.e., the first assertion of Proposition~\ref{p1} gives $\mathlarger{\int} \left\langle A_{\mu_t}(E_t),E_t-Y\right\rangle\,d\Theta\geq 0$ and therefore, for almost every $t\in[t_0,+\infty)$,
\begin{align*}
\frac12\dot q(t)
&=
-\int\left\langle A_{\mu_t}(E_t)+\varepsilon(t)E_t,E_t-Y\right\rangle\,d\Theta \\
& \leq
-\varepsilon(t)\int\left\langle E_t,E_t-Y\right\rangle\,d\Theta = -\frac{\varepsilon(t)}2
\bigl(q(t)+\mathcal M_2(\mu_t)-\mathcal M_2(\nu)\bigr)
\leq
-\frac{\varepsilon(t)}2q(t)+\frac{\varepsilon(t)}2\mathcal M_2(\nu),
\end{align*}
where we have used the identity $2\langle a,a-b\rangle=|a-b|^2+|a|^2-|b|^2$ and the
fact that $\|E_t\|^2_{L^2(\Theta)^d}=\mathcal M_2(\mu_t)$ and
$\|Y\|^2_{L^2(\Theta)^d}=\mathcal M_2(\nu)$. Thus
$\dot q(t)+\varepsilon(t)q(t)\leq\varepsilon(t)\mathcal M_2(\nu)$ for almost every
$t\in[t_0,+\infty)$.
Since $\dot\gamma_\varepsilon=\varepsilon\gamma_\varepsilon$
almost everywhere and $\gamma_\varepsilon(t_0)=1$, multiplying by
$\gamma_\varepsilon$ gives
\[
\frac{d}{dt} \Bigl(\gamma_\varepsilon(t)\bigl(q(t)-\mathcal M_2(\nu)\bigr)\Bigr)\leq 0
\qquad\text{for a.e. }t\in[t_0,+\infty),
\]
whence
\[
q(t) \leq \mathcal M_2(\nu) + \frac{q(t_0)-\mathcal M_2(\nu)}{\gamma_\varepsilon(t)}
\quad \text{ for every } t\geq t_0.
\]
Since $\varepsilon>0$, the function $\gamma_\varepsilon$ is nondecreasing and
$\gamma_\varepsilon\geq1$, while the optimality of $\gamma_0$ gives
\[
    q(t_0)= \mathlarger{\int}_{\mathbb R^{2d}}|x-y|^2\,d\gamma_0(x,y)=\mathcal W_2^2(\mu_0,\nu).
\]
Discussing the sign of $q(t_0)-\mathcal M_2(\nu)$, we conclude that
\[
q(t)\leq Q:=\max\bigl\{\mathcal W_2^2(\mu_0,\nu),\mathcal M_2(\nu)\bigr\}
\quad \text{ for every } t\geq t_0,
\]
and therefore
$\|e_t\|_{L^2(\eta)^d}=\|E_t\|_{L^2(\Theta)^d}
\leq\|E_t-Y\|_{L^2(\Theta)^d}+\|Y\|_{L^2(\Theta)^d}
\leq\sqrt Q+\sqrt{\mathcal M_2(\nu)}=:D$. In other words,
\begin{equation}\label{eq.thm2_moment_bound}
\sqrt{\mathcal M_2(\mu_t)}=\|e_t\|_{L^2(\eta)^d}\leq D
\quad \text{ for every } t\geq t_0.
\end{equation}

\medskip\noindent
\textbf{Step 2:} \emph{Estimate of the velocity field.}
Fix $h>0$ and set
\[
q_h(t):=\left\|\frac{e_{t+h}-e_t}h\right\|_{L^2(\eta)^d}^2,
\]
which is locally absolutely continuous, both $t\longmapsto e_{t+h}$ and
$t\longmapsto e_t$ being locally absolutely continuous with values in
$L^2(\eta)^d$; moreover, its derivative is
$2h^{-2}\langle\dot e_{t+h}-\dot e_t,e_{t+h}-e_t\rangle_{L^2(\eta)^d}$ for almost
every $t$. Writing \eqref{eq.thm2_characteristic} at $t+h$ and at $t$ and
applying the first assertion of Proposition~\ref{p1} to the coupling
$(e_{t+h},e_t)_\#\eta\in\Gamma(\mu_{t+h},\mu_t)$, we obtain
\[
\frac12\frac{d}{dt}\|e_{t+h}-e_t\|_{L^2(\eta)^d}^2
\leq
-\int
\left\langle
e_{t+h}-e_t,\varepsilon(t+h)e_{t+h}-\varepsilon(t)e_t
\right\rangle\,d\eta
\]
for almost every $t\in[t_0,+\infty)$. Since
\[
2\left\langle a-b,\alpha a-\beta b\right\rangle
=
(\alpha+\beta)|a-b|^2+(\alpha-\beta)\bigl(|a|^2-|b|^2\bigr)
\quad \text{ for every } a,b\in\mathbb R^d \text{ and } \alpha,\beta\in\mathbb R,
\]
and since $\|e_r\|_{L^2(\eta)^d}^2=\mathcal M_2(\mu_r)$ for every $r\geq t_0$,
dividing by $h^2$ yields
\begin{equation}\label{eq.thm2_qh}
\dot q_h(t)+\bigl(\varepsilon(t+h)+\varepsilon(t)\bigr)q_h(t)
\leq
-\,\frac{\varepsilon(t+h)-\varepsilon(t)}h\,
\frac{\mathcal M_2(\mu_{t+h})-\mathcal M_2(\mu_t)}h
\end{equation}
for almost every $t\in[t_0,+\infty)$. On the other hand, by
\eqref{eq.thm2_moment_bound},
\[
\left|\frac{\mathcal M_2(\mu_{t+h})-\mathcal M_2(\mu_t)}h\right|
=
\left|
\left\langle\frac{e_{t+h}-e_t}h,e_{t+h}+e_t\right\rangle_{L^2(\eta)^d}
\right|
\leq
2D\sqrt{q_h(t)},
\]
so that \eqref{eq.thm2_qh} becomes
\begin{equation}\label{eq.thm2_qh2}
\dot q_h(t)+\bigl(\varepsilon(t+h)+\varepsilon(t)\bigr)q_h(t)
\leq
2D\left|\frac{\varepsilon(t+h)-\varepsilon(t)}h\right|\sqrt{q_h(t)}
\end{equation}
for almost every $t\in[t_0,+\infty)$.

\medskip\noindent
Set $r_h:=\sqrt{q_h}$, i.e.,
$r_h(t)=\|h^{-1}(e_{t+h}-e_t)\|_{L^2(\eta)^d}$; being the norm of a locally
absolutely continuous $L^2(\eta)^d$-valued curve, $r_h$ is locally absolutely
continuous. At almost every $t$ at which $r_h(t)>0$ we have
$\dot q_h(t)=2r_h(t)\dot r_h(t)$ and \eqref{eq.thm2_qh2} gives
\[
\dot r_h(t)+\frac{\varepsilon(t+h)+\varepsilon(t)}2r_h(t)
\leq
D\left|\frac{\varepsilon(t+h)-\varepsilon(t)}h\right|;
\]
the same inequality holds at almost every $t$ at which $r_h(t)= 0$, because the
derivative of a nonnegative absolutely continuous function vanishes almost everywhere on
its zero set. Thus it holds for almost every $t\in[t_0,+\infty)$. Fix $s>t_0$, set
\[
G_h(t):=\exp\left(\frac12\int_s^t\bigl(\varepsilon(\tau+h)+\varepsilon(\tau)\bigr)\,d\tau
\right),
\]
multiply the preceding inequality by $G_h$ and integrate over $[s,t]$, so as to
obtain
\begin{equation}\label{eq.thm2_sqrt_qh}
G_h(t)r_h(t)
\leq
r_h(s)
+
D\int_s^tG_h(\tau)\left|\frac{\varepsilon(\tau+h)-\varepsilon(\tau)}h\right|\,d\tau
\quad \text{ for every } t\geq s.
\end{equation}
\smallbreak\noindent
Let $E\subseteq(t_0,+\infty)$ be a set of full Lebesgue measure such that
$t\longmapsto e_t$ is differentiable and \eqref{eq.thm2_characteristic} holds at every
$t\in E$ and let $s,t\in E$ with $s<t$. Then $r_h(s)\to\|\dot e_s\|_{L^2(\eta)^d}$
and $r_h(t)\to\|\dot e_t\|_{L^2(\eta)^d}$ as $h\searrow 0$. Moreover, since
$\varepsilon$ is locally absolutely continuous,
\[
\frac{\varepsilon(\tau+h)-\varepsilon(\tau)}h=\frac1h\int_0^h\dot\varepsilon(\tau+r)\,dr
\]
for almost every $\tau$, whence, by the continuity of translations in $L^1$, the
difference quotients converge to $\dot\varepsilon$ in $L^1(s,t)$. Similarly
$\varepsilon(\cdot+h)\to\varepsilon$ in $L^1(s,t)$ and therefore
\[
G_h(\tau)\longrightarrow
\exp\left(\int_s^\tau\varepsilon(r)\,dr\right)
=
\frac{\gamma_\varepsilon(\tau)}{\gamma_\varepsilon(s)}
\]
uniformly for $\tau\in[s,t]$. We may thus let $h\searrow  0$ in
\eqref{eq.thm2_sqrt_qh} and multiply by $\gamma_\varepsilon(s)$, which gives
\begin{equation}\label{eq.thm2_velocity_energy}
\gamma_\varepsilon(t)\|\dot e_t\|_{L^2(\eta)^d}
\leq
\gamma_\varepsilon(s)\|\dot e_s\|_{L^2(\eta)^d}
+
D\int_s^t|\dot\varepsilon(\tau)|\gamma_\varepsilon(\tau)\,d\tau
\quad \text{ for every } s,t\in E \text{ with } s < t.
\end{equation}
\smallbreak\noindent
For $t\in E$, the identity $(e_t)_\#\eta=\mu_t$ and
\eqref{eq.thm2_characteristic} give
\begin{equation}\label{eq.thm2_B}
\|\dot e_t\|_{L^2(\eta)^d}
=
\|A_{\mu_t}+\varepsilon(t)\id\|_{L^2(\mu_t)^d}.
\end{equation}
Fix $s\in E$, set $C:=\gamma_\varepsilon(s)\|A_{\mu_s}+\varepsilon(s)\id\|_{L^2(\mu_s)^d}$,
and rewrite \eqref{eq.thm2_velocity_energy} as
\begin{equation}\label{eq.thm2_B_est_E}
\|A_{\mu_t}+\varepsilon(t)\id\|_{L^2(\mu_t)^d}
\leq
\frac C{\gamma_\varepsilon(t)}
+
\frac D{\gamma_\varepsilon(t)}\int_s^t|\dot\varepsilon(\tau)|\gamma_\varepsilon(\tau)\,d\tau
\quad \text{ for every } t\in E \text{ with } t\geq s.
\end{equation}
\smallbreak\noindent
We now extend this estimate to every $t\geq s$ and to this end we claim that the
mapping $t\longmapsto\|A_{\mu_t}+\varepsilon(t)\id\|_{L^2(\mu_t)^d}$ is lower
semicontinuous. Let $t_n\to t$. Since $\mu_{t_n}\to\mu_t$ in $\mathcal W_2$,
demicontinuity gives
\[
\lambda_n:=(\id,A_{\mu_{t_n}})_\#\mu_{t_n} \rightharpoonup \lambda:=(\id,A_{\mu_t})_\#\mu_t
\]
narrowly. 
Set $T_n(x,z):=(x,z+\varepsilon(t_n)x)$ and $T(x,z):=(x,z+\varepsilon(t)x)$. 
Since $\varepsilon(t_n)\to\varepsilon(t)$, the mappings $T_n$ converge to $T$ uniformly on compact subsets of $\mathbb R^{2d}$, and the narrow convergence of $(\lambda_n)_n$, which entails its tightness, yields
\[
(\id,A_{\mu_{t_n}}+\varepsilon(t_n)\id)_\#\mu_{t_n}
=
(T_n)_\#\lambda_n \rightharpoonup T_\#\lambda
=
(\id,A_{\mu_t}+\varepsilon(t)\id)_\#\mu_t
\]
narrowly. 
Recall that
\[
    \mathcal{M}_2((\id,A_{\mu_{t}}+\varepsilon(t)\id)_\#\mu_{t}) = \mathcal{M}_2(\mu_t) + \| A_{\mu_t} + \varepsilon(t) \id \|_{L^2(\mu_t)^d}^2.
\]
Also, since $\mu_{t_n} \longrightarrow \mu_t$ in $\mathcal{W}_2$, $\mathcal{M}_2(\mu_{t_n}) \longrightarrow \mathcal{M}_2(\mu_t)$.
Therefore, the lower semicontinuity of the second moment under narrow convergence yields
\[
    \| A_{\mu_t} + \varepsilon(t) \id \|_{ L^2(\mu_t)^d} \leq \liminf_{n \to \infty} \|A_{ \mu_{t_n} } + \varepsilon(t_n) \id \|_{L^2(\mu_{t_n} )^d}.
\]

\medskip\noindent
Since $E$ has full measure and hence is dense, and since
the right-hand side of \eqref{eq.thm2_B_est_E} is continuous in $t$, the estimate
\eqref{eq.thm2_B_est_E} holds for every $t\geq s$. Enlarging the domain of integration
from $[s,t]$ to $[t_0,t]$, we conclude that
\begin{equation}\label{eq.thm2_main_est}
\|A_{\mu_t}+\varepsilon(t)\id\|_{L^2(\mu_t)^d}
\leq
R(t) := \frac C{\gamma_\varepsilon(t)}
+
\frac D{\gamma_\varepsilon(t)} \int_{t_0}^t|\dot\varepsilon(\tau)|\gamma_\varepsilon(\tau)\,d\tau \quad \text{for every } t\geq s.
\end{equation}

\medskip\noindent
\textbf{Step 3:} \emph{The estimates for $|\mu'|$ and $A_{\mu_t}$.}
The function $R$ defined in \eqref{eq.thm2_main_est} is continuous. For
$s\leq r\leq t$, inequality \eqref{lift_estimate}, the local absolute continuity of
$t\longmapsto e_t$ and \eqref{eq.thm2_B} give
\[
\mathcal W_2(\mu_r,\mu_t)
\leq
\|e_t-e_r\|_{L^2(\eta)^d}
\leq
\int_r^t\|\dot e_\tau\|_{L^2(\eta)^d}\,d\tau
=
\int_r^t\|A_{\mu_\tau}+\varepsilon(\tau)\id\|_{L^2(\mu_\tau)^d}\,d\tau
\leq
\int_r^tR(\tau)\,d\tau,
\]
the equality in the third integral holding for almost every $\tau$. Hence
$\limsup_{h\to0}|h|^{-1}\mathcal W_2(\mu_{t+h},\mu_t)\leq R(t)$ for every $t>s$ and
in particular
\[
|\mu'|(t)
\leq
\frac C{\gamma_\varepsilon(t)}
+
\frac D{\gamma_\varepsilon(t)}
\int_{t_0}^t|\dot\varepsilon(\tau)|\gamma_\varepsilon(\tau)\,d\tau
\]
at every $t>s$ at which the metric derivative exists. Moreover, by
\eqref{eq.thm2_moment_bound} we have
$\|\id\|_{L^2(\mu_t)^d}=\sqrt{\mathcal M_2(\mu_t)}\leq D$ and therefore
\[
\|A_{\mu_t}\|_{L^2(\mu_t)^d}
\leq
\|A_{\mu_t}+\varepsilon(t)\id\|_{L^2(\mu_t)^d}
+
\varepsilon(t)\|\id\|_{L^2(\mu_t)^d}
\leq
\frac C{\gamma_\varepsilon(t)}
+
\frac D{\gamma_\varepsilon(t)}
\int_{t_0}^t|\dot\varepsilon(\tau)|\gamma_\varepsilon(\tau)\,d\tau
+
D\varepsilon(t)
\]
for every $t\geq s$. Since $C$ and $D$ do not depend on $t$, these are the two
asserted asymptotic estimates.

\medskip\noindent

\textbf{Step 4}: \emph{A second estimate for $\big| \mu' \big|$ and $A_{\mu_t}$.}
For the sake of brevity, denote $m(t) := \mathcal{M}_2(\mu_t)$ and $V(t) := \| A_{\mu_t} + \varepsilon(t) \id \|_{L^2(\mu_t)^d}$.
In what follows, we use the notations in Step 2.
Notice that $m$ is locally absolutely continuous and since $m(t) = \| e_t \|^2_{L^2(\eta)^d}$, we infer that 
\[
    \dot{m}(t)  = 2 \langle e_t, \dot{e}_t \rangle_{L^2(\eta)^d} \quad \text{for a.e. } t.
\]
Fix $s, t \in E$ and $t > s$.
Multiplying~\eqref{eq.thm2_qh} by $G_h(t)^2$, integrating with respect to $t$ and then letting $h \searrow 0$, we obtain, with $C := \gamma_\varepsilon(s)^2V(s)^2\geq  0$,
\begin{equation}\label{2Rate}
    \gamma_\varepsilon^2(t) V(t)^2  \leq C - \int_s^t \dot{\varepsilon}(\tau) \gamma_\varepsilon^2(\tau) \dot{m}(\tau) \,d\tau \quad \text{ for every } t \in E, \, t \geq s.
\end{equation}
We write
\begin{equation}\label{2Rate.01}
   - \int_s^t \dot{\varepsilon}(\tau) \gamma_\varepsilon^2(\tau) \dot{m}(\tau) \,d\tau =  \int_s^t \varepsilon^2(\tau) \gamma_\varepsilon^2(\tau) \dot{m}(\tau) \, d\tau  - \int_s^t (\dot\varepsilon(\tau) + \varepsilon^2(\tau)) \gamma_\varepsilon^2(\tau) \dot{m}(\tau) \, d\tau.
\end{equation}
Since $\dot{\gamma_\varepsilon} = \varepsilon \gamma_\varepsilon$, we note that
\[
    \dfrac{d}{d\tau} \big( \varepsilon^2 \gamma_\varepsilon^2 \big) = 2 \varepsilon \dot{\varepsilon} \gamma_\varepsilon^2 + 2 \varepsilon^3 \gamma_\varepsilon^2 = 2 \varepsilon (\dot{\varepsilon} + \varepsilon^2) \gamma_\varepsilon^2.
\]
Integrating by parts then yields
\begin{equation}\label{2Rate.02}
    \int_s^t \varepsilon^2(\tau) \gamma_\varepsilon^2(\tau) \dot{m}(\tau) \, d\tau  =  \varepsilon^2(t)\gamma_\varepsilon^2(t) m(t) - \varepsilon^2(s)\gamma_\varepsilon^2(s) m(s) - \int_s^t 2 \varepsilon(\tau) (\dot{\varepsilon}(\tau) + \varepsilon^2(\tau)) \gamma_\varepsilon^2(\tau) m(\tau)\, d\tau.
\end{equation}
Combining~\eqref{2Rate},~\eqref{2Rate.01} and~\eqref{2Rate.02}, we infer that
\begin{align*}
    \gamma_\varepsilon^2(t) V(t)^2 \leq C + \varepsilon^2(t) \gamma_\varepsilon^2(t)m(t) - \int_s^t (\dot{\varepsilon}(\tau) + \varepsilon^2(\tau)) \gamma_\varepsilon^2(\tau) \big( \dot{m}(\tau) + 2 \varepsilon(\tau) m(\tau) \big) \, d\tau.
\end{align*}
Set $\theta(t) := \dot{\varepsilon}(t) + \varepsilon^2(t)$.
Thanks to Step 1, there is a constant $D > 0$ such that $m(\tau) \leq D^2$ and $ |\dot{m}(\tau)| \leq 2 D V(\tau)$.
Therefore,
\begin{equation}\label{2Rate.03}
    \gamma_\varepsilon^2(t) V(t)^2 \leq C + \varepsilon^2(t) \gamma_\varepsilon^2(t)m(t) + 2D \int_s^t |\theta(\tau)| \gamma_\varepsilon^2(\tau) V(\tau) \, d\tau + 2 D^2 \int_s^t |\theta(\tau)| \varepsilon(\tau) \gamma_\varepsilon^2(\tau) \, d\tau.
\end{equation}
Further, since $\dfrac{d}{dt}(\varepsilon \gamma_\varepsilon) = \theta \gamma_\varepsilon$, we infer that
\[
    \varepsilon(t) \gamma_\varepsilon(t) \leq \varepsilon(s) \gamma_\varepsilon(s) + \int_s^t |\theta(\tau)| \gamma_\varepsilon(\tau) \, d\tau.
\]
Consequently, there exists a constant $C'> 0$ depending on $C$, $D$, $\gamma_\varepsilon(s)$ and $\varepsilon(s)$ such that
\[  
    \varepsilon^2(t) \gamma_\varepsilon^2(t)m(t) + 2 D^2\int_s^t |\theta(\tau)| \varepsilon(\tau) \gamma_\varepsilon^2(\tau) \, d\tau \leq C' \left( 1 + \int_s^t |\theta(\tau)| \gamma_\varepsilon(\tau) \, d\tau \right)^2.
\]
Then, setting $S(t) := \gamma_\varepsilon(t) V(t)$, the estimate~\eqref{2Rate.03} becomes
\begin{equation}\label{2Rate.04}
     S(t)^2  \leq C'  \left( 1 + \int_s^t |\theta(\tau)| \gamma_\varepsilon(\tau) \, d\tau \right)^2 + 2D \int_s^t |\theta(\tau)| \gamma_\varepsilon(\tau) S(\tau) \, d\tau, \quad \text{ for a.e. } t \in E, t \geq s.
\end{equation}
Thanks to Step 3, $S = \gamma_\varepsilon V$ is locally bounded.
Since $|\theta|\gamma_\varepsilon \in L^1_{\mathrm{loc}}((t_0, + \infty))$, the mapping
\[
    t \longmapsto \int_s^t |\theta(\tau)| \gamma_\varepsilon(\tau) S(\tau) \, d\tau
\]
is continuous.
Therefore, since $E$ has full measure and $t \longmapsto S(t)$ is lower semicontinuous, the estimate~\eqref{2Rate.04} holds for every $t \geq s$.
Set
\[
    R(t)^2 :=   C'  \left( 1 + \int_s^t |\theta(\tau)| \gamma_\varepsilon(\tau) \, d\tau \right)^2 + 2D \int_s^t |\theta(\tau)| \gamma_\varepsilon(\tau) S(\tau) \, d\tau.
\]
Note that $R$ is locally absolutely continuous.
Hence, a direct computation yields, for a.e. $t \geq s$
\[
    \dot{R}(t) \leq \big( \sqrt{C'} + D \big) |\theta(t)|\gamma_\varepsilon(t).
\]
Therefore, integrating with respect to $t$, there exists a constant $C'' > 0$ such that
\[
    S(t) \leq R(t) \leq C'' \left( 1 + \int_{t_0}^t |\theta(\tau)| \gamma_\varepsilon(\tau) \, d\tau \right).
\]
Consequently, 
\[
    \| A_{\mu_t} + \varepsilon(t) \id  \|_{L^2(\mu_t)^d} = \dfrac{S(t)}{\gamma_\varepsilon(t)} \leq \dfrac{C''}{\gamma_\varepsilon(t)} \left( 1 + \int_{t_0}^t |\dot{\varepsilon}(\tau) + \varepsilon^2(\tau)| \gamma_\varepsilon(\tau) \, d\tau \right).
\]
An argument analogous to that in Step 3 gives, for a.e. $t \geq s$,
\[
    |\mu'|(t) \leq  \dfrac{C''}{\gamma_\varepsilon(t)} \left( 1 + \int_{t_0}^t |\dot{\varepsilon}(\tau) + \varepsilon^2(\tau)| \gamma_\varepsilon(\tau) \, d\tau \right)
\]
and
\[
    \| A_{\mu_t} \|_{L^2(\mu_t)^d} \leq  \dfrac{C''}{\gamma_\varepsilon(t)} \left( 1 + \int_{t_0}^t |\dot{\varepsilon}(\tau) + \varepsilon^2(\tau)| \gamma_\varepsilon(\tau) \, d\tau \right) + D \varepsilon(t).
\]

\medskip\noindent
\textbf{Step 5:} \emph{Vanishing of the terms and convergence of the trajectory.}
Recall that 
\begin{align*}
    R_1(t) = & \, \dfrac{1}{\gamma_\varepsilon(t) } \left( 1 + \int_{t_0}^t |\dot{\varepsilon}(\tau)| \gamma_\varepsilon(\tau) \, d\tau \right) \\
    R_2(t) = & \,\dfrac{1}{\gamma_\varepsilon(t) } \left( 1 + \int_{t_0}^t |\dot{\varepsilon}(\tau) + \varepsilon^2(\tau)| \gamma_\varepsilon(\tau) \, d\tau \right).
\end{align*}
Combining Step 3 and Step 4, we infer that
\[
    |\mu'|(t) = O( \min \{ R_1(t), R_2(t) \} ) \quad \text{ and } \quad \| A_{\mu_t} \|_{L^2(\mu_t)^d} = O \big(  \min \{ R_1(t), R_2(t) \} + \varepsilon(t)  \big).
\]
Since $\mathlarger{\int}_{t_0}^{+\infty}\varepsilon(t)\,dt=+\infty$, we have
$\gamma_\varepsilon(t)\to+\infty$ and hence $1/\gamma_\varepsilon(t)\to0$. We claim
that also
\begin{equation}\label{eq.thm2_weighted_term}
\frac1{\gamma_\varepsilon(t)}
\int_{t_0}^t|\dot\varepsilon(\tau)|\gamma_\varepsilon(\tau)\,d\tau
\longrightarrow 0
\qquad\text{as }t\to+\infty.
\end{equation}
Fix $T>t_0$. Since $\varepsilon>0$, the function $\gamma_\varepsilon$ is
nondecreasing, so that $\gamma_\varepsilon(\tau)\leq\gamma_\varepsilon(t)$ for
$\tau\leq t$ and therefore
\[
\frac1{\gamma_\varepsilon(t)}
\int_{t_0}^t|\dot\varepsilon(\tau)|\gamma_\varepsilon(\tau)\,d\tau
\leq
\frac1{\gamma_\varepsilon(t)}
\int_{t_0}^T|\dot\varepsilon(\tau)|\gamma_\varepsilon(\tau)\,d\tau
+
\int_T^{+\infty}|\dot\varepsilon(\tau)|\,d\tau
\quad \text{ for every } t\geq T.
\]
Letting first $t\to+\infty$ and then $T\to+\infty$ and using
$\dot\varepsilon\in L^1(t_0,+\infty)$, we obtain
\eqref{eq.thm2_weighted_term}. Since $\varepsilon(t)\to0$ as well, the estimates of
Step~3 give $\|A_{\mu_t}\|_{L^2(\mu_t)^d}\to0$ as $t\to+\infty$ and
$|\mu'|(t)\to0$ as $t\to+\infty$ through the set of points at which the metric
derivative is defined.

\medskip\noindent
It remains to prove the convergence of the trajectory.
Set $p(t):=\mathrm{Var}(\mu_t)$.
Since $m_{\mu_t} = \mathlarger{\int} e_t\,d\eta$, the mapping $t\longmapsto m_{\mu_t}$ is locally absolutely continuous and \eqref{eq.thm2_characteristic} gives
\begin{align*}
\dot p(t)
&=
-2\int\left\langle A_{\mu_t}(e_t),e_t-m_{\mu_t}\right\rangle\,d\eta-2\varepsilon(t)p(t)\\
&=
-\int_{\mathbb R^{2d}}\left\langle A_{\mu_t}(x)-A_{\mu_t}(y),x-y\right\rangle\,d(\mu_t\otimes\mu_t)(x,y)-2\varepsilon(t)p(t)
\leq -2\varepsilon(t)p(t)
\end{align*}
for almost every $t\in[t_0,+\infty)$, where the last inequality follows from the first assertion of Proposition~\ref{p1}.
Multiplication by $\gamma_\varepsilon^2$ and integration yield
\[
\mathrm{Var}(\mu_t)\leq\frac{\mathrm{Var}(\mu_0)}{\gamma_\varepsilon(t)^2}\longrightarrow 0
\quad\text{as }t\longrightarrow+\infty.
\]

\medskip\noindent
Let $t_n\to+\infty$.
By \eqref{eq.thm2_moment_bound}, the sequence $(m_{\mu_{t_n}})_n$ is bounded, so that, after passing to a subsequence, $m_{\mu_{t_n}}\to\bar x$ for some $\bar x\in\mathbb R^d$.
By \eqref{dirac_distance},
\[
\mathcal W_2^2(\mu_{t_n},\delta_{\bar x})
=\mathrm{Var}(\mu_{t_n})+|m_{\mu_{t_n}}-\bar x|^2\longrightarrow 0.
\]
Demicontinuity gives $(\id,A_{\mu_{t_n}})_\#\mu_{t_n}\rightharpoonup(\id,A_{\delta_{\bar x}})_\#\delta_{\bar x}$ narrowly, while Step 5 gives $\|A_{\mu_{t_n}}\|_{L^2(\mu_{t_n})^d}\to0$.
Thus $A_{\delta_{\bar x}}(\bar x)= 0$ and $\bar x\in K$.
Consequently, if $z:=\mathrm{proj}_K(0)$, the characterization of the projection onto $K$ yields
\[
\limsup_{t\to+\infty}\left\langle m_{\mu_t}-z,-z\right\rangle\leq 0.
\]

\medskip\noindent
Set $q(t):=\mathcal W_2^2(\mu_t,\delta_z)=\|e_t-z\|_{L^2(\eta)^d}^2$.
Since $A_{\delta_z}(z)= 0$, \eqref{eq.thm2_characteristic} and the first assertion of Proposition~\ref{p1} give
\[
\frac12\dot q(t)
\leq-\varepsilon(t)\int\left\langle e_t,e_t-z\right\rangle\,d\eta
=-\varepsilon(t)q(t)+\varepsilon(t)\left\langle m_{\mu_t}-z,-z\right\rangle
\]
for almost every $t\in[t_0,+\infty)$.
Given $\delta>0$, choose $T>t_0$ such that $\langle m_{\mu_t}-z,-z\rangle\leq\delta$ for every $t\geq T$.
Multiplying the preceding inequality by $2\gamma_\varepsilon^2$ and integrating over $[T,t]$, we obtain
\[
q(t)\leq\delta+\bigl(q(T)-\delta\bigr)\frac{\gamma_\varepsilon(T)^2}{\gamma_\varepsilon(t)^2}
\quad\text{for every }t\geq T.
\]
Letting first $t\to+\infty$ and then $\delta\searrow0$ gives $\mathcal W_2(\mu_t,\delta_z)\to0$.
This completes the proof.

\subsection{Proof of Corollary~\ref{mccr}}
\noindent\emph{Case $t_0>0$ and $\varepsilon(t)=\alpha/t$.}
It is not difficult to check that the function $\varepsilon(t) = \alpha/t$ satisfies the assumptions in Theorem~\ref{thm2}.
Also, $ \gamma_\varepsilon(t) = (t/t_0)^\alpha$ for every $t > t_0$ and
\[
    \dot{\varepsilon}(t)  + \varepsilon^2(t) = \dfrac{\alpha(\alpha  - 1)}{t^2}.
\]
Consequently, we get
\[
    R_2(t) := \dfrac{1}{\gamma_\varepsilon(t)} \left( 1 + \int_{t_0}^t |\dot{\varepsilon}(\tau) + \varepsilon^2(\tau)| \gamma_\varepsilon(\tau) \, d\tau \right) = \dfrac{t_0^\alpha}{t^\alpha} \left( 1 + \dfrac{\alpha |\alpha - 1|}{t_0^\alpha} \int_{t_0}^t \tau^{\alpha - 2} \,d\tau \right).
\]
\noindent
In the case $\alpha \in (0, 1)$, one has $R_2(t) = O(t^{- \alpha})$ since
\[
    \int_{t_0}^{+\infty} \tau^{\alpha - 2} \,d\tau < + \infty.
\] 
In the case $\alpha = 1$, we have $\dot{\varepsilon}(t) + \varepsilon^2(t) = 0$ for all $t \geq t_0$ and hence $R_2(t) = O(t^{-1})$. 
In the case $\alpha > 1$, one gets $R_2(t) = O(t^{-1})$ since 
\[
    \int_{t_0}^t \tau^{\alpha - 2} \, d\tau = O(t^{\alpha - 1}) \quad \text{ as } t \to + \infty.
\]
Consequently, for any $\alpha > 0$, one has $R_2(t) = O(t^{- \min \{ \alpha, 1\} })$ as $t \to + \infty$.
Theorem~\ref{thm2} then yields
\[
    |\mu'|(t) = O\big( t^{- \min \{1, \alpha\} }  \big) \quad \text{ and } \quad \| A_{\mu_t} \|_{L^2(\mu_t)^d} = O\big( t^{- \min \{1, \alpha \} }  \big).
\]
This proves the assertion in this case.

\medskip\noindent
\emph{Case $t_0>0$ and $\varepsilon(t)=\alpha/t^q$.}
The function $\varepsilon$ is smooth and positive on $[t_0,+\infty)$, it tends to zero and it is nonincreasing, so that
$\mathlarger{\int}_{t_0}^{+\infty}|\dot\varepsilon(t)|\,dt=\varepsilon(t_0)<+\infty$; moreover,  $\mathlarger{\int}_{t_0}^{+\infty}\varepsilon(t)\,dt=+\infty$ because $q<1$. 
Hence Theorem~\ref{thm2} applies and
\[
\gamma_\varepsilon(t)
=
\exp\left(\frac{\alpha}{1-q}\left(t^{1-q}-t_0^{1-q}\right)\right)
\quad \text{ for every } t\geq t_0.
\]
Since $1-q>0$, the function $\gamma_\varepsilon$ grows faster than any power of
$t$ and in particular $1/\gamma_\varepsilon(t)=O(1/t)$.
\smallbreak\noindent
We claim that
\[
\Xi(t):=
\frac1{\gamma_\varepsilon(t)}
\int_{t_0}^t|\dot\varepsilon(\tau)|\gamma_\varepsilon(\tau)\,d\tau
=
O\left(\frac1t\right).
\]
Indeed, $|\dot\varepsilon(\tau)|=\alpha q\tau^{-q-1}=q\varepsilon(\tau)/\tau$ and
$\dot\gamma_\varepsilon=\varepsilon\gamma_\varepsilon$, so that
$|\dot\varepsilon(\tau)|\gamma_\varepsilon(\tau)=q\dot\gamma_\varepsilon(\tau)/\tau$.
Let $t\geq2t_0$. Splitting the integral at $t/2$ and using that
$\tau\longmapsto q/\tau$ is nonincreasing and that $\gamma_\varepsilon$ is
nondecreasing, we obtain
\[
\int_{t_0}^t\frac{q}{\tau}\,\dot\gamma_\varepsilon(\tau)\,d\tau
\leq
\frac{q}{t_0}\int_{t_0}^{t/2}\dot\gamma_\varepsilon(\tau)\,d\tau
+
\frac{2q}{t}\int_{t/2}^{t}\dot\gamma_\varepsilon(\tau)\,d\tau
\leq
\frac{q}{t_0}\gamma_\varepsilon\!\left(\frac t2\right)
+
\frac{2q}{t}\gamma_\varepsilon(t),
\]
whence
\[
\Xi(t)
\leq
\frac{q}{t_0}\,
\frac{\gamma_\varepsilon(t/2)}{\gamma_\varepsilon(t)}
+
\frac{2q}{t}
=
\frac{q}{t_0}
\exp\left(-\frac{\alpha\bigl(1-2^{q-1}\bigr)}{1-q}\,t^{1-q}\right)
+
\frac{2q}{t}.
\]
As $q<1$, the exponential above tends to zero faster than any power of $1/t$,
which proves the claim.

\medskip\noindent
The two estimates of Theorem~\ref{thm2} now give
$|\mu'|(t)=O(1/\gamma_\varepsilon(t))+O(\Xi(t))=O(1/t)$ and
$\|A_{\mu_t}\|_{L^2(\mu_t)^d}=O(1/t)+O(\varepsilon(t))=O(t^{-q})$, the last
identity because $t^{-1}\leq t^{-q}$ for $t\geq1$. This proves the assertion in this case.

\medskip\noindent
\emph{Case $t_0>1$ and $\varepsilon(t)=1/(t\log t)$.}
The function $\varepsilon$ is smooth and positive on $[t_0,+\infty)$ and it tends
to zero; since
$\dot\varepsilon(t)=-(\log t+1)/(t^2\log^2t)<0$ for $t>1$, it is nonincreasing and
$\mathlarger{\int}_{t_0}^{+\infty}|\dot\varepsilon(t)|\,dt = \varepsilon(t_0)<+\infty$. Moreover
$\mathlarger{\int}_{t_0}^t\varepsilon(s)\,ds=\log\log t-\log\log t_0\to+\infty$. Hence
Theorem~\ref{thm2} applies and
\[
\gamma_\varepsilon(t)=\frac{\log t}{\log t_0}
\quad \text{ for every } t\geq t_0,
\]
so that $1/\gamma_\varepsilon(t)=O(1/\log t)$.
\smallbreak\noindent
On the other hand,
\[
|\dot\varepsilon(\tau)|\gamma_\varepsilon(\tau)
=
\frac{\log\tau+1}{\tau^2\log^2\tau}\,\frac{\log\tau}{\log t_0}
=
\frac1{\log t_0}
\left(\frac1{\tau^2}+\frac1{\tau^2\log\tau}\right)
\leq
\frac1{\log t_0}
\left(1+\frac1{\log t_0}\right)\frac1{\tau^2}
\quad \text{ for every } \tau\geq t_0,
\]
because $\log\tau\geq\log t_0>0$. The right-hand side being integrable on
$[t_0,+\infty)$, the integral
$\mathlarger{\int}_{t_0}^t|\dot\varepsilon(\tau)|\gamma_\varepsilon(\tau)\,d\tau$ is bounded
uniformly in $t$ and therefore
\[
\frac1{\gamma_\varepsilon(t)}
\int_{t_0}^t|\dot\varepsilon(\tau)|\gamma_\varepsilon(\tau)\,d\tau
=
O\left(\frac1{\log t}\right).
\]
Since $\varepsilon(t)=1/(t\log t)=O(1/\log t)$ as well, the two estimates of
Theorem~\ref{thm2} give the assertion in this case.

\subsection{Proof of Proposition~\ref{prop:wass_nonexpansive}.}
Let $\mu,\nu\in\mathcal P_2(\mathbb R^d)$ and $\gamma\in\Gamma_o(\mu,\nu)$.
Since $A_\mu=\id-B_\mu$, for every $(x,y)\in\mathbb R^{2d}$, a direct computation gives
\[
2\left\langle A_\mu(x)-A_\nu(y),x-y\right\rangle-\left|A_\mu(x)-A_\nu(y)\right|^2 = |x-y|^2-\left|B_\mu(x)-B_\nu(y)\right|^2.
\]
Integrating with respect to $\gamma$, we obtain
\[
\int_{\mathbb R^{2d}}\left|B_\mu(x)-B_\nu(y)\right|^2\,d\gamma(x,y)
\leq
\int_{\mathbb R^{2d}}|x-y|^2\,d\gamma(x,y)
\]
if and only if
\[
\int_{\mathbb R^{2d}}\left\langle A_\mu(x)-A_\nu(y),x-y\right\rangle\,d\gamma(x,y)
\geq
\frac12\int_{\mathbb R^{2d}}\left|A_\mu(x)-A_\nu(y)\right|^2\,d\gamma(x,y).
\]
Thus $B$ is nonexpansive if and only if $A$ is ($1/2$)-cocoercive.
Moreover, $A_\mu = 0$ $\mu$--a.e. if and only if $B_\mu = \id$ $\mu$--a.e. and hence $\mathrm{fix}(B)  = \mathrm{zer}(A) $. 

\medskip

\noindent
Assume that $B$ is nonexpansive.
Let $\mu,\nu\in\mathcal P_2(\mathbb R^d)$ and let $\gamma\in\Gamma_o(\mu,\nu)$.
Then
\[
(B_\mu\circ\mathrm{proj}_1,B_\nu\circ\mathrm{proj}_2)_\#\gamma
\in
\Gamma\bigl((B_\mu)_\#\mu,(B_\nu)_\#\nu\bigr).
\]
Consequently,
\begin{align*}
\mathcal W_2\bigl(\mathcal B(\mu),\mathcal B(\nu)\bigr)^2 \leq \int_{\mathbb R^{2d}} \left|B_\mu(x)-B_\nu(y)\right|^2\,d\gamma(x,y) \leq
\int_{\mathbb R^{2d}}|x-y|^2\,d\gamma(x,y)
= \mathcal W_2(\mu,\nu)^2.
\end{align*}
Hence $\mathcal B$ is nonexpansive.
This completes the proof. 
\hfill$\square$

%%%%%%%%%%%%%%%%%%%%% Wass Halpern

\subsection{Proof of Lemma~\ref{lemfixmoment} and Theorem~\ref{thm:wass_halpern}}

Let us start by recording some properties related to nonexpansive maps.

\begin{lemma}\label{lem:FixT}
Let $B$ be a nonexpansive map on $\mathcal{P}_2(\R^d)$ such that $\mathrm{fix}(B) \neq \varnothing$.
Set $K_B := \{ x \in \R^d \, : \, \delta_x \in \mathrm{fix}(B) \}$.
Then, the following assertions hold:
\begin{itemize}
    \item[(i)] If $\mu \in \mathrm{fix}(B)$, then $m_\mu \in K_B$, where $m_\mu$ is the barycenter of $\mu$.
    
    \item[(ii)] $K_B$ is a nonempty closed convex subset of $\R^d$.

     \item[(iii)] For every $v \in \R^d$, it holds
     \[
        \{ \delta_{\mathrm{proj}_{K_B}(v)} \} = \mathrm{argmin}_{\mu \in \mathrm{fix}(B)} \mathcal{W}_2(\mu, \delta_v)^2.
     \]
\end{itemize}
\end{lemma}

\begin{proof}
Set $A := \id - B$.
Thanks to Proposition~\ref{prop:wass_nonexpansive}, $A$ is $(1/2)$--cocoercive, and moreover,  $\mathrm{fix}(B) = \mathrm{zer}(A)$.
As shown in the proof of Theorem~\ref{thm1}, $A$ is demicontinuous and so, by Proposition~\ref{p1}, $\mathrm{fix}(B)$ is convex along the interpolations induced by arbitrary couplings.

\medskip

\noindent\emph{(i)}
Fix $\mu \in \mathrm{fix}(B)$. 
Note that $\delta_{m_\mu} \otimes \mu$ is the unique coupling between $\delta_{m_\mu}$ and $\mu$. 
Hence, by the nonexpansivity of $B$ and $\mu \in \mathrm{fix}(B)$, one has
\[
\int_{\mathbb R^d}|B_{\delta_{m_\mu}}(m_\mu)-y|^2\,d\mu(y)\leq\int_{\mathbb R^d}|m_\mu-y|^2\,d\mu(y) = \mathrm{Var}(\mu).
\]
On the other hand, by the definition of $m_\mu$, we get
\[
\int_{\mathbb R^d}|B_{\delta_{m_\mu}}(m_\mu)-y|^2\,d\mu(y)=|B_{\delta_{m_\mu}}(m_\mu)-m_\mu|^2+\mathrm{Var}(\mu).
\]
Thus, $B_{\delta_{m_\mu}}(m_\mu)=m_\mu$ and so $m_\mu\in K_B$.

\medskip

\noindent\emph{(ii)}
It follows from \emph{(i)} that $K_B \neq \varnothing$, since $\mathrm{fix}(B) \neq \varnothing$. 
Let $x,y\in K_B$ and $t\in[0,1]$. 
Since $\delta_x,\delta_y\in\mathrm{fix}(B)$ and $\mathrm{fix}(B)$ is convex along the interpolations induced by arbitrary couplings, we infer that $\delta_{(1-t)x+ty}\in\mathrm{fix}(B)$.
Hence, $(1-t)x+ty\in K_B$ and so $K_B$ is convex.

\medskip

\noindent 
Let now $(x_k)_{k \in \mathbb N}\subseteq K_B$ with $x_k\to x$ as $k \to +\infty$.
Applying \eqref{eq:wass_nonexpansive} to the case $\mu = \delta_x$ and $\nu = \delta_{x_k}$ and $\gamma = \delta_x \otimes \delta_{x_k}$ and using $B_{\delta_{x_k}}(x_k) = x_k$, we get
\[
|B_{\delta_x}(x)-x_k|^2\leq|x-x_k|^2 \quad \text{for every } k \geq 1.
\]
Letting $k\to+\infty$ yields $B_{\delta_x}(x)=x$. 
Thus, $x\in K_B$ and $K_B$ is closed.

\medskip

\noindent\emph{(iii)}
Fix $v\in\mathbb R^d$.
For every $\mu\in\mathrm{fix}(B)$, (i) gives $m_\mu\in K_B$ and therefore
\[
\mathcal W_2^2(\mu,\delta_v) = \int_{\mathbb R^d}|x-v|^2\,d\mu(x) = \mathrm{Var}(\mu) + |m_\mu-v|^2\geq |\mathrm{proj}_{K_B}(v) - v|^2 = \mathcal W_2^2(\delta_{\mathrm{proj}_{K_B}(v)},\delta_v).
\]
This shows that
\[
    \delta_{\mathrm{proj}_{K_B}(v)} \in \mathrm{argmin}_{\mu \in \mathrm{fix}(B)} \mathcal{W}_2^2(\mu, \delta_v).
\]
To prove the reverse, assume that there is $\mu_\ast \in \mathrm{fix}(B)$ such that
\[
    \mathcal W_2^2(\mu_\ast,\delta_v) \leq \mathcal W_2^2(\mu,\delta_v) \quad \text{ for all } \mu \in \mathrm{fix}(B).
\]
It follows that $\mathrm{Var}(\mu_\ast) = 0$ and $|m_{\mu_\ast} - v | = |\mathrm{proj}_{K_B}(v) - v|$.
Consequently, $\mu_\ast = \delta_{ m_{\mu_\ast}  }$ and, by the uniqueness of the projection onto $K_B$, $m_{\mu_\ast} = \mathrm{proj}_{K_B}(v)$.
Therefore, $\mu_\ast = \delta_{\mathrm{proj}_{K_B}(v)}$, completing the proof.
\end{proof}

Since $Q = K_B$, assertions \emph{(ii)} and \emph{(iii)} give Lemma~\ref{lemfixmoment}.

\noindent We are now ready to prove Theorem~\ref{thm:wass_halpern}.
%%%%%%%%% Proof of Halpern typ algorithm %%%%%%%%%%%%%%%%%%%

\begin{proof}[Proof of Theorem~\ref{thm:wass_halpern}]
Set $K_B := \{ x \in \R^d : \delta_x \in \mathrm{fix}(B) \}$ and $z := \mathrm{proj}_{K_B}(v)$.
Thanks to Lemma~\ref{lem:FixT}, $K_B$ is nonempty, closed and convex and $\delta_z$ is the unique minimizer of the following problem: 
\[
    \min_{\mu \in \mathrm{fix}(B)} \, \mathcal{W}_2^2(\mu, \delta_v).
\]

\medskip

\noindent\textbf{Step 1:} \emph{Nonexpansivity and boundedness.}
Define $\mathcal{B} : \mathcal{P}_2(\R^d) \to \mathcal{P}_2(\R^d)$ by $\mathcal{B}(\mu) := (B_\mu)_\# \mu$.
As we noted in Remark~\ref{rmk.nonexp}, $\mathcal B$ is nonexpansive on $(\mathcal P_2(\mathbb R^d),\mathcal W_2)$.
Moreover, since $z\in K_B$, we have $\mathcal B(\delta_z)=\delta_z$.
From \eqref{eq:wass_halpern} and the Minkowski inequality in $L^2(\mathcal{B}(\mu_k))^d$, we have
\[
\mathcal W_2(\mu_{k+1},\delta_z) \leq \alpha_k|v-z|+(1-\alpha_k)\mathcal W_2(\mathcal B(\mu_k),\delta_z) \leq \alpha_k|v-z|+(1-\alpha_k)\mathcal W_2(\mu_k,\delta_z) \quad \text{for every } k \geq 0.
\]
Hence, for $C_0 := \max\{\mathcal W_2(\mu_0,\delta_z),|v-z|\}$,
we obtain
\[
\mathcal W_2(\mu_k,\delta_z)\leq C_0 \quad \text{for every } k \geq 0.
\]
Consequently, there exists a constant $C > 0$ such that
\begin{equation}\label{MuDel_bound}
\mathcal W_2(\mathcal B(\mu_k),\delta_v)\leq C \quad \text{for every } k \geq 0.
\end{equation}

\medskip

\noindent\textbf{Step 2:} \emph{Vanishing of the variance.}
For every $\mu\in\mathcal P_2(\mathbb R^d)$, applying \eqref{eq:wass_nonexpansive} with $\mu$, $\delta_{m_\mu}$ and $\mu\otimes\delta_{m_\mu}$ gives
\begin{equation}\label{V_Est01}
\mathrm{Var}(\mathcal B(\mu))\leq\int_{\mathbb R^d}|B_\mu(x)-B_{\delta_{m_\mu}}(m_\mu)|^2\,d\mu(x)\leq\mathrm{Var}(\mu).
\end{equation}
If $S_k(x) := \alpha_k  v + (1 - \alpha_k)x$ for $x \in \R^d$, then for $\mu \in \mathcal{P}_2(\R^d)$, $m_{(S_k)_\# \mu} = \alpha_k v + (1 - \alpha_k) m_{\mu}$ and hence
\begin{equation}\label{V_Est02}
    \mathrm{Var}((S_k)_\# \mu) = \int_{\R^d} |S_k(x) - m_{(S_k)_\#\mu}|^2 \, d\mu(x) = (1 - \alpha_k)^2 \mathrm{Var}(\mu) \quad \text{for every } k \geq 0.
\end{equation}
Combining~\eqref{V_Est01} and~\eqref{V_Est02}, we obtain
\[
\mathrm{Var}(\mu_{k+1})=(1-\alpha_k)^2\mathrm{Var}(\mathcal B(\mu_k))\leq(1-\alpha_k)^2\mathrm{Var}(\mu_k) \quad \text{for every } k \geq 0.
\]
Consequently,
\[
\mathrm{Var}(\mu_k)^{1/2}\leq \mathrm{Var}(\mu_0)^{1/2}\prod_{j=0}^{k-1}(1-\alpha_j) \longrightarrow 0 \quad \text{ as } k \longrightarrow +\infty,
\]
where the convergence follows from \eqref{eq:halpern_parameters}.
Further, $|m_{\mu_k}-z|\leq\mathcal W_2(\mu_k,\delta_z)\leq C_0$, hence $(m_{\mu_k})_{k\in\mathbb N}$ is bounded.

\medskip

\noindent\textbf{Step 3:} \emph{Asymptotic regularity.}
Set $d_k:=\mathcal W_2(\mu_k,\mu_{k-1})$ for every $k \geq 1$.
We claim that 
\begin{equation}\label{St3-01}
    d_{k + 1} \leq (1 - \alpha_k) d_k + C | \alpha_k - \alpha_{k - 1}| \quad \text{ for every } k \geq 1.
\end{equation}
Indeed, let $\gamma \in \Gamma_o \big( \mathcal{B}(\mu_k), \mathcal{B}(\mu_{k - 1}) \big)$.
Since
\[
    \mu_{k + 1} = (\alpha_k v + (1 - \alpha_k)\id)_\# \mathcal{B}(\mu_k) \quad \text{ and } \quad \mu_k = (\alpha_{k - 1} v + (1 - \alpha_{k - 1})\id)_\# \mathcal{B}(\mu_{k - 1}),
\]
we infer that
\[
    \pi := \big( \alpha_k v + (1 - \alpha_k)\mathrm{proj}_1, \alpha_{k - 1} v + (1 - \alpha_{k - 1})\mathrm{proj}_2 \big)_\# \gamma \in \Gamma(\mu_{k + 1}, \mu_k) \quad \text{for every } k \geq 1.
\]
Therefore, using Minkowski's inequality, we get
\begin{align*}
    d_{k + 1} \leq & \, \left( \int_{\R^{2d}} |x - y|^2 \, d\pi(x, y) \right)^{1/2} \\
    = & \, \left( \int_{\R^{2d}} |\alpha_k v + (1 - \alpha_k)x - \alpha_{k - 1} v - (1 - \alpha_{k - 1})y|^2 \, d\gamma(x, y) \right)^{1/2} 
    \\
    \leq & \, (1 - \alpha_k) \left( \int_{\R^{2d}} |x - y|^2 \, d\gamma(x, y) \right)^{1/2} + |\alpha_k - \alpha_{k - 1}| \left( \int_{\R^{2d}} |v - y|^2 \, d\gamma(x, y) \right)^{1/2} \\
    = & \, (1 - \alpha_k) \mathcal{W}_2(\mathcal{B}(\mu_k), \mathcal{B}(\mu_{k - 1})) + |\alpha_k - \alpha_{k - 1}| \mathcal{W}_2(\delta_v, \mathcal{B}(\mu_{k - 1}) ) \\
    \leq & \, (1 - \alpha_k) d_k + C |\alpha_k - \alpha_{k - 1}| \quad \text{for every } k \geq 1,
\end{align*}
where in the last estimate, we have used the estimate~\eqref{MuDel_bound} and the fact that $\mathcal{B}$ is nonexpansive in $(\mathcal{P}_2(\R^d), \mathcal{W}_2)$.
Next, by an argument analogous to that used in the proof of \cite[Theorem 4.3]{BN_2026}, we deduce
\begin{equation}\label{dk_lim}
d_k \longrightarrow 0 \quad \text{ as } \quad k \longrightarrow +\infty.
\end{equation}

\medskip

\noindent
From~\eqref{eq:wass_halpern}, notice that $\eta := (\alpha_k v + (1 - \alpha_k)\id, \id)_\# \mathcal{B}(\mu_k) \in \Gamma(\mu_{k + 1}, \mathcal{B}(\mu_k))$.
Therefore, 
\[
    \mathcal{W}_2^2(\mu_{k + 1}, \mathcal{B}(\mu_k)) \leq \int_{\R^{2d}} |x - y|^2 \, d\eta(x, y) = \alpha_k^2 \mathcal{W}_2^2(\mathcal{B}(\mu_k), \delta_v) \quad \text{for every } k \geq 0.
\]
Thanks to~\eqref{MuDel_bound}, we get
\begin{equation}\label{St3-02}
    \mathcal{W}_2(\mu_{k + 1}, \mathcal{B}(\mu_k)) \leq C \alpha_k \xrightarrow{k \to + \infty} 0.
\end{equation}
Combining~\eqref{dk_lim} and~\eqref{St3-02}, we deduce 
\[
\mathcal W_2(\mu_k,\mathcal B(\mu_k))\leq d_{k+1}+C\alpha_k \xrightarrow{k \to + \infty} 0.
\]

\medskip

\noindent\textbf{Step 4:} \emph{Convergence to $\delta_z$.}
Since $(m_{\mu_{k} })_{k \in \mathbb N}$ is bounded, there exist a subsequence $( m_{\mu_{k_n}} )_{n \in \mathbb N}$ and $\bar z \in \R^d$ such that
\[
    m_{\mu_{k_n}} \xrightarrow{n \longrightarrow +\infty} \bar z \quad \text{ and } \quad \limsup_{k \to \infty} \, \langle m_{\mu_k} - z, v - z \rangle = \lim_{n \to \infty}\, \langle m_{\mu_{k_n}} - z, v - z \rangle = \langle \bar z - z, v - z \rangle.
\]
It follows from Step~2 that
\[
\mathcal W_2^2(\mu_{k_n},\delta_{\bar z})=\mathrm{Var}(\mu_{k_n})+|m_{\mu_{k_n}}-\bar z|^2 \longrightarrow 0 \quad \text{ as } n \longrightarrow +\infty.
\]
By the nonexpansivity of $\mathcal B$ and Step~3, we get
\[
\mathcal W_2(\mathcal B(\delta_{\bar z}),\delta_{\bar z})\leq 2 \mathcal W_2(\mu_{k_n},\delta_{\bar z})+\mathcal W_2(\mu_{k_n},\mathcal B(\mu_{k_n})) \longrightarrow 0 \quad \text{ as } \quad n \longrightarrow +\infty.
\]
Thus, $\mathcal B(\delta_{\bar z})=\delta_{\bar z}$ and consequently $\bar z\in K_B$.
Since $z=\mathrm{proj}_{K_B}(v)$, it follows that $\langle \bar z-z,v-z\rangle\leq0$ and so
\[
\limsup_{k\to+\infty} \, \langle m_{\mu_k}-z,v-z\rangle\leq0.
\]
Moreover, using \eqref{eq:wass_halpern}, the identity $\mathcal B(\delta_z)=\delta_z$ and the nonexpansivity of $\mathcal B$, we obtain
\[
\mathcal W_2^2(\mu_{k+1},\delta_z)\leq(1-\alpha_k)\mathcal W_2^2(\mu_k,\delta_z)+2\alpha_k\langle m_{\mu_{k+1}}-z,v-z\rangle \quad \text{for every } k \geq 0.
\]
Thus, \eqref{eq:halpern_parameters} and \cite[Lemma A.9]{BN_2026} yield $\mathcal W_2(\mu_k,\delta_z)\longrightarrow 0$ as $k \longrightarrow + \infty$.
Theorem~\ref{thm:wass_halpern} is proven.
\end{proof}

\subsection{Proof of Corollary~\ref{cor:halpern_cocoercive}}
Set $B:=\id-\tau A$.
For every $\mu,\nu\in\mathcal P_2(\R^d)$ and $\gamma\in\Gamma_o(\mu,\nu)$, the $\beta$--cocoercivity of $A$ yields
\[
\begin{aligned}
\int_{\R^{2d}}|B_\mu(x)-B_\nu(y)|^2\,d\gamma(x,y)
&=\int_{\R^{2d}}|x-y-\tau(A_\mu(x)-A_\nu(y))|^2\,d\gamma(x,y)\\
&\leq \int_{\R^{2d}}|x-y|^2\,d\gamma(x,y)-\tau(2\beta-\tau)\int_{\R^{2d}}|A_\mu(x)-A_\nu(y)|^2\,d\gamma(x,y)\\
&\leq \int_{\R^{2d}}|x-y|^2\,d\gamma(x,y),
\end{aligned}
\]
where we have used $\tau\in(0,2\beta]$.
Thus, $B$ is nonexpansive.
Moreover, $\mathrm{fix}(B)=\mathrm{zer}(A)$.
Since \eqref{eq:halpern_cocoercive} coincides with the iteration \eqref{eq:wass_halpern} applied to $B$, the conclusion follows from Theorem~\ref{thm:wass_halpern}.
\hfill$\square$

%%%%%%%%%%%%%%%%%%%%%%%%%%%%
%%%%%%%%%%%%%%%%%%%%%%%%%%%%

%%%%%%%%%%%%%%%%%%%%%%%%%%%%
%%%%%%%%%%%%%%%%%%%%%%%%%%%%
\section{Examples}\label{sec:example}

\subsection{Examples for cocoercive operators}\label{sec:example1}

\begin{example}\label{Exm1}\normalfont
Let us begin with a standard potential--interaction energy on $\mathcal P_2(\R^d)$ of the form
\begin{equation}\label{eq.PoIn.E}
    \mathcal E(\mu) := \int_{\R^d} V(x) \,d\mu(x) + \dfrac{1}{2} \int_{\R^d \times \R^d} W(x - y) \, d\mu(x) \, d\mu(y) \quad \text{ for } \mu \in \mathcal{P}_2(\R^d),
\end{equation}
where $V, W \in C^1(\R^d)$  are convex and $W(x) = W(-x)$ for every $x \in \R^d$.
Such energies are classical in optimal transport, see \cite{M_1997}.
Whenever the first variation belongs to $L^2(\mu)$, it is given by (see e.g. \cite[Examples 4.2.4--4.2.5]{FG_2021})
\[
    \dfrac{\delta \mathcal E [ \mu]}{\delta \mu}(x) = V(x) + (W \ast \mu)(x) = V(x) + \int_{\R^d} W(x - y) \,d\mu(y).
\]
We will show that under suitable regularity assumptions, the velocity field
\[
    A_\mu := \nabla \dfrac{\delta \mathcal{E}[ \mu]}{\delta \mu} = \nabla V + (\nabla W \ast \mu)
\]
is cocoercive.
\end{example}

\begin{proposition}\label{prop.PoIn.E}
    Let $V, W \in C^1(\R^d)$ be convex such that $W$ is even, i.e., $W(x) = W(-x)$ for every $x \in \R^d$.
    Assume that there exist $L_V, L_W > 0$ such that $\nabla V$ and $\nabla W$ are $L_V$-- and $L_W$--Lipschitz, respectively.
    Then, the velocity field of the potential--interaction energy~\eqref{eq.PoIn.E},
    \[
        A_{ \mu} := \nabla \dfrac{\delta \mathcal E[ \mu]}{\delta \mu} = \nabla V + \nabla W \ast \mu
    \]
    is $(L_V + L_W)^{-1}$--cocoercive.
\end{proposition}

\begin{proof}
Setting $B_\mu := \nabla W \ast \mu$, we first claim that
\begin{equation}\label{Bmu.coco}
    \int_{\R^{2d}} \langle B_\mu(x) - B_\nu(y), x - y \rangle \, d\gamma(x, y) \geq \dfrac{1}{L_W} \int_{\R^{2d}} |B_\mu(x) - B_\nu(y)|^2 \, d\gamma(x, y), 
\end{equation}
for every $\mu, \nu \in \mathcal{P}_2(\R^d)$ and $\gamma \in \Gamma(\mu, \nu)$.
Fix $\mu, \nu \in \mathcal P_2(\R^d)$ and $\gamma \in \Gamma(\mu, \nu)$.
Using Fubini's theorem, we write, with $z = (x, y)$ and $z' = (x', y')$,
\begin{equation}\label{Bmu.01}
    \int_{\R^{2d}} \langle B_\mu(x) - B_\nu(y), x - y \rangle \, d\gamma(x, y) = \int_{\R^{2d} \times \R^{2d}} \langle \nabla W(x - x') - \nabla W(y - y'), x - y \rangle \, d(\gamma \otimes \gamma)(z, z').
\end{equation}
Notice that $\nabla W$ is odd since $W$ is even.
Exchanging $z$ and $z'$ then gives
\begin{equation}\label{Bmu.02}
\begin{split}
     & \, \int_{\R^{2d}} \langle B_\mu(x) - B_\nu(y), x - y \rangle \, d\gamma(x, y)\\
     = & \, \dfrac{1}{2} \int_{\R^{2d} \times \R^{2d}} \langle \nabla W(x - x') - \nabla W(y - y'), (x - x') - (y - y') \rangle \, d(\gamma \otimes \gamma)(z, z') \\
     \geq & \, \dfrac{1}{2L_W} \int_{\R^{2d} \times \R^{2d}} | \nabla W(x - x') - \nabla W(y - y') |^2 \, d(\gamma \otimes \gamma)(z, z'),
\end{split}
\end{equation}
where we have applied the Baillon--Haddad theorem to the function $W$.
Define 
\[
    b(z) := \int_{\R^{2d}} \nabla W(x - x') - \nabla W(y - y') \, d\gamma(z') = B_\mu(x) - B_\nu(y).
\]
Using the fact that 
\[
    \int_{\R^{2d}} b(z) \, d\gamma(z) = 0
\]
and expanding the standard estimate
\[
    \int_{\R^{2d} \times \R^{2d}} |\nabla W(x - x') - \nabla W(y - y') - b(z) + b(z')|^2 \, d(\gamma \otimes \gamma)(z, z') \geq 0,
\]
we obtain that
\begin{equation}\label{Bmu.03}
    \int_{\R^{2d} \times \R^{2d}} |\nabla W(x - x') - \nabla W(y - y')|^2 \, d (\gamma \otimes \gamma)(z, z') \geq 2 \int_{\R^{2d}}|b(z)|^2 \, d\gamma(z).
\end{equation}
Combining~\eqref{Bmu.02} and~\eqref{Bmu.03}, we get the  desired estimate~\eqref{Bmu.coco}.

\medskip

\noindent
To proceed, applying the Baillon--Haddad theorem to the function $V$, we observe that 
\begin{equation}\label{V.coco}
    \int_{\R^{2d}} \langle \nabla V(x) - \nabla V(y), x - y \rangle \, d\gamma(x, y) \geq \dfrac{1}{L_V} \int_{\R^{2d}} |\nabla V(x) - \nabla V(y) |^2 \, d\gamma (x, y),
\end{equation}
for every $\mu, \nu \in \mathcal{P}_2(\R^d)$ and $\gamma \in \Gamma(\mu, \nu)$.
Combining~\eqref{Bmu.coco} and~\eqref{V.coco}, a direct computation yields
\[
    \int_{\R^{2d}} \langle A_\mu(x) - A_\nu(y), x - y \rangle \, d\gamma(x, y) \geq \dfrac{1}{L_V + L_W} \int_{\R^{2d}} |A_\mu(x) - A_\nu(y)|^2 \, d\gamma(x, y),
\]
for every $\mu, \nu \in \mathcal{P}_2(\R^d)$ and $\gamma \in \Gamma(\mu, \nu)$.
This completes the proof.
\end{proof}

\begin{example}\label{Exm2}\normalfont
One can also consider higher--order interaction energies, see e.g. \cite[Example 3.6]{AS_2007} and \cite[Chapter 7]{S_2015}.
Fix $k \in \mathbb{N}$ with $k \geq 2$.
Let $W \in C^1(\R^{dk})$ be a convex $L_W$--smooth function such that $W$ is invariant under permutations of its $k$ variables.
Consider the interaction energy
\[
    \mathcal W_k : \mu \longmapsto \dfrac{1}{k} \int_{\R^d} \cdots \int_{\R^d} W(x_1, \cdots, x_k) \, d\mu^{\otimes k}(x_1, \cdots, x_k),
\]
where $\mu^{\otimes k} := \mu \times \cdots \times \mu$.
An argument analogous to that in the proof of Proposition~\ref{prop.PoIn.E} implies that the velocity field
\[
    A_\mu(x) := \nabla \dfrac{\delta \mathcal W_k[\mu]}{\delta \mu}(x) = \int_{\R^{d(k - 1)}} \nabla_{x_1} W(x, x_2, \cdots, x_k) \, d\mu^{\otimes (k - 1)}(x_2, \cdots, x_k)
\]
is $L_W^{-1}$--cocoercive.
\end{example}

\begin{example}\normalfont
We consider a simple class of functionals depending on the barycenter of the measure.
Such functionals are closely related to the so--called cylindrical functionals on Wasserstein spaces, see e.g. \cite{FSS_2023, S_2023}.
For $\mu \in \mathcal P_2(\R^d)$, set 
\[
    m_\mu := \int_{\R^d} x \, d\mu(x).
\]
Consider $\Phi(\mu) := \varphi(m_\mu)$, where $\varphi : \R^d \to \R$ is a convex $L$--smooth function.
A direct computation gives 
\[
    \dfrac{\delta \Phi[\mu]}{\delta \mu}(x) = \langle \nabla \varphi(m_\mu), x \rangle.
\]
Hence, it follows from the Baillon--Haddad theorem that the operator $A_\mu := \nabla \dfrac{\delta \Phi[\mu]}{\delta \mu} = \nabla \varphi(m_\mu)$ is cocoercive.
\end{example}

\begin{example}\normalfont
We now consider a functional arising in min--max problems and competitive games over spaces of probability measures, see e.g. \cite{CHMR_2024}.
Let $d_1, d_2 \geq 1$ and let $U: \R^{d_1} \to \R$, $V: \R^{d_2} \to \R$ and $K : \R^{d_1} \times \R^{d_2} \to \R$ be of class $C^2$.
Under the assumptions of Proposition~\ref{prop.saddle}, for any $\mu \in \mathcal P_2(\R^{d_1})$ and $\nu \in \mathcal P_2(\R^{d_2})$, consider the zero--sum payoff
\begin{equation}
    \mathcal L(\mu, \nu) = \int_{\R^{d_1} \times \R^{d_2}} K(x, y) \, d\mu(x) d\nu(y) + \int_{\R^{d_1}} U(x) \, d\mu(x) - \int_{\R^{d_2}} V(y) \, d\nu(y).
\end{equation}
To formulate the associated saddle system, we consider the operator $A = (A_\rho)$ defined by 
\begin{equation}\label{A_saddle}
    A_\rho(x, y) :=
    \begin{pmatrix}
        \nabla U(x) + \displaystyle\int_{\R^{d_2}} \nabla_x K(x, v) \, d (\mathrm{proj}_2)_\# \rho \\
        \\
        \nabla V(y) - \displaystyle\int_{\R^{d_1}} \nabla_y K(u, y) \, d(\mathrm{proj}_1)_\#\rho
    \end{pmatrix},
\end{equation}
for every $\rho \in \mathcal P_2(\R^{d_1} \times \R^{d_2})$.
Under natural assumptions, it can be shown that the operator $A_\rho$ is cocoercive, see Proposition~\ref{prop.saddle}.
\end{example}

\begin{proposition}\label{prop.saddle}
Let $d_1, d_2 \geq 1$ and let $U: \R^{d_1} \to \R$, $V: \R^{d_2} \to \R$ and $K : \R^{d_1} \times \R^{d_2} \to \R$ be of class $C^2$.
Assume that the Hessians of $U, V$ and $K$ are uniformly bounded and there exists $\kappa > 0$ such that $\nabla^2 U \succeq \kappa \id$ and $\nabla^2 V \succeq \kappa \id$.
Assume in addition that 
\[
    \nabla^2_{xx} K(x, y) \succeq 0 \quad \text{ and } \quad \nabla^2_{yy} K(x, y) \preceq 0 \quad \text{ for every } (x, y) \in \R^{d_1} \times \R^{d_2}.
\]
Then, there exists $\beta > 0$ such that the operator $A_\rho$ defined by~\eqref{A_saddle} is $\beta$--cocoercive.
\end{proposition}

\begin{proof}
Set $\hat d := d_1 + d_2$.
Fix  $\rho, \rho' \in \mathcal{P}_2(\R^{\hat d})$ and $\gamma \in \Gamma(\rho, \rho')$.
Notice first that the interaction term is monotone since
\[
    (x, y) \longmapsto (\nabla_x K(x, y), - \nabla_y K(x, y))
\]
is monotone.
Consequently, combining with the strong convexity of $U$ and $V$, we have
\[
    \int_{\R^{2\hat d}} \langle A_\rho(z) - A_{\rho'}(z'), z - z' \rangle \, d\gamma(z, z') \geq \kappa \int_{\R^{2\hat d}} |z - z'|^2 \, d\gamma(z, z').
\]
On the other hand, the uniform boundedness of the Hessians of $U, V$ and $K$ together with Jensen's inequality implies that, for some $L > 0$
\[
    \int_{\R^{2\hat d}} |A_\rho(z) - A_{\rho'}(z')|^2 \, d\gamma(z, z') \leq L^2 \int_{\R^{ 2 \hat d}} |z - z'|^2 \, d\gamma(z, z').
\]
Combining the above estimates, we arrive at
\[
    \int_{\R^{2\hat d}} \langle A_\rho(z) - A_{\rho'}(z'), z - z' \rangle \, d\gamma(z, z') \geq \dfrac{\kappa}{L^2} \int_{\R^{2\hat d}} |A_\rho(z) - A_{\rho'}(z')|^2 \, d\gamma(z, z').
\]
Therefore, $A_\rho$ is $(\kappa/L^2)$--cocoercive, completing the proof.
\end{proof}

\subsection{Examples for nonexpansive operators}\label{sec:example2}

The previous examples of cocoercive operators yield directly a large class of nonexpansive operators.
Proposition~\ref{prop:wass_nonexpansive} implies that $B := \id-\tau A$ is nonexpansive for $\tau \in (0, 2\beta]$.

\begin{example}\normalfont
Under the assumptions of Proposition~\ref{prop.PoIn.E}, for every $\tau\in(0,2/(L_V+L_W)]$, the operator
\[
    B_\mu(x):=x-\tau\left(\nabla V(x)+\int_{\R^d}\nabla W(x-y)\,d\mu(y)\right)
\]
is nonexpansive.
\end{example}

\begin{example}\normalfont
Consider the higher--order interaction energy from Example~\ref{Exm2}.
For every $\tau\in(0,2/L_W]$, the operator
\[
    B_\mu(x):=x-\tau\int_{\R^{d(k-1)}}\nabla_{x_1}W(x,x_2,\cdots,x_k)\,d\mu^{\otimes(k-1)}(x_2,\cdots,x_k)
\]
is nonexpansive.
\end{example}

\begin{example}\normalfont
Let $\varphi:\R^d\to\R$ be a convex $L$--smooth function and set $A_\mu:=\nabla\varphi(m_\mu)$.
As shown above, $A$ is $L^{-1}$--cocoercive.
Hence, for every $\tau\in(0,2/L]$, the operator
\[
    B_\mu(x):=x-\tau\nabla\varphi(m_\mu)
\]
is nonexpansive.
\end{example}

\begin{example}\normalfont
Under the assumptions of Proposition~\ref{prop.saddle}, the operator $A$ defined in~\eqref{A_saddle} is $(\kappa/L^2)$--cocoercive for some $L>0$.
Consequently, for every $\tau\in(0,2\kappa/L^2]$, the operator $B:=\id-\tau A$ is nonexpansive.
\end{example}
\medskip\noindent
\textbf{On the use of artificial intelligence.}
LLMs (Opus 5.5 and GPT-6) were used to assist the authors in editing the manuscript, in particular in preparing the figures. The main ideas of the proofs build on previous work of the authors; though some questions were posed to these models in order to check the correctness of certain arguments and the sharpness of some of the assumptions. The authors take full responsibility for the contents of the manuscript.

\nocite{*}
\bibliography{references}{}
\bibliographystyle{siam}

\end{document}